\documentclass[10pt,journal,compsoc,twocolumn]{IEEEtran}
\usepackage{times}
\usepackage{epsfig}
\usepackage{graphicx}
\usepackage{amsmath,mathrsfs}
\usepackage{amssymb}
\usepackage{enumerate}
\usepackage{threeparttable}
\usepackage{multirow}
\usepackage{booktabs}
\usepackage{algorithm}
\usepackage{algorithmic}
\usepackage{amsmath,bm,amsfonts}
\usepackage{mcode}
\usepackage{subcaption}
\usepackage{amsthm}
\usepackage{url}
\usepackage{caption,color}
\usepackage{pifont}

\newtheorem{theorem}{Theorem}
\newtheorem{proposition}{Proposition}

\newtheorem{lemma}{Lemma}
\newtheorem{definition}{Definition}
\renewcommand{\u}{\mathbf{u}}

\renewcommand{\v}{\mathbf{v}}

\newcommand{\alli}{i=1,\cdots,n}
\newcommand{\sumi}{\sum_{i=1}^n}

\newcommand{\x}{\mathbf{x}}

\newcommand{\G}{\mathbf{G}}
\newcommand{\Gkk}{\mathbf{G}^{k+1}}
\newcommand{\I}{\mathbf{I}}
\newcommand{\xk}{\mathbf{x}^{k}}

\newcommand{\betak}{\beta^{(k)}}
\newcommand{\betakk}{\beta^{(k+1)}}
\newcommand{\gammak}{\gamma^{(k)}}

\newcommand{\xkk}{\mathbf{x}^{k+1}}

\newcommand{\y}{\mathbf{y}}

\newcommand{\HH}{\mathbf{H}}
\newcommand{\Hk}{\mathbf{H}^{k}}

\newcommand{\Hkk}{\mathbf{H}^{k+1}}
\newcommand{\e}{\mathbf{e}}
\newcommand{\hatr}{\hat{r}}
\newcommand{\hatf}{\hat{f}}

\newcommand{\bL}{\mathbf{L}}

\newcommand{\C}{\mathbf{C}}

\newcommand{\Q}{\mathbf{Q}}
\newcommand{\bP}{\mathbf{P}}

\newcommand{\K}{\mathbf{K}}
\newcommand{\U}{\mathbf{U}}
\newcommand{\T}{\mathbf{T}}

\newcommand{\blambda}{\bm{\lambda}}
\renewcommand{\b}{\mathbf{b}}
\newcommand{\bzero}{\mathbf{0}}
\newcommand{\Axk}{\A\xk}

\newcommand{\Axkk}{\A\xkk}
\newcommand{\xbone}{\x_{B_1}}
\newcommand{\xbonek}{\x^k_{B_1}}
\newcommand{\xbonekk}{\x^{k+1}_{B_1}}
\newcommand{\xbtwo}{\x_{B_2}}
\newcommand{\sumione}{\sum_{i\in B_1}}
\newcommand{\sumitwo}{\sum_{i\in B_2}}
\newcommand{\xbtwok}{\x^k_{B_2}}
\newcommand{\xbtwokk}{\x^{k+1}_{B_2}}
\newcommand{\Abone}{\A_{B_1}}
\newcommand{\Abtwo}{\A_{B_2}}

\newcommand{\lambdak}{\bm{\lambda}^{k}}
\newcommand{\lambdakk}{\bm{\lambda}^{k+1}}

\newcommand{\hatlambdakk}{\hat{\bm{\lambda}}^{k+1}}

\newcommand{\A}{\mathbf{A}}
\newcommand{\J}{\mathbf{J}}
\newcommand{\B}{\mathbf{B}}
\newcommand{\X}{\mathbf{X}}

\newcommand{\E}{\mathbf{E}}
\newcommand{\Z}{\mathbf{Z}}

\newcommand{\lbar}{\left\|}
\newcommand{\rbar}{\right\|}

\newcommand{\bS}{\mathbf{S}}
\newcommand{\st}{\text{s.t.}}

\newcommand{\norm}[1]{\lVert#1\rVert}

\newcommand{\Pomega}{\mathcal{P}_{\Omega}}

\newcommand{\Kk}{\mathbf{K}^k}
\newcommand{\Gk}{\mathbf{G}^k}

\DeclareMathOperator*{\argmin}{argmin}

\font\fontforurl=cmtt10

\newcommand{\tX}{\bm{\mathcal{X}}}
\newcommand{\tL}{\bm{\mathcal{L}}}
\newcommand{\tS}{\bm{\mathcal{S}}}
\newcommand{\tE}{\bm{\mathcal{E}}}
\newcommand{\tM}{\bm{\mathcal{M}}}
\newcommand{\M}{\bm{\mathcal{M}}}

\begin{document}
	
	\title{A Unified Alternating Direction Method of Multipliers by Majorization Minimization}

	\author{Canyi~Lu,~\IEEEmembership{Student~Member,~IEEE,}
		Jiashi~Feng,~
		Shuicheng Yan,~\IEEEmembership{Senior~Member,~IEEE,}
		and~Zhouchen Lin,~\IEEEmembership{Senior~Member,~IEEE}
		
		\IEEEcompsocitemizethanks{\IEEEcompsocthanksitem C. Lu, J. Feng and S. Yan are with the Department of Electrical and Computer Engineering, National University of Singapore, Singapore (e-mail: canyilu@gmail.com; elefjia@nus.edu.sg; eleyans@nus.edu.sg).\protect\\
			\IEEEcompsocthanksitem Z. Lin is with the Key Laboratory of Machine Perception (MOE), School of EECS, Peking University, China (e-mail: zlin@pku.edu.cn).}
	}
	
	\markboth{IEEE TRANSACTIONS ON PATTERN ANALYSIS AND MACHINE INTELLIGENCE}%
	{Shell \MakeLowercase{\textit{et al.}}: Bare Advanced Demo of IEEEtran.cls for Journals}

	
	\IEEEtitleabstractindextext{%
		\begin{abstract}
Accompanied with the rising popularity of compressed sensing, the Alternating Direction Method of Multipliers (ADMM) has become the most widely used solver for   linearly constrained convex problems with separable objectives. In this work, we observe that many previous variants of ADMM  update the primal variable by minimizing different majorant functions with their convergence proofs given case by case. Inspired by the principle of majorization minimization, we respectively present the unified frameworks and convergence analysis for the Gauss-Seidel ADMMs and Jacobian ADMMs, which use different historical information for the current updating. Our frameworks further generalize previous ADMMs to the ones capable of solving  the problems with non-separable objectives by minimizing their separable majorant surrogates. We also show that the bound which measures the convergence speed of ADMMs depends on the tightness of the used majorant function. Then several techniques are introduced  to improve the efficiency of ADMMs by tightening the majorant functions. In particular, we propose  the Mixed Gauss-Seidel and Jacobian ADMM (M-ADMM) which alleviates the slow convergence issue of Jacobian ADMMs by absorbing   merits of the Gauss-Seidel ADMMs. M-ADMM can be further improved by using backtracking, wise variable partition and fully exploiting the structure  of the constraint.  Beyond the guarantee in theory, numerical experiments on both synthesized and real-world data further demonstrate the superiority of our new ADMMs in practice. Finally, we release a toolbox at \textcolor{red}{{\fontforurl https://github.com/canyilu/LibADMM}} that implements efficient ADMMs for many   problems in compressed sensing.
		\end{abstract}
		
		}
		\maketitle		
		\IEEEdisplaynontitleabstractindextext
		\IEEEpeerreviewmaketitle
		
		\ifCLASSOPTIONcompsoc
		\IEEEraisesectionheading{\section{Introduction}\label{sec:introduction}}
		\else
		
		\section{Introduction}
		\label{sec:introduction}
		\fi
		\IEEEPARstart{T}{his}  work aims to solve the following convex problem
		\begin{equation}\label{multipro}
		\min_{\x} f(\x)=f(\x_1,\cdots,\x_n), \ \text{s.t.} \ \A\x=\sum_{i=1}^{n}\A_i\x_i=\b,
		\end{equation}
		where $f: \mathbb{R}^{p_1\times\cdots\times p_n}\rightarrow\mathbb{R}$ is convex and $n \ (\geq2)$ denotes  the block number of variables. We denote $\x=[\x_1;\cdots;\x_n]$ with $\x_i\in\mathbb{R}^{p_i}$, and $\A=[\A_1,\cdots,\A_n]$ with $\A_i\in\mathbb{R}^{d\times p_i}$.
		Problem (\ref{multipro}) has drawn increasing attention recently for the emerging applications of compressive sensing in computer vision and signal processing, e.g., sparsity based face recognition \cite{SRC,elhamifar2011robust}, saliency detection \cite{shen2012unified}, motion segmentation \cite{liu2011latent,lu2013correlation,6482137}, image denoising \cite{goldstein2009split,LADMPS}, video denoising \cite{ji2010robust}, texture repairing \cite{liang2012repairing} and many others \cite{cheng2011multi,zhuang2012non,vinayak2014graph,xia2014robust,HMWG15} .

		To solve (\ref{multipro}), the popular Augmented Lagrangian Method (ALM) \cite{hestenes1969multiplier} updates the primal variable $\x$ by
		\begin{align}
		\xkk=&\arg\min_{\x} \mathcal{L}(\x,\lambdak,\betak)
		=\arg\min_{\x} f(\x)+r^k(\x),\label{updateallx}
		\end{align}
		where $\mathcal{L}$ is the augmented Lagrangian function defined as 
		$$\mathcal{L}(\x,\blambda,\beta)=f(\x)+\langle\blambda,\A\x-\b\rangle+\frac{\beta}{2}\|\A\x-\b\|^2,$$
		and
		\begin{equation}\label{almr}
		r^k(\x)=\frac{\betak}{2}\left\|\A\x-\b+\frac{\lambdak}{\betak}\right\|^2.
		\end{equation}
		Then the dual variable $\blambda$ is updated to minimize $-\mathcal{L}$ by gradient descent with the step size $\betak$, i.e.,
		\begin{equation}\label{updatelambdaall}
		\lambdakk=\lambdak+\betak(\A\xkk-\b).
		\end{equation}
		However, (\ref{updateallx}) may not be easily solvable, since $r^k$ is non-separable. The Alternating Direction Method of Multipliers (ADMM) \cite{gabay1976dual} instead solves (\ref{updateallx}) inexactly by updating $\x_i$'s in an alternating way and thus the per-iteration cost can be much lower. Many variants of ADMM have been proposed by using different properties of $f$ and $\A$. We will review the most related works in Section \ref{secreview}, and claim our contributions in Section \ref{sec_contribution}.
		
		\textbf{Notations.} The $\ell_2$-norm of a vector and Frobenius norm of a matrix are denoted as $\|\cdot\|$. The spectral norm and the smallest singular value of a matrix $\A$ are denoted as $\|\A\|_2$ and $\sigma_{\min}(\A)$, respectively. The identity matrix is denoted as $\I$ without specifying its size. The all-one vector is denoted as $\bm{1}$. We denote $\mathbb{S}$ and $\mathbb{S}_+$ as the set of symmetry and positive semidefinite matrices respectively and  define $\langle \bm{a},\bm{a}\rangle_{\A}=\|\bm{a}\|_{\A}^2=\bm{a}^\top\A\bm{a}$ for  $\A\in\mathbb{S}$. If $\A-\B$ is positive   semi-definite, then we denote $\A\succeq\B$. The block diagonal matrix $\text{Diag}\{\A_i,\alli\}$ has $\A_i$ as its $i$-th block on the diagonal. A function $f: \mathbb{R}^p\rightarrow\mathbb{R}$ is said to be $L$-smooth (or $\nabla f$ is Lipschitz continuous), if
		\begin{equation}\label{lipschitzone}
		\|\nabla f(\x)-\nabla f(\y)\|\leq L\|\x-\y\|, \ \forall\x,\y\in\mathbb{R}^p.
		\end{equation}

		\subsection{Review of ADMMs}\label{secreview}
		Most of  ADMMs are only able to solve (\ref{multipro}) with separable $f$; i.e., there exist $f_i$'s such that $f(\x)=\sumi f_i(\x_i)$. They can be categorized into Gauss-Seidel ADMMs and Jacobian ADMMs. The Gauss-Seidel ADMMs update $\x_i$'s in a sequential way, i.e., update $\xkk_i$ by fixing others as their latest versions, while the Jacobian ADMMs update $\x_i$'s in a parallel way, i.e., update each $\xkk_i$ by fixing $\x_j=\xk_j$, for all $j\neq i$. We review these two types of ADMMs respectively. The difference between ADMMs lies in the updating of $\x_i$'s, while $\blambda$ is updated in the same way by (\ref{updatelambdaall}).
		
		Gauss-Seidel ADMMs solve (\ref{multipro}) with $n=2$ blocks. The standard ADMM \cite{boyd2011distributed} solves (\ref{updateallx}) inexactly by updating $\x_1$ and $\x_2$ in a sequential way, i.e.,
		\begin{align}
		\xkk_1=&\arg\min_{\x_1}\mathcal{L}([\x_1;\xk_2],\lambdak,\betak)\notag\\
		=&\arg\min_{\x_1}f_1(\x_1)+r_1^k(\x_1),\label{admmx1}\\
		\xkk_2=&\arg\min_{\x_2}\mathcal{L}([\xkk_1;\x_2],\lambdak,\betak)\notag \\
		=&\arg\min_{\x_2}f_2(\x_2)+r^k_2(\x_2),\label{admmx2}
		\end{align}
		where
		\begin{align}
		r_1^k(\x_1)=&\frac{\betak}{2}\left\|\A_1\x_1+\A_2\xk_2-\b+\frac{\lambdak}{\betak}\right\|^2, \label{admr1}\\
		r^k_2(\x_2)=&\frac{\betak}{2}\left\|\A_1\xkk_1+\A_2\x_2-\b+\frac{\lambdak}{\betak}\right\|^2. \label{admr2}
		\end{align}

		By using different properties of $f_1$ and $\A_1$, $\x_1$ (the same discussion is also applicable to $\x_2$) can be updated more efficiently than solving (\ref{admmx1}). If $f_1$ is $L_1$-smooth, then $\x_1$ can be updated by
		\begin{equation}\label{padmx1}
		\xkk_1=\argmin_{\x_1}\hat{f}_1(\x_1)+r_1^k(\x_1),
		\end{equation}
		where $\hatf_1(\x_1)=f(\xk_1)+\langle\nabla f_1(\xk_1),\x_1-\xk_1\rangle+\frac{L_1}{2}\|\x_1-\xk_1\|^2$. The motivation is that $\hatf_1$ is a majorant (upper bound) function of $f_1$, i.e., $\hatf_1\geq f_1$ \cite{beck2009fast}. If $f_1=g_1+h_1$, where $g_1$ is convex and $h_1$ is convex and $L_1$-smooth, then $\x_1$ can be updated by (\ref{padmx1}) with $\hatf_1(\x_1)=g(\x_1)+h(\xk_1)+\langle\nabla h_1(\xk_1),\x_1-\xk_1\rangle+\frac{L_1}{2}\|\x_1-\xk_1\|^2$. In this case, $\hatf_1\geq f_1$. We name the method using (\ref{padmx1}) as Proximal ADMM (P-ADMM) for these two cases. Similar techniques have been used in  \cite{beck2009fast,ouyang2013stochastic}.
		
		If the columns of $\A_1$ are not orthogonal, solving (\ref{admmx1}) is usually very expensive especially when $f_1$ is nonsmooth. Then Linearized ADMM (L-ADMM) \cite{LADMAP} instead updates $\x_1$ by
		\begin{align}
		\xkk_1=\arg\min_{\x_1}f_1(\x_1)+\hat{r}_1^k(\x_1),\label{ladmmx1}
		\end{align}
		where $\hat{r}_1^k(\x_1)=r_1^k(\xk_1)+\langle\nabla r_1^k(\xk_1),\x_1-\xk_1\rangle+\frac{\eta_1}{2}\|\x_1-\xk_1\|^2$ with $\eta_1>\|\A_1\|_2^2$.
		Note that  $\hat{r}_1^k\geq r_1^k$ since $r^k_1$ is $\|\A_1\|_2^2$-smooth.
		For some nonsmooth $f_1$, e.g., the $\ell_1$-norm, (\ref{ladmmx1}) can be solved efficiently with a closed form solution.
		
		If $f_1$ is a sum of a nonsmooth function and an $L_1$-smooth function, then we can simultaneously use the majorant function $\hatf_1$ of $f_1$ as P-ADMM and $\hat{r}_1^k$ of $r_1^k$ as L-ADMM. Thus $\hatf_1+\hat{r}_1^k\geq f_1+r_1^k$. This motivates the Proximal Linearized ADMM (PL-ADMM) which updates $\x_1$ by
		\begin{align}
		\xkk_1=\arg\min_{\x_1}\hatf_1(\x_1)+\hat{r}_1^k(\x_1).\label{pladmmx1}
		\end{align}

		For (\ref{multipro}) with $n>2$ blocks of variables, the naive extension of Gauss-Seidel ADMMs may diverge \cite{chen2013direct}. To address this issue, several Jacobian ADMMs have been proposed by using different properties of $f_i$ and $\A_i$. The Linearized ADMM with Parallel Splitting (L-ADMM-PS) \cite{liu2013linearized} solves (\ref{updateallx}) inexactly by linearizing $r^k$ in (\ref{almr}) at $\xk_i$'s and updates $\x_i$'s in parallel by
		\begin{align}\label{ladmmpsxi}
		\xkk_i=\arg\min_{\x_i}& \ f_i(\x_i)+\left\langle\A_i^\top(\betak(\A\xk-\b)+\blambda^k),\x_i\right\rangle\notag\\
		&+\frac{\betak\eta_i}{2}\|\x_i-\xk_i\|^2,
		\end{align}
		where $\eta_i>n\|\A_i\|_2^2$. A more general method proposed in the Algorithm 4 of \cite{deng20131} updates $\x_i$'s in parallel by
		\begin{align}\label{pjadmmxi}
		\xkk_i=\arg\min_{\x_i}& f_i(\x_i)+\frac{\betak}{2}\left\|\A_i\x_i+\sum_{j\neq i}\A_j\xk_j-\b-\frac{\lambdak}{\betak}\right\|^2\notag\\
		&+\frac{\betak}{2}\|\x_i-\xk_i\|^2_{\G_i},
		\end{align}
		where $\G_i\succ (n-1)\A_i^\top\A_i$. Actually (\ref{ladmmpsxi}) is a special case of (\ref{pjadmmxi}) when $\G_i=\eta_i\I-\A_i^\top\A_i$ with $\eta_i>n\|\A_i\|_2^2$. So we name the method using (\ref{pjadmmxi}) as Generalized Linearized ADMM with Parallel Splitting (GL-ADMM-PS) in this work. If $f_i=g_i+h_i$, where $g_i$ is convex and $h_i$ is convex and $L_i$-smooth, then
		the Proximal Linearized ADMM with Parallel Splitting (PL-ADMM-PS) \cite{LADMPS} updates $\x_i$'s in parallel by
		\begin{align}\label{pladmmpsxi}
		\xkk_i=&\arg\min_{\x_i} \hatf_i(\x_i)+\left\langle\A_i^\top(\betak(\A\xk-\b)+\lambdak),\x_i\right\rangle\notag\\
		&+\frac{\betak\eta_i}{2}\|\x_i-\xk_i\|^2,
		\end{align}
		where $\hatf_i(\x_i)=g(\x_i)+h(\xk_i)+\langle\nabla h_i(\xk_i),\x_i-\xk_i\rangle+\frac{L_i}{2}\|\x_i-\xk_i\|^2$ and $\eta_i>n\|\A_i\|_2^2$. As we will show later, the updating rules (\ref{pjadmmxi}) and (\ref{pladmmpsxi}) are equivalent to minimizing different majorant functions of $f(\x)+r^k(\x)$ in (\ref{updateallx}).
		
		For the convergence guarantee, all the above ADMMs own the convergence rate $O(1/K)$ \cite{he20121,liu2013linearized,LADMPS}, where $K$ is the number of iterations. There are also some other works which consider different special cases of our problem (\ref{multipro}) and give different convergence rates of ADMMs. For example, the works \cite{goldstein2014fast,lu2016fast} propose fast ADMMs with better convergence rate. But their considered problems are quite specific  and their convergence guarantees require several additional assumptions. For problem (\ref{multipro}) with separable objective and $n>2$,  the works \cite{hong2012linear,lin2014sublinear,lin2015iteration} prove the convergence of the naive multi-blocks extension of ADMM  under various assumptions, 
		 e.g., full column rank of $\A_i$, strong convexity or Lipschitz continuity of some $f_i$ and some others which may be hard to be verified in practice.
		 The work \cite{wang2013solving} reformulates the multi-blocks problem into a two-block one by variable splitting and solves it by ADMM. But it is verified to be slower than  GL-ADMM-PS in \cite{deng20131} since the variable splitting substantially increases the number of variables and constraints,  especially when $n$ is large.


		\subsection{Contributions}\label{sec_contribution}
		From the above discussions, we observe that different ADMMs  can be regarded as variants of inexact ALM  in the sense that the primal variable $\xkk$ in ADMMs is updated by solving (\ref{updateallx}) in ALM approximately. This actually slows the convergence, but the per-iteration cost is lower. So there is a trade-off between the exactness of the subproblem optimization and the convergence speed. In practice, we balance both to choose the proper solver. 
		Generally, if $f$ is not very simple, e.g., sum of several nonsmooth functions, ADMMs are much more efficient than ALM. ADMMs use two main techniques for approximation and update   $\xkk$ in an easier way than ALM: Alternating Minimization (AM) and Majorization Minimization (MM)   \cite{langeoptimization}. AM, which updates one block each time when fixing others, makes the subproblems easier to solve.   For example, the updating of $[\xkk_1;\xkk_2]$ in  ADMM  (\ref{admmx1})-(\ref{admmx2}) is easier than the one  in ALM (\ref{updateallx}). But the cost of the one block updating may be still high and it can be further reduced by using MM, which minimizes a majorant function instead  of the original objective to find an approximated solution.  For example, as reviewed in Section \ref{secreview}, different Gauss-Seidel ADMMs update $\x_1$ by minimizing different majorant functions of the objective in standard ADMM (\ref{admmx1}), while different Jacobian ADMMs update $\x_i$'s by minimizing different majorant functions of the objective in ALM (\ref{updateallx}). Actually, Gauss-Seidel ADMMs first use AM and then apply  MM to update each block, while Jacobian ADMMs first use  MM  and then   AM to update each block (though this is equivalent to updating all blocks simultaneously). 
		Besides the primal variables, the    dual variable $\lambdakk$ updating in (\ref{updatelambdaall}) is also equivalent to  minimizing a majorant function of $-\mathcal{L}(\xkk,\blambda,\betak)$, i.e.,
		\begin{equation}\label{eqlambdamajor}
		\blambda^{k+1} = \argmin_{\blambda} -\mathcal{L}(\xkk,\blambda,\betak) + \frac{1}{2\betak}\norm{\blambda-\blambda^k}^2.
		\end{equation} 		
		These observations suggest that  MM provides a new insight to interpret ADMMs. The convergences of ADMMs which use different majorant functions are guaranteed, but they are proved case by case. It is not clear what is the role of MM in ADMMs. Another issue is that, in practice, one   can develop many ADMMs for the same problem. But it is generally difficult to see which one converges faster.  The proved same rate $O(1/K)$  in the worst case fails to characterize  the  different speeds of ADMMs in practice. We lack practical principles and guidelines for designing efficient ADMMs.

%
%

		In this work, we raise several crucial questions:
		\begin{itemize}
			\item What kind of majorant functions can be used in ADMMs?  
			\item Is that possible to give a unified convergence analysis of  ADMMs which use different majorant functions by using certain common properties of majorant functions?
			\item What is the connection between the convergence speed of ADMMs and the used majorant functions?
			\item How to choose the proper majorant functions   for designing efficient ADMMs?  
		\end{itemize}	
		In this work, we show  many interesting findings about ADMMs through the lens of MM.   We aim to address the above questions and in particular we make the following contributions. First, for a multivariable function $f$,  we propose the majorant first-order surrogate  function $\hatf$, which requires three conditions to be satisfied:  majorization, proximity and separability. The first two guarantee that  $\hatf$ is  a reasonable approximation of $f$, while the last one  makes the minimizing of $\hatf$ easy. Note that  the objective  $f$ in (\ref{multipro}) can be non-separable since we only need to minimize $\hatf$. Second, we present the unified frameworks of Gauss-Seidel ADMMs and Jacobian ADMMs based on our majorant first-order surrogate and give the unified convergence guarantee. They not only draw connections with existing ADMMs, but also extend them to solve new problems with non-separable objective. Third, we show that  the bound which measures the convergence speed of ADMMs depends on the tightness of the used majorant function. 
		The tighter, the faster.
		This explains our previous intuitive observation that ADMMs converge faster when   (\ref{updateallx}) in ALM is solved more accurately. Fourth, we develop several useful techniques to tighten the majorant surrogates and thus improve the efficiency of ADMMs. Consider (\ref{multipro}) with $n>2$, we propose the Mixed Gauss-Seidel and Jacobian ADMM (M-ADMM) algorithm. It divides $n$ blocks of variables into two super blocks, and then updates them in a sequential way as Gauss-Seidel ADMMs, while the variables in each super block are updated in a parallel way as Jacobian ADMMs. 
		M-ADMM  takes the     structure of $\A$, e.g., $\frac{1}{2}\norm{\A\x-\b}^2$ that may be partially separable, into account to compute a tighter majorant surrogate, while previous Jacobian ADMMs fail to do so. In addition, we   show how to partition $n$ blocks of  variables into two super blocks wisely, which is crucial in the efficient implementation of ADMMs.  The last contribution is the developed toolbox which implements efficient ADMMs for many popular problems in compressed sensing. See  	
		\centerline{\textcolor{red}{{\fontforurl https://github.com/canyilu/LibADMM}}.}
	   Though there are already many toolboxes in compressed sensing, the solved problems are more or less limited due to the applicability of the used solvers, e.g., SPAMS \cite{mairal2011sparse} and SLEP \cite{liuslep} focus more on sparse models and non-constrained problems.
		We instead focus on the constrained  problem ({\ref{multipro}}), which is much more general. See  a list of problems  in our toolbox in the supplementary material.
		
		\section{Majorant First-Order Surrogate of a Multivariable Function}
		\label{sec2}
		In this section, we propose the majorant first-order surrogate of the multivariable functions which enjoy some ``good" properties.

		\begin{definition}\label{Lem_lips} \textit{\em{\textit{\textbf{(Lipschitz Continuity)}}}}
			Let $f: \mathbb{R}^{p_1}\times\cdots\times\mathbb{R}^{p_n}\rightarrow\mathbb{R}$ be differentiable. Then $\nabla f$ is called Lipschitz continuous if there exist $\bL_i\succeq\bm{0}, i=1,\cdots,n$, such that
			\begin{equation}\label{lipmuv}
			|f(\x)-f(\y)-\langle\nabla f(\y),\x-\y\rangle|\leq\frac{1}{2}\sum_{i=1}^{n}\|\x_i-\y_i\|_{\bL_i}^2,
			\end{equation}
			for any $\x=[\x_1;\cdots;\x_n]$ and $\y=[\y_1;\cdots;\y_n]$ with $\x_i, \y_i\in\mathbb{R}^{p_i}$. In this case, we say that $f$ is $\{\bL_i\}_{i=1}^n$-smooth.
		\end{definition}
		The Lipschitz continuity of the multivariable function is crucial in this work. It is different from the single variable case defined in (\ref{lipschitzone}).
		For $n=1$, (\ref{lipmuv}) holds if (\ref{lipschitzone}) holds (Lemma 1.2.3 in \cite{nesterov2004introductory}), but not vice versa. This motivates the above definition.
		
		\begin{definition}\label{Lem_strongcon} \textit{\em{\textit{\textbf{(Strong Convexity)}}}}
			A function $f: \mathbb{R}^{p_1}\times\cdots\times\mathbb{R}^{p_n}\rightarrow\mathbb{R}$ is called $\{\bP_i\}_{i=1}^n$-strongly convex if there exist $\bP_i\succeq\bzero$, $i=1,\cdots,n$, such that for any $\y_i\in\mathbb{R}^{p_i}$, the function $\x\rightarrow f(\x)-\frac{1}{2}\sum_{i=1}^{n}\|\x_i-\y_i\|_{\bP_i}^2$ is convex.
		\end{definition}

		\begin{definition}\label{defmfos} \textit{\em{\textit{\textbf{(Majorant First-Order Surrogate)}}}}
			A function $\hatf: \mathbb{R}^{p_1}\times\cdots\times\mathbb{R}^{p_n}\rightarrow\mathbb{R}$ is a majorant first-order surrogate of $f: \mathbb{R}^{p_1}\times\cdots\times\mathbb{R}^{p_n}\rightarrow\mathbb{R}$ near $\bm\kappa=[\bm\kappa_1;\cdots;\bm\kappa_n]$ with $\bm\kappa_i\in\mathbb{R}^{p_i}$ when the following conditions are satisfied:
			\begin{itemize}
				\item \textbf{Majorization}: $\hatf$ is a majorant function of $f$, i.e., $\hatf(\x)\geq f(\x)$ for any $\x$.  
				\item \textbf{Proximity}: there exists $\bL_i\succeq \bm{0}$ such that the approximation error $h(\x):=\hatf(\x)-f(\x)$ satisfies
				\begin{align}\label{keylm11}
				|h(\x)|\leq\frac{1}{2}\sum_{i=1}^{n}\|\x_i-\bm\kappa_i\|_{\bL_i}^2.
				\end{align}
				\item \textbf{Separability}: $\hatf$ is separable w.r.t. $\x_i$'s; i.e., there exist $\hatf_i$'s such that $\hatf(\x)=\sumi \hatf_i(\x_i)$.
			\end{itemize}	
			We denote by $\mathcal{S}_{\{\bL_i,\bP_i\}_{i=1}^n}(f,\bm{\kappa})$ the set of $\{\bP_i\}_{i=1}^n$-strongly convex surrogates.
		\end{definition}
		In MM, one aim to find an approximated solution to $\min_{\x} f(\x)$ by solving $\min_\x \hatf(\x)$, which is easier. To this end,  the above three conditions on $\hatf$ look reasonable. Majorization  guarantees that $f(\x)$ tends to be minimized when $\hatf(\x)$ is minimized. Proximity means that $\hatf(\x)$ cannot be too loose and this guarantees a controllable  approximation to $f(\x)$. The separability makes the optimization on $\hatf(\x)$ easier than $f(\x)$, which can be non-separable.  This is important for multi-blocks optimization.

		Note that $\bL_i$ measures the difference $\hatf-f$, or the tightness of the majorant surrogate $\hatf$. If $\norm{\bL_i}_2$ is smaller, then the majorant surrogate is tighter. This   plays an important role in this work.
		\begin{lemma}  \label{lemmasmmoth}
			If the approximation error $h(\x)=\hatf(\x)-f(\x)$ satisfies the following \textbf{Smoothness} assumption, i.e., 
			\begin{align}
			\text{$h(\x)$ is $\{\bL_i\}_{i=1}^n$-smooth, $h(\bm\kappa)=0$ and $\nabla h(\bm\kappa)=0$,} \label{eqsmoothenessh}
			\end{align}
			then the \textbf{Proximity} assumption in (\ref{keylm11}) holds.
		\end{lemma}
		Lemma \ref{lemmasmmoth} can be obtained by using (\ref{lipmuv}) for $h$ at $\bm\kappa$. Lemma \ref{lemmasmmoth} is useful to verify the Proximity assumption. 
		Some widely used majorant first-order surrogates are (see Lemma \ref{lemma5} in   Appendix):
		\begin{itemize}
			\item \textbf{Proximal Surrogates.}
			For any $f$ and $\bL\succeq\bzero$,  $\hatf\in\mathcal{S}_{\{\bL,\bL\}}(f,\bm{\kappa})$, where $\hatf(\x)=f(\x)+\frac{1}{2}\|\x-\bm{\kappa}\|_{\bL}^2$.
			\item \textbf{Lipschitz Gradient Surrogates.}
			Let   $f$ be  $\{\bL_i\}_{i=1}^n$-smooth. Then $\hatf\in\mathcal{S}_{\{\bL_i,\bL_i\}_{i=1}^n}(f,\bm{\kappa})$, where $\hat{f}(\x)=f(\bm\kappa)+\langle \nabla f(\bm\kappa),\x-\bm\kappa\rangle+\frac{1}{2}\sum_{i=1}^{n}\|\x_i-\bm{\kappa_i}\|_{\bL_i}^2$.
			\item \textbf{Proximal Gradient Surrogates.}
			Let   $f=f_1+f_2$, where $f_1$ is $\{\bL_i\}_{i=1}^n$-smooth.  Then  $\hatf\in\mathcal{S}_{\{\bL_i,\bL_i\}_{i=1}^n}(f,\bm{\kappa})$, where $\hat{f}(\x)=f_1(\bm\kappa)+\langle \nabla f_1(\bm\kappa),\x-\bm\kappa\rangle+\frac{1}{2}\sum_{i=1}^{n}\|\x_i-\bm{\kappa_i}\|_{\bL_i}^2+f_2(\x)$.
		\end{itemize}
		Note that if $f$ is separable, then $\hatf=f$ is also  a majorant first-order surrogate of $f$. Some other examples, e.g., DC programming surrogates, can be found in \cite{Mairal13optimizationwith}.

		\begin{lemma}\label{lemma44} \textit{\em{\textit{\textbf{(Key Property of the Majorant First-Order Surrogate)}}}}
			Let $\hat{f}\in\mathcal{S}_{\{\bL_i,\bP_i\}_{i=1}^n}(f,\bm{\kappa})$. Then, we have
			\begin{align}
			&f(\x)+\langle\mathbf{u},\y-\x\rangle-f(\y)\notag\\
			\leq&\frac{1}{2}\sum_{i=1}^{n}\left(\|\y_i-\bm\kappa_i\|^2_{\bL_i}-\|\y_i-\x_i\|^2_{\bP_i}\right), \ \forall \x,\y,\label{keylm22}
			\end{align}
			where $\mathbf{u}\in\partial \hat{f}(\x)$ is any subgradient of the convex $\hatf$.
		\end{lemma}
		
		The majorant first-order surrogate given in Definition \ref{defmfos} is motivated by \cite{Mairal13optimizationwith}. However, they have many key differences:
		\begin{itemize}
			\item Our majorant first-order surrogate is defined based on the multivariable function and thus it is much more general than the single variable case considered in \cite{Mairal13optimizationwith}. For example, the Lipschitz continuity of the multivariable function is different; the \textbf{Separability} of $\hatf$ is new.
			\item For approximation error $h=\hatf-f$, we use the \textbf{Proximity} assumption in (\ref{keylm11}) which is less restricted than of the \textbf{Smoothness} assumption in (\ref{eqsmoothenessh}). We only require the error $h$ to be bounded, and it is not necessary to be smooth.
			\item Our Lemma \ref{lemma44} is new and it plays a central role in our convergence analysis.  Lemma 2.1 in \cite{Mairal13optimizationwith} also introduces some   properties of the majorant first-order surrogate. But their bounds are too loose and   are not applicable to our proofs  due to the   constraint of (\ref{multipro}) considered in this work.
			\item The considered constrained problem in this work is   different from the non-constrained problem in \cite{Mairal13optimizationwith}. When proving Proposition 2.3 in \cite{Mairal13optimizationwith}, they use a key property $f(\xkk)\leq f(\xk)$, while this does not hold in ADMMs. 
		\end{itemize}

		At the end of this section, we discuss some properties of $\frac{1}{2}\norm{\A\x-\b}^2$ which are important  for designing efficient ADMMs.  


		\begin{lemma}\label{lem5}
			Let $r(\x)=\frac{1}{2}\|\A\x-\b\|^2$, where $\x=[\x_1;\cdots;\x_n]$, $\A=[\A_1,\cdots,\A_n]$ and $\b$ are of compatible sizes.  We have
			\begin{enumerate}[(1)]
				\item $r(\x)$ is $\{\bL'_i\}_{i=1}^n$-smooth. The choice of $\bL'_i$ depends on $\A_i^\top\A_i$.
				\item $r(\x)\leq \hatr(\x)$, where
				\begin{align}
				\hatr(\x)=&\frac{1}{2}\sumi\left\|\A_i\x_i+\sum_{j\neq i}\A_j\y_j-\b\right\|^2\notag \\
				&+\frac{1}{2}\sumi\|\x_i-\y_i\|^2_{\G_i}   +\frac{1-n}{2}\|\A\y-\b\|^2,\label{sqrsurr2}
				\end{align}
				for any  $\y=[\y_1;\cdots;\y_n]$ and  $\G_i\succeq \bL'_i-\A_i^\top\A_i$.
				\item If $\G_i=\eta_i\I-\A_i^\top\A_i$ with $\eta_i\geq \|\bL_i'\|_2$, (\ref{sqrsurr2}) reduces to
				\begin{align}\label{sqrsurr1}
				\hatr(\x)=&\sumi\left\langle\x_i-\y_i,\A_i^\top(\A\y-\b)\right\rangle \notag \\
				&+\sumi\frac{\eta_i}{2}\|\x_i-\y_i\|^2+\frac{1}{2}\|\A\y-\b\|^2.
				\end{align}
			\end{enumerate}			
			
		\end{lemma}
To guarantee that $\hatr\geq r$, it is required to choose $\bL'_i$ with $\norm{\bL'_i}_2$ sufficiently large. 
Without any additional assumption on $\A$, we can choose $\bL'_i=n\A_i^\top\A_i$. This explains the choice of $\eta_i>\norm{\bL'_i}_2=n\norm{\A_i}_2^2$   in L-ADMM-PS (\ref{ladmmpsxi}). However, such a choice of $\bL'_i$ may not be good since it does not make fully use of the structure of $\A$, and thus $\hatr$ may not be a tight surrogate of $r$. For example, let $\A_1 = [\C_1;\bm{0}]$, $\A_2 = [\C_2;\bm{0}]$, $\A_3 = [\bm{0}; \C_3]$,  $\A_4= [\bm{0}; \C_4]$, and $\b=[\b_1;\b_2]$ of compatible sizes. Then $r(\x)=\frac{1}{2}\norm{\sum_{i=1}^{2}\C_i\x_i-\b_1}^2+\frac{1}{2}\norm{\sum_{i=3}^{4}\C_i\x_i-\b_2}^2$. We can choose  $\bL'_i = 2\A_i^\top\A_i$, which is much better than   $4\A_i^\top\A_i$. Actually, the choice of $\bL'_i$ depends on the separability of $r$. In practice, it is easy to compute $\bL'_i$ when given $\A$. 
A good choice of $\bL'_i$ gives a tight surrogate $\hatr$, and  this may significantly improve the efficiency of Jacobian ADMMs (see Section \ref{sec4}).

%
%
%
%
%

		\section{Unified Gauss-Seidel ADMMs}\label{sec3}
		In this section, we consider   solving   (\ref{multipro}) with $n=2$ blocks  by a 
		unified framework of Gauss-Seidel ADMMs.
		  In the $(k+1)$-th iteration, we compute the  majorant surrogate $\hatf^k$ of $f$ near $\xk$, i.e., $\hat{f}^k\in\mathcal{S}_{\{\bL_i, \bP_i\}_{i=1}^2}(f,\xk)$ and $\hatf^k$ is separable, i.e., $\hatf^k(\x)=\hatf_1^k(\x_1)+\hatf_2^k(\x_2)$. For $r_1^k$ and $r_2^k$ in (\ref{admr1}) and (\ref{admr2}), we construct their proximal surrogates respectively as follows\footnote{Note that the definitions of $\hatr_i^k$ in Section \ref{sec3}, \ref{sec4} and \ref{sec5} are different.}
		\begin{align}
		\hat{r}_1^k(\x_1)=&r_1^k(\x_1)+\frac{\betak}{2}\|\x_1-\xk_1\|^2_{\G_1},\label{hatr1admms}\\
		\hat{r}_2^k(\x_2)=&r_2^k(\x_2)+\frac{\betak}{2}\|\x_2-\xk_2\|^2_{\G_2},\label{hatr2admms}
		\end{align}
		where $\G_1\succeq\bzero$ and $\G_2\succ\bzero$.
		Then we  update $\x_1$ and $\x_2$ by
		\begin{align}
		\xkk_1  = & \arg\min_{\x_1} \hat{f}_1^k(\x_1)+\hatr_1^k(\x_1),\label{updatex1}\\
		\xkk_2  = & \arg\min_{\x_2} \hat{f}_2^k(\x_2)+\hatr_2^k(\x_2)\label{updatey1}.
		\end{align}
		Finally, $\blambda$ is updated by (\ref{updatelambdaall}).
		This leads to the unified framework of Gauss-Seidel ADMMs, as shown in Algorithm \ref{alg1}. 
		
		Note that in Algorithm \ref{alg1}, $f$ is not necessarily separable. In this case, our algorithm and the convergence
		guarantee shown later are completely new.
		If $f$ is already separable, then the objectives in (\ref{updatex1}) and (\ref{updatey1}) are majorant surrogates of the ones in (\ref{admmx1}) and (\ref{admmx2}), respectively. Many previous Gauss-Seidel ADMMs are special cases by using different majorant surrogates $\hatf_1$ and $\hat{r}_1^k$ (depending on $\Gk_1$) in Algorithm \ref{alg1}. See Table \ref{tab1} for a summary.

		Assume that there exists an KKT point $(\x^*,\blambda^*)$ of (\ref{multipro}), i.e.,
		$\A\x^*=\b$ and $-\A^\top\blambda^*\in\partial f(\x^*)$. Previous works prove that ADMMs converge to the KKT point at the rate $O(1/K)$ ($K$ is the number of iterations) in different ways. The works \cite{he20121,ouyang2013stochastic} give the same rate of ADMM, L-ADMM, and P-ADMM. But they require that both the primal and dual feasible sets should be bounded. The work \cite{LADMPS} removes the above assumptions and shows that the convergence rates of L-ADMM-PS and PL-ADMM-PS are
		\begin{align}\label{lem3eq12}
		&f(\bar{\x}^{K})-f(\x^*)+\langle\A^\top\blambda^*,\bar{\x}^{K}-\x^*\rangle+\frac{\alpha}{2}\lbar\A\bar{\x}^{K}-\bm{b}\rbar^2 \notag \\
		\leq & O(1/K),
		\end{align}
		where $\bar{\x}^K$ is a weighted sum of $\x^k$'s and $\alpha>0$. Now we give the  convergence bound of Algorithm \ref{alg1} as (\ref{lem3eq12}).
		
		\begin{algorithm}[t]
			\caption{A Unified Framework of Gauss-Seidel ADMMs}
			\textbf{For} $k=0,1,2,\cdots$ \textbf{do}
			\begin{enumerate}
				\item Compute a majorant first-order surrogate $\hat{f}^k\in\mathcal{S}_{\{\bL_i, \bP_i\}_{i=1}^2}(f,\xk)$ with $\hatf^k(\x)=\hatf_1^k(\x_1)+\hatf_2^k(\x_2)$.		\item Update $\x_1$ by solving (\ref{updatex1}).
				\item Update $\x_2$ by solving (\ref{updatey1}).
				\item Update $\bm{\lambda}$ by $\blambda^{k+1}=\lambdak+\betak(\A\xkk-\b)$.
				\item Choose $\betakk\geq\betak$.
			\end{enumerate}
			\textbf{end}
			\label{alg1} 
		\end{algorithm}


		\begin{theorem}\label{them1}
			In Algorithm \ref{alg1}, assume that $\hatf^k\in\mathcal{S}_{\{\bL_i,\bP_i\}_{i=1}^2}(f,\xk)$ with $\bP_i\succeq \bL_i\succeq \bzero$,  $i=1,2$, $\G_1\succeq\bzero$ in (\ref{hatr1admms}), and $\G_2\succ\bzero$ in  (\ref{hatr2admms}).
			For any $K>0$, let $\bar{\x}^{K}=\sum_{k=0}^{K}\gammak\xkk$ with $\gammak={(\betak)^{-1}}/{\sum_{k=0}^{K}(\betak)^{-1}}$. Then
			\begin{align}
			&f(\bar{\x}^{K})-f(\x^*)+\langle \A^\top\blambda^*,\bar{\x}^{K}-\x^*\rangle+\frac{\beta^{(0)}\alpha}{2}\|\A\bar{\x}^{K}-\b\|^2\notag\\
			\leq&\frac{\sum_{i=1}^2\|\x_i^*-\x_i^0\|^2_{\HH^0_i}+\|\blambda^*-\blambda^0\|^2_{\HH^0_3}}{2\sum_{k=0}^{K}\left(\betak\right)^{-1}}, \label{thm1con}
			\end{align}
			where $\alpha=\min\left\{\frac{1}{2},\frac{\sigma^2_{\min}(\G_2)}    {2\|\A_2\|^2_2}\right\}$, $\HH^0_1=\frac{1}{\beta^{(0)}}\bL_1+\G_1$, $\HH^0_2=\frac{1}{\beta^{(0)}}\bL_2+\A_2^\top\A_2+\G_2$, and $\HH^0_3=\left({1/\beta^{(0)}}\right)^2\I$.
		\end{theorem}
		Consider $\HH^0_i$, $i=1,2$, at the RHS of (\ref{thm1con}), it can be seen that they depend on $\bL_i$ and $\G_i$, which control the difference $\hatf-f$ and $\hatr^k_i-r_i^k$, respectively.   This suggests a faster convergence when using tighter majorant surrogates, though the convergence rate of Gauss-Seidel ADMMs in Algorithm \ref{alg1} is $O(1/K)$ when $\betak$'s are bounded. 
		
		Note that the assumption $\G_2\succ\bzero$ guarantees that $\alpha>0$. Such an assumption is also used in \cite{he20121,ouyang2013stochastic} which prove the same convergence rate in different ways.
		It suggests that using $\G_2\succ\bzero$ instead of $\G_2=\bzero$ in the traditional ADMM can achieve the $O(1/K)$ convergence rate. 

		
		\renewcommand{\arraystretch}{1.3}
		\begin{table}[t]
			\center
			\caption{\small Previous Gauss-Seidel ADMMs are special cases of Algorithm \ref{alg1} with different $\hatf_1$ and $\G_1$. In this table, $\eta_1>\|\A_1\|_2^2$.} 
			\small
			\begin{tabular}{r|c|c}\hline
				\multicolumn{1}{c|}{}   & \multicolumn{1}{c|}{$\hatf_1^k(\x_1)$} & $\G_1$                                 \\\hline
				ADMM                     & $f_1(\x_1)$                      & $\bzero$                               \\\hline
				\multirow{2}{*}{P-ADMM}  & Lipschitz Gradient Surrogate or  & \multirow{2}{*}{$\bzero$}              \\
				& Proximal Gradient Surrogate     &                                        \\\hline
				L-ADMM                   & $f_1(\x_1)$                      & $\eta_1\I-\A_1^\top\A_1$                 \\\hline
				\multirow{2}{*}{PL-ADMM} & Lipschitz Gradient Surrogate or  & \multirow{2}{*}{$\eta_1\I-\A_1^\top\A_1$  } \\
				& Proximal Gradient Surrogate     &                                       \\\hline		
			\end{tabular}\vspace{-0.2cm}\label{tab1}
		\end{table}

		\section{Unified Jacobian ADMMs}\label{sec4}
		In this section, we consider   solving (\ref{multipro}) with $n>2$ by a unified framework of Jacobian ADMMs. 	
		The motivation is to solve (\ref{updateallx}) inexactly by minimizing a majorant surrogate of $f(\x)+r^k(\x)$. In the $(k+1)$-th iteration, we first compute the  majorant surrogate of $f$ near $\xk$, i.e., $\hat{f}^k\in\mathcal{S}_{\{\bL_i, \bP_i\}_{i=1}^n}(f,\xk)$, and $\hatf^k$ is separable, $\hatf^k(\x)=\sumi\hatf_i(\x_i)$. Assume that $\frac{1}{2}\norm{\A\x}^2$ is $\{\bL'_i\}_{i=1}^n$-smooth. For $r^k$ in (\ref{updateallx}), we define its majorant surrogate $\hat{r}^k$ by using (\ref{sqrsurr2}), i.e., $\hat{r}^k(\x)=\sumi\hat{r}_i^k(\x_i)$, where
		\begin{align}
		\frac{\hat{r}_i^k(\x_i)}{\betak}=&\frac{1}{2}\left\|\A_i\x_i+\sum_{j\neq i}\A_j\xk_j-\b+\frac{\blambda^k}{\betak}\right\|^2\notag \\&+\frac{1}{2}\|\x_i-\xk_i\|^2_{\G_i}+c_i^k,\label{jadmmri}
		\end{align}
		with $\G_i\succ \bL'_i-\A_i^\top\A_i$  and $c_i^k$'s are constants satisfying $\sumi c_i^k=\frac{1-n}{2}\|\A\x^k-\b\|^2$. Thus $\hatf^k(\x)+\hatr^k({\x})$ is a majorant surrogate of $f(\x)+ r^k(\x)$ in (\ref{updateallx}). Now we minimize $\hatf^k(\x)+\hatr^k({\x})$ instead to update $\x$, i.e.,
		\begin{equation}\label{updatexmul}
		\xkk=\arg\min_{\x}\hatf^k(\x)+\hatr^k({\x}).
		\end{equation}
		Note that both $\hat{f}$ and $\hat{r}^k$ are separable. Thus solving (\ref{updatexmul}) is equivalent to updating each $\x_i$ in parallel, i.e.,
		\begin{equation}\label{updatexmult}
		\xkk_i=\arg\min_{\x_i}\hatf_i^k(\x_i)+\hatr_i^k({\x_i}).
		\end{equation}
		Finally $\blambda$ is updated by (\ref{updatelambdaall}). This leads to the unified framework of Jacobian ADMMs, as shown in Algorithm \ref{alg2}.
		
		If $f$ is non-separable, then our algorithm and convergence guarantee shown later are completely new. If $f$ is separable, several previous Jacobian ADMMs are special cases by using different majorant surrogates $\hatf_i$ and $\hat{r}_i^k$ (depending on $\G_i$) in Algorithm \ref{alg2}. See Table \ref{tab2} for a summary.

		\begin{theorem}\label{them2}
			In Algorithm \ref{alg2}, assume that $\hatf^k\in\mathcal{S}_{\{\bL_i,\bP_i\}_{i=1}^n}(f,\xk)$ with $\bP_i\succeq \bL_i\succeq \bzero$, $\frac{1}{2}\norm{\A\x}^2$ is $\{\bL'_i\}_{i=1}^n$-smooth, and $\G_i\succ \bL_i'-\A_i^\top\A_i$ in (\ref{jadmmri}).
			For any $K>0$, let $\bar{\x}^K=\sum_{k=0}^{K}\gammak\xkk$ with $\gammak={(\betak)^{-1}}/{\sum_{k=0}^{K}(\betak)^{-1}}$. Then
			\begin{align}
			&f(\bar{\x}^K)-f(\x^*)+\langle \A^\top\blambda^*,\bar{\x}^K-\x^*\rangle+\frac{\beta^{(0)}\alpha}{2}\|\A\bar{\x}^K-\b\|^2\notag\\
			\leq&\frac{\sum_{i=1}^n\|\x_i^*-\x_i^0\|^2_{\HH^0_i}+\|\blambda^*-\blambda^0\|^2_{\HH^0_{n+1}}}{2\sum_{k=0}^{K}\left(\betak\right)^{-1}}, \label{thm2con}
			\end{align}
			where $\alpha=\min\left\{\frac{1}{2},\frac{\sigma^2_{\min}\left(\text{Diag}\{\A_i^\top\A_i+\G_i,i=1,\cdots,n\}-\A^\top\A\right)}{2\|\A\|^2_2}\right\}$, $\HH^0_i=\frac{1}{\beta^{(0)}}\bL_i+\A_i^\top\A_i+\G_i$, $\alli$, and $\HH^0_{n+1}=\left({1/\beta^{(0)}}\right)^2\I$.
		\end{theorem}
		The above bound implies an interesting connection between the convergence speed and the tightness of the   majorant surrogates.
		For simplicity, let $\betak=\beta$. Then   (\ref{thm2con}) reduces to
		\begin{align}
		&\frac{\sum_{i=1}^n\|\x_i^*-\x_i^0\|^2_{\beta\HH^0_i}+\frac{1}{\beta}\|\blambda^*-\blambda^0\|^2}{2(K+1)} \notag\\
		\leq & \frac{\frac{1}{2}\sum_{i=1}^n\|\x_i^*-\x_i^0\|^2_{\bL_i+\beta\bL'_i}+\frac{1}{2\beta}\|\blambda^*-\blambda^0\|^2}{(K+1)}, \label{eqjacobabound2}
		\end{align}
		where (\ref{eqjacobabound2}) uses $\G_i\succ \bL'_i-\A_i^\top\A_i$. Now consider the two constant terms in the numerator of (\ref{eqjacobabound2}). The first term controls the tightness of the used majorant surrogate for the $\x$ updating, i.e.,
		$|\hatf^0(\x^*)+\hatr^0(\x^*)-f(\x^*)-r(\x^*)| \leq \frac{1}{2}\sum_{i=1}^{n} \norm{\x_i^*-\x^0_i}_{\bL_i+\beta\bL'_i}^2$, which uses (\ref{keylm11}) with $\x=\x^*$ and $k=0$. The second term is actually the difference function $\frac{1}{2\beta}\norm{\blambda-\lambdak}^2$ between $-\mathcal{L}(\xkk,\blambda,\betak)$ and its   majorant surrogate in (\ref{eqlambdamajor}) when $\blambda=\blambda^*$ and $k=0$. So the convergence bound depends on the tightness of the used majorant surrogates for both the primal and dual variables updating.  If $\hatf^k + \hatr^k$ is tighter (associated to the $\x$ updating) or  $\beta$ is larger (associated to the $\blambda$ updating), the algorithm converges faster. In practice, ADMMs stop based on certain criteria induced by the KKT conditions. If  $\beta$ is relatively larger, the algorithm seems to converge faster but the objective function value may be larger. How to choose   the best $\beta $ or $\betak$ is still an open issue. In this work, we focus the discussion on how to improve the tightness of the majorant surrogate for the primal variable updating.

		Note that   Algorithm \ref{alg2} improves previous Jacobian ADMMs which use $\bL'_i=n\A_i^\top\A_i$. Such a choice of $\bL'_i$ does not fully use the structure of $\A$ or $r(\x)$ (see the discussions after Lemma \ref{lem5}). Our Algorithm \ref{alg2} instead uses  the $\{\bL'_i\}_{i=1}^n$-smooth property of $r(\x)$. This may make the surrogate $\hatr^k(\x)$  tighter and thus the algorithm converges faster. In Section \ref{sec5}, we  discuss how to further improve the tightness of $\hatr^k(\x)$ by introducing alternating minimization in Jacobian ADMMs  
		and the backtracking technique.
		
 \begin{algorithm}[t]
 	\caption{A Unified Framework of Jacobian ADMMs}
 	\textbf{For} $k=0,1,2,\cdots$ \textbf{do}
 	\begin{enumerate}
 		\item Compute a majorant first-order surrogate $\hat{f}^k\in\mathcal{S}_{\{\bL_i, \bP_i\}_{i=1}^n}(f,\xk)$ with $\hatf^k(\x)=\sumi \hatf^k_i(\x_i)$.
 		\item Update $\x_i$, $\alli$, in parallel by solving (\ref{updatexmult}).
 		\item Update $\bm{\lambda}$ by $\blambda^{k+1}=\lambdak+\betak(\A\x^{k+1}-\b)$.
 		\item Choose $\betakk\geq\betak$.
 	\end{enumerate} 
 	\textbf{end}
 	\label{alg2}%
 \end{algorithm}
		

%

		\section{Mixed Gauss-Seidel and Jacobian ADMM}
		\label{sec5}
		Consider   solving   (\ref{multipro}) with $n=2$ by Gauss-Seidel ADMMs   and 
		Jacobian ADMMs,   the former one will converge faster.    The reason is that Jacobian ADMMs require $\G_i\succ \A_i^\top\A_i$, while Gauss-Seidel ADMMs only require $\G_i\succ\bzero$. Thus the bound in (\ref{thm1con}) is expected to be tighter than the one in (\ref{thm2con}). The superiority of Gauss-Seidel ADMMs over Jacobian ADMMs is   that the former   first use  alternating minimization to reduce the complexity of the problem (fewer   variables) and then the used majorant surrogate can be tighter when using majorization minimization. 
		
		In this section, we consider problem (\ref{multipro}) with $n>2$ blocks. 
		We propose the Mixed Gauss-Seidel and Jacobian ADMM (M-ADMM), which introduces the alternating minimization before using majorization minimization. M-ADMM first divides these $n$ blocks $\x=[\x_1;\cdots;\x_n]$ into two super blocks, i.e., $\xbone=[\x_i, i\in B_1]$ with $n_1$ blocks of variables, and $\xbtwo=[\x_i, i\in B_2]$ with $n_2$ blocks of variables, where $B_1$ and $B_2$ satisfy $B_1\cap B_2=\varnothing$ and $B_1 \cup B_2=\{1,\cdots,n\}$. Then $\xbone$ and $\xbtwo$ are updated in a sequential way as Gauss-Seidel ADMMs,  while $\x_i$'s in each super block are updated in a parallel way as Jacobian ADMMs. As shown later, M-ADMM owns a tighter bound than (\ref{thm2con}), and thus it will be faster than Jacobian ADMMs. In the following, we first introduce  M-ADMM, and then discuss the variable partition and backtracking  technique which are crucial for the efficient implementation of M-ADMM in practice. 
		
		\begin{table}[t]
			\centering
			\caption{\small Previous Jacobian ADMMs are special cases of Algorithm \ref{alg2} with different $\hatf_i$ and $\G_i$. In this table, $\eta_i>n\|\A_i\|_2^2$.}
			\small
			\begin{tabular}{c|c|c}\hline
				\multicolumn{1}{c|}{}   & \multicolumn{1}{c|}{$\hatf_i^k(\x_i)$} & $\G_i$                                 \\\hline
				L-ADMM-PS                     & $f_i(\x_i)$                      & $\eta_i\I-\A_i^\top\A_i$,                               \\\hline
				{PL-ADMM-PS}  & Proximal Gradient Surrogate   & {$\eta_i\I-\A_i^\top\A_i$ }    \\\hline
				{GL-ADMM-PS }                    & {$f_i(\x_i)$}                        & {$\succ (n-1)\A_i^\top\A_i$}                 \\\hline
			\end{tabular}\vspace{-0.4cm}\label{tab2}
		\end{table}
		
		\subsection{M-ADMM}\label{subsectionmadmm}
		Assume that we are given a partition of $n$ blocks, denoted as $\{B_1,B_2\}$. 
		We accordingly partition $\A$ into $\Abone=[\A_i, i\in B_1]$ and $\Abtwo=[\A_i, i\in B_2]$. Then (\ref{multipro}) is equivalent to
		\begin{equation}\label{multipro2}
		\min_{\xbone,\xbtwo} f(\x), \ \ \text{s.t.} \ \ \Abone\xbone+\Abtwo\xbtwo=\b.
		\end{equation}
		In the $(k+1)$-th iteration, we first compute the  majorant surrogate of $f$ near $\xk$, i.e., $\hatf^k\in\mathcal{S}_{\{\bL_i,\bP_i\}_{i=1}^n}(f,\xk)$, and $\hatf^k$ is separable,  $\hatf^k(\x)=\hatf^k_{B_1}(\xbone)+\hatf^k_{B_2}(\xbtwo)$, where $\hatf^k_{B_i}(\x_{B_i})=\sum_{j\in B_i}\hatf_j^k(\x_j)$, $i=1,2$. Then (\ref{multipro2}) can be solved by updating $\xbone$ and $\xbtwo$ as the traditional ADMM, i.e.,
		\begin{align}
		\xbonekk=&\argmin_{\xbone} \hatf^k_{B_1}(\xbone)+r^k_{B_1}(\xbone),\label{updatexbone}\\
		\xbtwokk=&\argmin_{\xbtwo} \hatf^k_{B_2}(\xbtwo)+r^k_{B_2 }(\xbtwo),\label{updatexbtwo}
		\end{align}
		where
		$$r_{B_1}^k(\xbone)=\frac{\betak}{2}\left\|\Abone\xbone+\Abtwo\xbtwok-\b+\frac{\lambdak}{\betak}\right\|^2,$$
		and
		$$r^k_{B_2}(\xbtwo)=\frac{\betak}{2}\left\|\Abone\xbonekk+\Abtwo\xbtwo-\b+\frac{\lambdak}{\betak}\right\|^2.$$
		However, the above problems are expensive to solve since they may not be separable w.r.t. $\x_i$'s in $B_1$ or $B_2$. Assume that $\frac{1}{2}\norm{\Abone\xbone}^2$ is $\{\bL'_i\}_{i\in B_1}$-smooth. By   (\ref{sqrsurr2}), we construct a majorant surrogate $\hatr^k_{B_1}$ of $r^k_{B_1}$ near $\xbonek$, i.e., $\hat{r}^k_{B_1}(\xbone)=\sum_{i\in B_1}\hat{r}^k_i(\x_i)$,
		where
		\begin{align}\label{rone}
		\frac{\hat{r}^k_i(\x_i)}{\betak}=&\frac{1}{2}\left\|\A_i\x_i+\sum_{j\in B_1\atop j\neq i }\A_j\x_j^k+\A_{B_2}\x_{B_2}^k-\b+\frac{\lambdak}{\betak}\right\|^2\notag \\&+\frac{1}{2}\|\x_i-\xk_i\|_{\G_i}^2+c_i^k, \ i\in B_1,
		\end{align}
		with $\G_i\succeq \bL'_i-\A_i^T\A_i$, $i\in B_1$, and  $c_i^k$'s satisfying $\sumione c_i^k=\frac{1-n_1}{2}\left\|\A\xk-\b+\frac{\lambdak}{\betak}\right\|^2$.
		Similarly, assume that $\frac{1}{2}\norm{\Abtwo\xbtwo}^2$ is $\{\bL'_i\}_{i\in B_2}$-smooth. Then a majorant surrogate $\hatr_{B_2}^k$ of $r^k_{B_2}$ near $\xbtwok$ is $\hat{r}^k_{B_2}(\xbtwo)=\sum_{i\in B_2}\hat{r}^k_i(\x_i)$, where
		\begin{align}\label{rtwo}
		\frac{\hat{r}^k_i(\x_i)}{\betak}=&\frac{1}{2}\left\|\A_i\x_i+\sum_{j\in B_2 \atop j\neq i }\A_j\xk_j+\Abone\xbonekk-\b+\frac{\lambdak}{\betak}\right\|^2\notag \\ &+\frac{1}{2}\|\x_i-\xk_i\|_{\G_i}^2+c_i^k,\ i\in B_2,
		\end{align}
		with $\G_i\succ  \bL'_i-\A_i^\top\A_i$, $i\in B_2$, and $c_i^k$'s satisfying $\sumitwo c_i^k=\frac{1-n_2}{2}\left\|\Abone\xbonekk+\Abtwo\xbtwok-\b+\frac{\lambdak}{\betak}\right\|^2$. By replacing ${r}^k_{B_1}(\xbone)$ and ${r}^k_{B_2}(\xbtwo)$ with their majorant surrogates $\hat{r}^k_{B_1}(\xbone)$ and $\hat{r}^k_{B_2}(\xbtwo)$ in (\ref{updatexbone}) and (\ref{updatexbtwo}) respectively, we update $\xbone$ and $\xbtwo$ by
		\begin{align*}
		\xbonekk=&\argmin_{\xbone} \hat{f}^k_{B_1}(\xbone)+\hat{r}^k_{B_1}(\xbone),\\
		\xbtwokk=&\argmin_{\xbtwo} \hat{f}^k_{B_2}(\xbtwo)+\hat{r}^k_{B_2}(\xbtwo).
		\end{align*}
		Note that the above two problems are separable for each $\x_i$ in $B_1$ and $B_2$. They are respectively equivalent to
		\begin{align}
		\xkk_i=&\argmin_{\x_i} \hat{f}^k_i(\x_i)+\hat{r}^k_i(\x_i), \ i\in B_1,\label{updatexbone2}\\
		\xkk_i=&\argmin_{\x_i} \hat{f}^k_i(\x_i)+\hat{r}^k_i(\x_i), \ i\in B_2.\label{updatexbtwo2}
		\end{align}
		Finally $\blambda$ is updated by (\ref{updatelambdaall}). This leads to the Mixed Gauss-Seidel and Jacobian ADMM (M-ADMM), as shown in Algorithm \ref{alg3}. Now we give its convergence bound as (\ref{lem3eq12}).

		\begin{algorithm}[t]
			\caption{Mixed Gauss-Seidel and Jacobian ADMM}
			\textbf{For} $k=0,1,2,\cdots$ \textbf{do}
			\begin{enumerate}
				\item Compute a majorant first-order surrogate $\hat{f}^k\in\mathcal{S}_{\{\bL_i, \bP_i\}_{i=1}^n}(f,\xk)$ with $\hatf^k(\x)=\sum_{i=1}^{n}\hatf_i^k(\x_i)$.
				\item For all $i\in B_1$, update $\x_i$'s in parallel by solving (\ref{updatexbone2}).
				\item For all $i\in B_2$, update $\x_i$'s in parallel by solving (\ref{updatexbtwo2}).
				\item Update $\bm{\lambda}$ by $\blambda^{k+1}=\lambdak+\betak(\A\x^{k+1}-\b)$.
				\item Choose $\betakk\geq\betak$.
			\end{enumerate}
			\textbf{end}
			\label{alg3}
		\end{algorithm}

		\begin{theorem}\label{them3}
			In Algorithm \ref{alg3}, assume that $\hatf^k\in\mathcal{S}_{\{\bL_i,\bP_i\}_{i=1}^n}(f,\xk)$ with $\bP_i\succeq \bL_i\succeq \bzero$, $\frac{1}{2}\norm{\Abone\xbone}^2$ is $\{\bL'_i\}_{i\in B_1}$-smooth, $\frac{1}{2}\norm{\Abtwo\xbtwo}^2$ is $\{\bL'_i\}_{i\in B_2}$-smooth,   $\G_i\succeq\bL'_i-\A_i^\top\A_i$,   $i\in B_1$ in (\ref{rone}) and $\G_i\succ \bL'_i-\A_i^\top\A_i$,   $i\in B_2$ in (\ref{rtwo}).
			For any $K>0$, let $\bar{\x}^K=\sum_{k=0}^{K}\gammak\xkk$ with $\gammak={(\betak)^{-1}}/{\sum_{k=0}^{K}(\betak)^{-1}}$. Then
			\begin{align}
			&f(\bar{\x}^K)-f(\x^*)+\langle \A^\top\blambda^*,\bar{\x}^K-\x^*\rangle+\frac{\beta^{(0)}\alpha}{2}\|\A\bar{\x}^K-\b\|^2\notag\\
			\leq&\frac{\sum_{j=1}^2\|\x_{B_j}^*-\x_{B_j}^0\|^2_{\HH^0_j}+\|\blambda^*-\blambda^0\|^2_{\HH^0_3}}{2\sum_{k=0}^{K}\left(\betak\right)^{-1}},\label{thm4con}
			\end{align}
			where $\alpha=\min\left\{\frac{1}{2},\frac{\sigma^2_{\min}\left(\text{Diag}\left\{\A_i^\top\A_i+\G_i, i\in B_2\right\}-\Abtwo^\top\Abtwo\right)}{2\|\Abtwo\|^2_2}\right\}$,
			$\HH^0_1=\text{Diag}\left\{\frac{1}{\beta^{(0)}}\bL_i+\A_i^\top\A_i+\G_i,i\in B_1\right\}-\A_{B_1}^\top\A_{B_1}$, 
			$\HH_2^0=\text{Diag}\left\{\frac{1}{\beta^{(0)}}\bL_i+\A_i^\top\A_i+\G_i, i\in B_2\right\}, \text{ and}$
			$\HH^0_3=\left({1/\beta^{(0)}}\right)^2\I$.
		\end{theorem}
		
		M-ADMM in Algorithm \ref{alg3} further unifies Gauss-Seidel ADMMs in Algorithm \ref{alg1} and Jacobian ADMMs in Algorithm \ref{alg2}. If $n=2$, let $B_1=\{1\}$ and $B_2=\{2\}$. Then M-ADMM degenerates into the Gauss-Seidel ADMMs, and (\ref{thm4con}) reduces to (\ref{thm1con}).
		If $n>2$, let $B_1=\varnothing$ and $B_2=\{1,\cdots,n\}$. Then M-ADMM degenerates into the Jacobian ADMMs, and (\ref{thm4con}) reduces to (\ref{thm2con}).
		More importantly, for the case of $n>2$ and other choices of $B_1$ and $B_2$, M-ADMM will be faster than Jacobian ADMMs, since the right hand side of (\ref{thm4con}) may be much smaller than the one of (\ref{thm2con}). This is due to their different choices of $\G_i$. Without the additional assumption on the structure of $\A$, Jacobian ADMMs require $\G_i\succ (n-1)\A_i^\top\A_i$ for all $\alli$, while M-ADMM only requires $\G_i\succ(n_1-1)\A_i^\top\A_i$ for $i\in B_1$ and $\G_i\succ(n_2-1)\A_i^\top\A_i$ for $i\in B_2$. Note that $n=n_1+n_2$. The improvement benefits from the sequential updating rules of $\x_{B_1}$ and $\x_{B_2}$ by using tighter majorant surrogates in M-ADMM. Indeed, M-ADMM only needs to majorize  $r^k_{B_1}(\x_{B_1})$ in (\ref{updatexbone}) and $r^k_{B_2}(\x_{B_2})$ in (\ref{updatexbtwo}) for $\x_{B_1}$ and $\x_{B_2}$ respectively, while Jacobian ADMMs need to majorize $r^k(\x)$ in (\ref{almr}) for all $\x_i$'s simultaneously. 
		
		Note that the work \cite{he2015blockk} proposes a block-wise ADMM which is another special case of our Algorithm \ref{alg3}. But their considered problem is more specific and the convergence guarantee requires much stronger assumptions, e.g., $\A_i$ that has a full column rank.
		
			\begin{algorithm}[t]
				\caption{M-ADMM with backtracking}
				\textbf{Initialization:} $k=0$, $\x_i^k$, $\Gk_i\succ \bm0$, $\blambda^k$,  $\beta^k>0$, $\tau>0$, $\mu>1$.\\
				\textbf{For} $k=0,1,2,\cdots$ \textbf{do}
				\begin{enumerate}
					\item Compute a majorant first-order surrogate $\hat{f}^k\in\mathcal{S}_{\{\bL_i, \bP_i\}_{i=1}^n}(f,\xk)$ with $\hatf^k(\x)=\sum_{i=1}^{n}\hatf_i^k(\x_i)$.
					\item Set $\G_i=\Gk_i$ and compute $\xkk_i$ by (\ref{updatexbone2})-(\ref{updatexbtwo2}).
					\item If (\ref{eqbteq2}) does not hold, set $\Gk_i = \mu \Gk_i$, $i\in B_1$. Go to 2).
					\\
					If  (\ref{eqbteq3}) does not hold, set $\Gk_i = \mu \Gk_i$, $i\in B_2$. Go to 2).
					\item Update $\bm{\lambda}$ by $\blambda^{k+1}=\lambdak+\betak(\A\x^{k+1}-\b)$.
					\item Choose $\betakk\geq\betak$. Set $\Gkk_i=\Gk_i$, $\alli$.
				\end{enumerate}
				\textbf{end}
				\label{alg4}
			\end{algorithm}

	\subsection{M-ADMM with Backtracking}
	We have given the convergence guarantee of M-ADMM when fixing  $\G_i$. In practice, we can estimate it by the backtracking technique  which will lead to tighter majorant surrogate. The effectiveness has been verified   in first-order optimization \cite{beck2009fast}. Now, we introduce the  backtracking technique into M-ADMM.
	
	To guarantee the convergence,  $\G_i$ can be replaced by $\Gk_i$ such that   $r^k_{B_1}(\xbone^{k+1}) \leq \hat{r}^k_{B_1}(\xbone^{k+1})$ and $r^k_{B_2}(\xbtwo^{k+1}) \leq \hat{r}^k_{B_2}(\xbtwo^{k+1})$. They are guaranteed when   
	\begin{equation}\label{eqbteq2}
	\norm{\Abone(\xbonekk-\xbonek)}^2 \leq \sumione\norm{\xkk_i-\xk_i}_{\Gk_i+\A_i^\top\A_i}^2,
	\end{equation}
	\begin{equation}\label{eqbteq3copy1}
	\norm{\Abtwo(\xbtwokk-\xbtwok)}^2 \leq \sumitwo\norm{\xkk_i-\xk_i}_{\Gk_i+\A_i^\top\A_i}^2.
	\end{equation}
	To achieve the $O(1/K)$ convergence rate, we replace (\ref{eqbteq3copy1}) as
	\begin{equation}\label{eqbteq3}
	\tau \norm{\xbtwo^{k+1}-\xbtwo^k}^2 \leq \norm{\xbtwo^{k+1}-\xbtwo^k}^2_{\Kk_2-\Abtwo^\top\Abtwo},
	\end{equation}		
	for some small constant $\tau>0$ and  $\Kk_2=\text{Diag}\{\A_i^\top\A_i+\Gk_i,i\in B_2\}$.
	In this case, we may be able to find $\Gk_i$ with relatively smaller $\norm{\Gk_i}_2$, and thus  $\hat{r}^k_{B_1}(\xbone^{k+1})$ and $\hat{r}^k_{B_2}(\xbtwo^{k+1})$ are tighter upper bounds of $r^k_{B_1}(\xbone^{k+1})$ and $r^k_{B_2}(\xbtwo^{k+1})$, respectively. This  leads to a better approximation solution and improves the efficiency.
	We summarize M-ADMM with backtracking in Algorithm \ref{alg4}. Note that Step 3) will only be performed for finitely many times.  Similarly, the convergence guarantee is given as follows.
	\begin{theorem}\label{them4}
		In Algorithm \ref{alg4}, assume that $\hatf^k\in\mathcal{S}_{\{\bL_i,\bP_i\}_{i=1}^n}(f,\xk)$ with $\bP_i\succeq \bL_i\succeq \bzero$.
		Then (\ref{thm4con}) holds with
		$\HH^0_1=\text{Diag}\left\{\frac{1}{\beta^{(0)}}\bL_i+\A_i^\top\A_i+\G^0_i,i\in B_1\right\}-\A_{B_1}^\top\A_{B_1}$, 
		$\HH_2^0=\text{Diag}\left\{\frac{1}{\beta^{(0)}}\bL_i+\A_i^\top\A_i+\G^0_i, i\in B_2\right\}$,
		$\HH^0_3=\left({1/\beta^{(0)}}\right)^2\I$, and $\alpha=\min\left\{\frac{1}{2},\frac{\tau}{2\norm{\Abtwo}_2^2}\right\}$.
	\end{theorem}
	Note that Algorithm \ref{alg4} reduces to Algorithm \ref{alg3} by choosing $\G_i$'s in Theorem \ref{them3}. Theorem \ref{them3} is a special case of Theorem \ref{them4} by setting $\tau=\sigma^2_{\min}\left(\text{Diag}\left\{\A_i^\top\A_i+\G_i, i\in B_2\right\}-\Abtwo^\top\Abtwo\right)$. So we only give the proof of Theorem \ref{them4} in   Appendix.  	It is worth mentioning that, when using backtracking, $\hatr^k_{B_j}$ is not a majorant first-order surrogate of $r^k_{B_j}$, since the  majorization condition may not hold. Actually,  $\hatr^k_{B_j}$ only majorizes $r^k_{B_j}$ locally at $\x_{B_j}^{k+1}$. But this is sufficient for the convergence proof, since the formulations of $r^k_{B_j}$ and $\hatr^k_{B_j}$ are known and we are able to use their  specific properties instead of (\ref{keylm22}) in the proofs.


	
		\subsection{Variable Partition}\label{sec_varpat}
			
			
			For (\ref{multipro}) with $n>2$, M-ADMM requires   partitioning variables into 2 super blocks $B_1$ and $B_2$. Different variable partitions lead to different choices of $\bL'_i$ which controls the tightness of the majorant surrogates, and thus the convergence behaviors of M-ADMM are different. 			
			Looking for an intelligent way of variable partition may significantly improve the efficiency of M-ADMM. We discuss how to partition variables in three cases by considering the property of $\A_i$'s in (\ref{multipro}). The principle is to find a partition such that the constructed surrogate $\hatr_{B_1}^k$ for $r_{B_1}^k$   in (\ref{updatexbone}) and $\hatr_{B_2}^k$ for  $r_{B_2}^k$ in (\ref{updatexbtwo}) can be as tight as possible. 
			
			 \textbf{Case I (complex case):}  $\A_i^\top\A_l\neq\bm{0}$ for any $i\neq l$. This case is complex since $r_{B_j}^k$, $j=1,2$   in (\ref{updatexbone})-(\ref{updatexbtwo})  are   non-separable for any partition.  Then the separable surrogates $\hatr_{B_j}^k$'s  may be loose when considering   the choices of $\G_i$ in (\ref{rone})-(\ref{rtwo}).
			As suggested by   Theorem \ref{them3}, to tighten the bound of (\ref{thm4con}), a reasonable partition is to make $L_{B_1}+L_{B_2}$, where $L_{B_1}=(n_1-1)\sum_{i\in B_1}\|\A_i\|^2_2-\norm{\Abone}_2^2$ and $L_{B_2}=(n_2-1)\sum_{i\in B_2}\|\A_i\|^2_2$, as small as possible\footnote{If $n_1$ is not very small, $\norm{\Abone}_2^2$ is usually much smaller than $(n_1-1)\sum_{i\in B_1}\|\A_i\|^2_2$. We can use $L_{B_1}=(n_1-1)\sum_{i\in B_1}\|\A_i\|^2_2$ in this case.}. 
			We have a heuristic approach to  this end:
			
			\noindent Step 1: Sort $\|\A_i\|_2^2$'s in a descending order.
			
			\noindent Step 2: Group the largest $n_1$ elements as the first block and the rest as the second block.
			
			\noindent Step 3: The best value of $n_1$ is the one which minimizes $L_{B_1}+L_{B_2}$ by a one-shot searching from 1 to $n$.

			\textbf{Case II (simple case):} there exists a partition such that 
			\begin{equation}\label{sepcondi}
			\text{$\A_i^\top \A_l = \bm{0}$, $i\neq l$, for any $i, l \in B_1$ and $i, l \in B_2$.}
			\end{equation}
			This case is simple since the above partition makes $r_{B_j}^k$, $j=1,2$   in (\ref{updatexbone})-(\ref{updatexbtwo})     separable. Then $\hatr_{B_j}^k$'s tend to be tight since we can compute each $\hatr_i^k$ independently and use $\G_i\succeq \bm{0}$, $i\in B_1$ in (\ref{rone}) and $\G_i\succ \bm{0}$, $i\in B_2$ in (\ref{rtwo}). Even, the per-iteration cost is cheap when using $\hatr^k_i=r^k_i$ for many problems  in practice. In this case, (\ref{multipro2}) can be solved by (\ref{updatexbone})-(\ref{updatexbtwo}), which is similar to the standard ADMM. For example, the Low-Rank Representation  model in \cite{robustlrr}  satisfies (\ref{sepcondi}),  
			\begin{equation}\label{lrrprobexmaple}
			\min_{\Z,\J,\E} \norm{\J}_*+\lambda\norm{\E}_{2,1}, \ \text{s.t.} \ \X=\A\Z+\E, \Z=\J,
			\end{equation}
			where $\lambda>0$. The augmented Lagrangian function is
			\begin{align*}
			\mathcal{L}&(\Z,\J,\E,\blambda_1,\blambda_2)=  \norm{\J}_*+\lambda\norm{\E}_{2,1} + \langle \blambda_1, \X-\A\Z-\E\rangle  \\
			&+ \langle\blambda_2,\Z-\J\rangle+\frac{\beta}{2}\left(\norm{\X-\A\Z-\E}^2+\norm{\Z-\J}^2\right). 
			\end{align*}
			Based on the partition $\{\J,\E\}$ and $\{\Z\}$, they can be updated by  
			\begin{equation}
			\left\{
			\begin{aligned}
			\{\J^{k+1},\E^{k+1}\}=&\arg\min_{\J,\E}  \mathcal{L}(\Z^k,\J,\E,\blambda^k_1,\blambda^k_2),    \\
			\Z^{k+1}=&\arg\min_{\Z} \mathcal{L}(\Z,\J^{k+1},\E^{k+1},\blambda^k_1,\blambda^k_2).
			\end{aligned}
			\right.\notag
			\end{equation}
			This is the standard ADMM and its convergence is guaranteed.
			Note that $\mathcal{L}(\Z^k,\J,\E,\blambda^k_1,\blambda^k_2)$ is separable w.r.t. $\J$ and $\E$ and thus $\J^{k+1}$ and $\E^{k+1}$ can be computed independently. The updates of the three blocks are  similar to the  naive multi-block extension of ADMM used in \cite{robustlrr}, but in different updating orders. Our simple modification fixes the convergence issue of the naive multi-block extension of ADMM  in \cite{robustlrr} for (\ref{lrrprobexmaple}). 
 			
			In computer vision and signal processing, there are a lot of multi-blocks problems, or their equivalent ones by introducing auxiliary variables, with the property (\ref{sepcondi}) and thus can be solved more efficiently by the Gauss-Seidel ADMMs than Jacobian ADMMs, e.g., sparse subspace clustering model (70) in \cite{6482137}, nonnegative matrix completion problem (143) in \cite{LADMPS}, multi-task low-rank affinity pursuit model (4) in  \cite{cheng2011multi}, sparse spectral clustering model (6) in \cite{lu2016convex}, nonnegative low-rank and sparse graph model (5) in \cite{zhuang2012non}, simultaneously structured models (3.3) in \cite{oymak2015simultaneously}, convex program (8) in \cite{chen2014improved} for  graph clustering, robust multi-view spectral clustering model (3) in \cite{xia2014robust}    and consolidated tensor recovery model (2.6) in \cite{HMWG15}. However, some of previous works do not use the property (\ref{sepcondi}) to implement the efficient ADMMs, and this is the reason why we release the toolbox.
			
			\textbf{Case III (other cases):} neither assumptions in Case I and Case II holds. It is generally difficult to find the best partition in this case. But one can combine the ideas in both Case I and II. For example, there exists one or  more subgroups $B_{S}$, such that $\A_i^\top\A_l=\bm{0}$, $i\neq l$, for any $i,l\in B_S$.  We can put the whole subgroup in one super block, i.e., $B_S\subset B_1$. 
			
			In practice, one usually needs to reformulate the original problem as an equivalent one by introducing   auxiliary variables such that the subproblem in ADMMs can be simple. When designing efficient ADMMs, the problem reformulation and the above variable partition strategies should be considered simultaneously. Some more examples can be found in our released toolbox.

		\section{Experiments}
		\label{sec6}
		In this section, we conduct several experiments to show the effectiveness of our new ADMMs. All the algorithms are implemented by Matlab and are tested on a PC with 8 GB of RAM and Intel Core 2 Quad CPU Q9550.  The details of the compared solvers can be found in the supplementary material.
		
		\subsection{Experimental Analysis of M-ADMM}
		Besides the unified analysis of several variants of ADMMs, another main contribution of this work is the proposed M-ADMM for multi-block problems. In this subsection, our purpose is to perform some analyses on M-ADMM. For the simplicity, we first consider  the following   
		nonnegative sparse coding problem 
		\begin{align}\label{nonsparse11}
			 \min_{\{\x_i\}}\sumi\|\x_i\|_1, \ \text{s.t. }  \y=\sumi\A_i\x_i, \x_i\geq\bzero,
		\end{align}
		
				\begin{figure}
					\centering
					\begin{subfigure}[b]{0.24\textwidth}
						\centering
						\includegraphics[width=\textwidth]{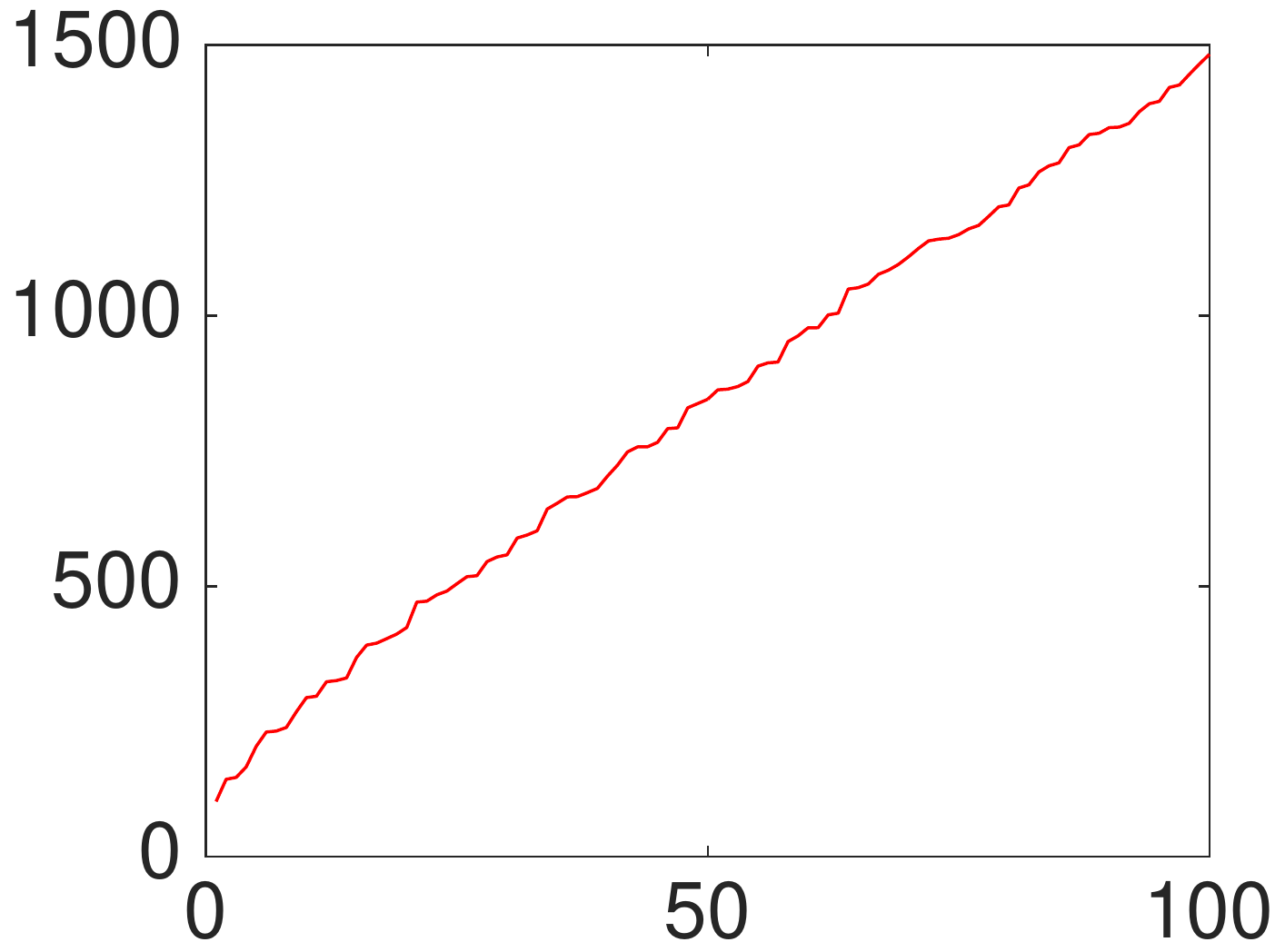} 
						\caption{}
					\end{subfigure}
					\begin{subfigure}[b]{0.23\textwidth}
						\centering
						\includegraphics[width=\textwidth]{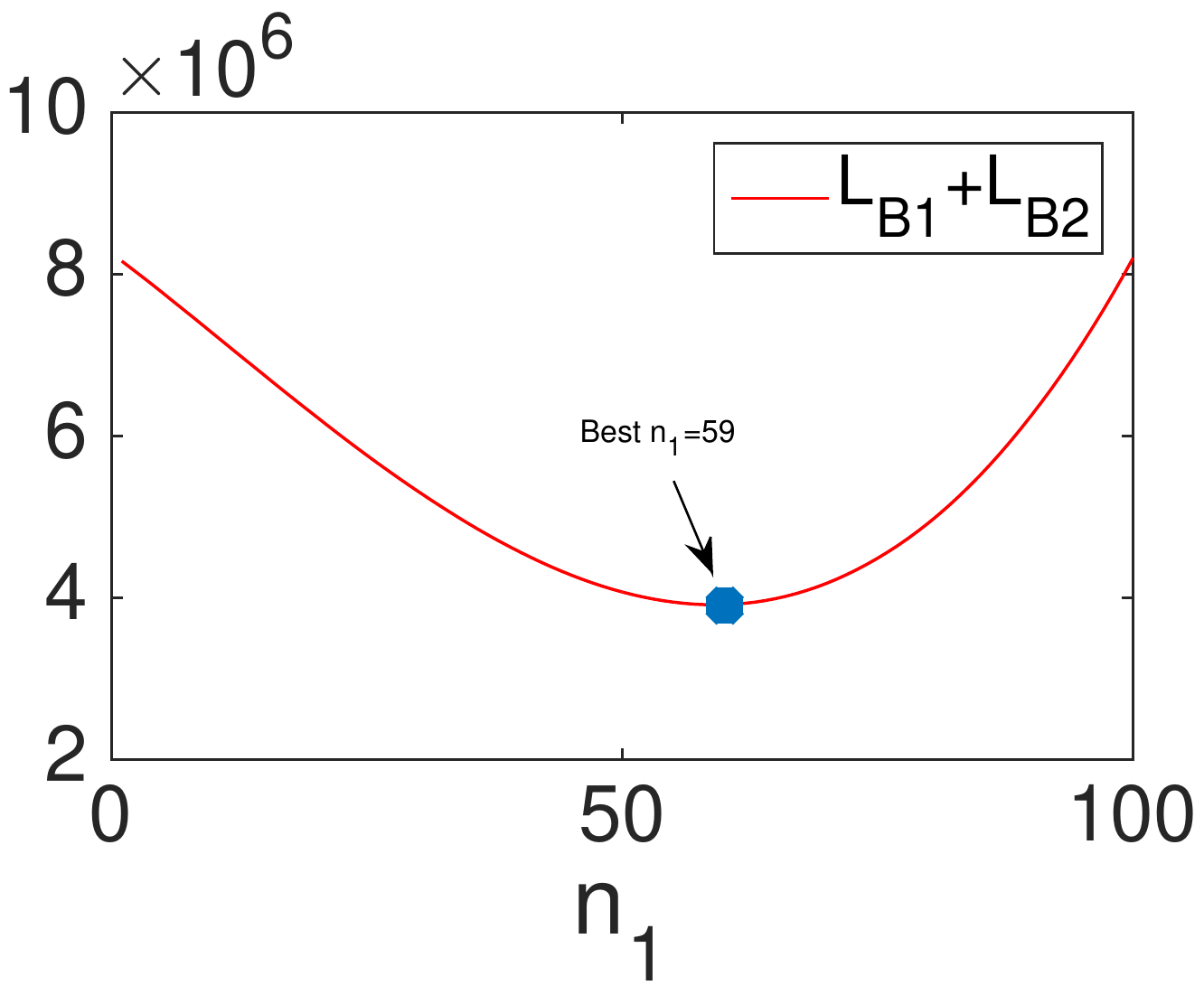} 
						\caption{}
					\end{subfigure}
					\begin{subfigure}[b]{0.24\textwidth}
						\centering
						\includegraphics[width=\textwidth]{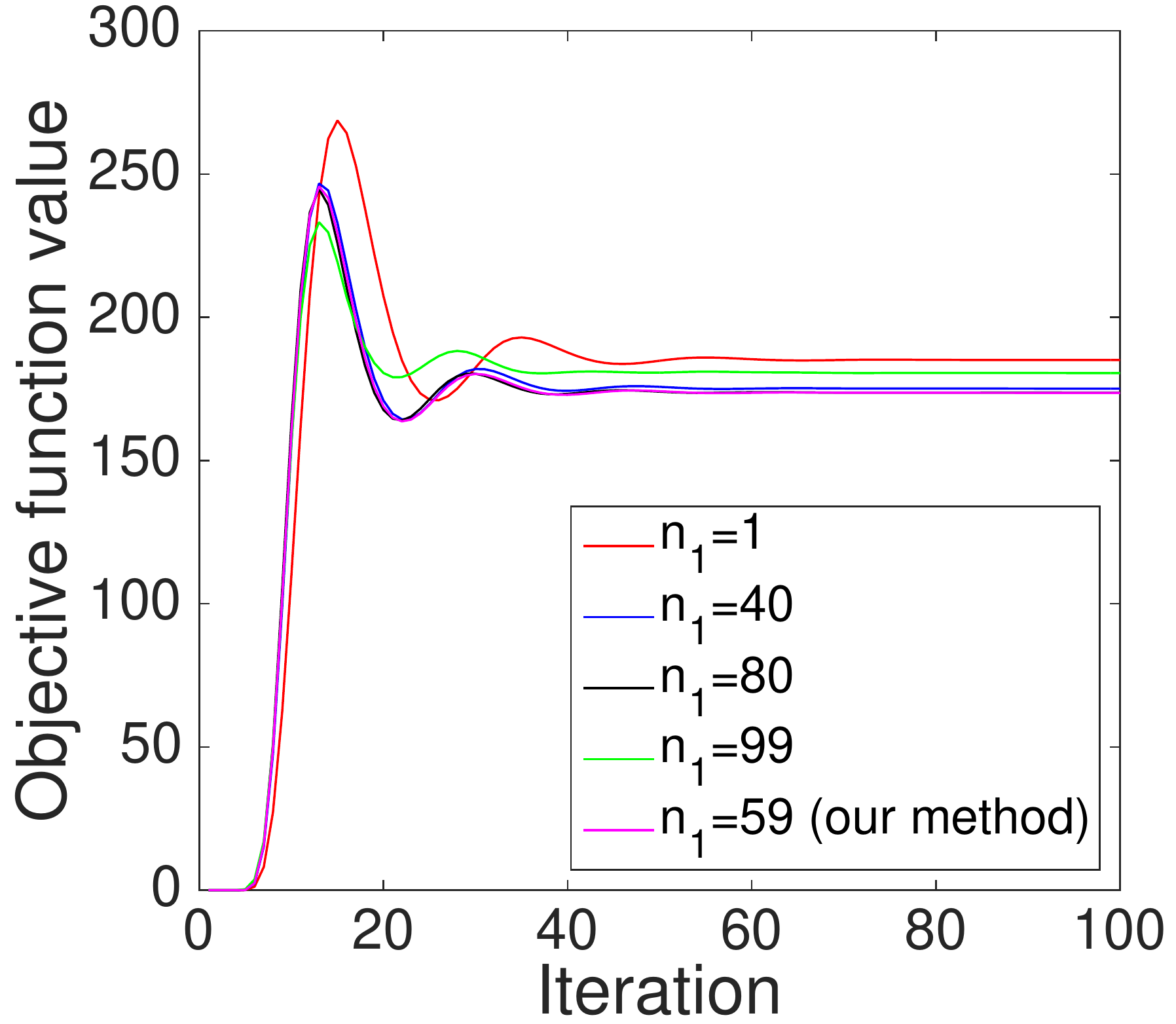} 
						\caption{\small }
					\end{subfigure} 
					\begin{subfigure}[b]{0.24\textwidth}
						\centering
						\includegraphics[width=\textwidth]{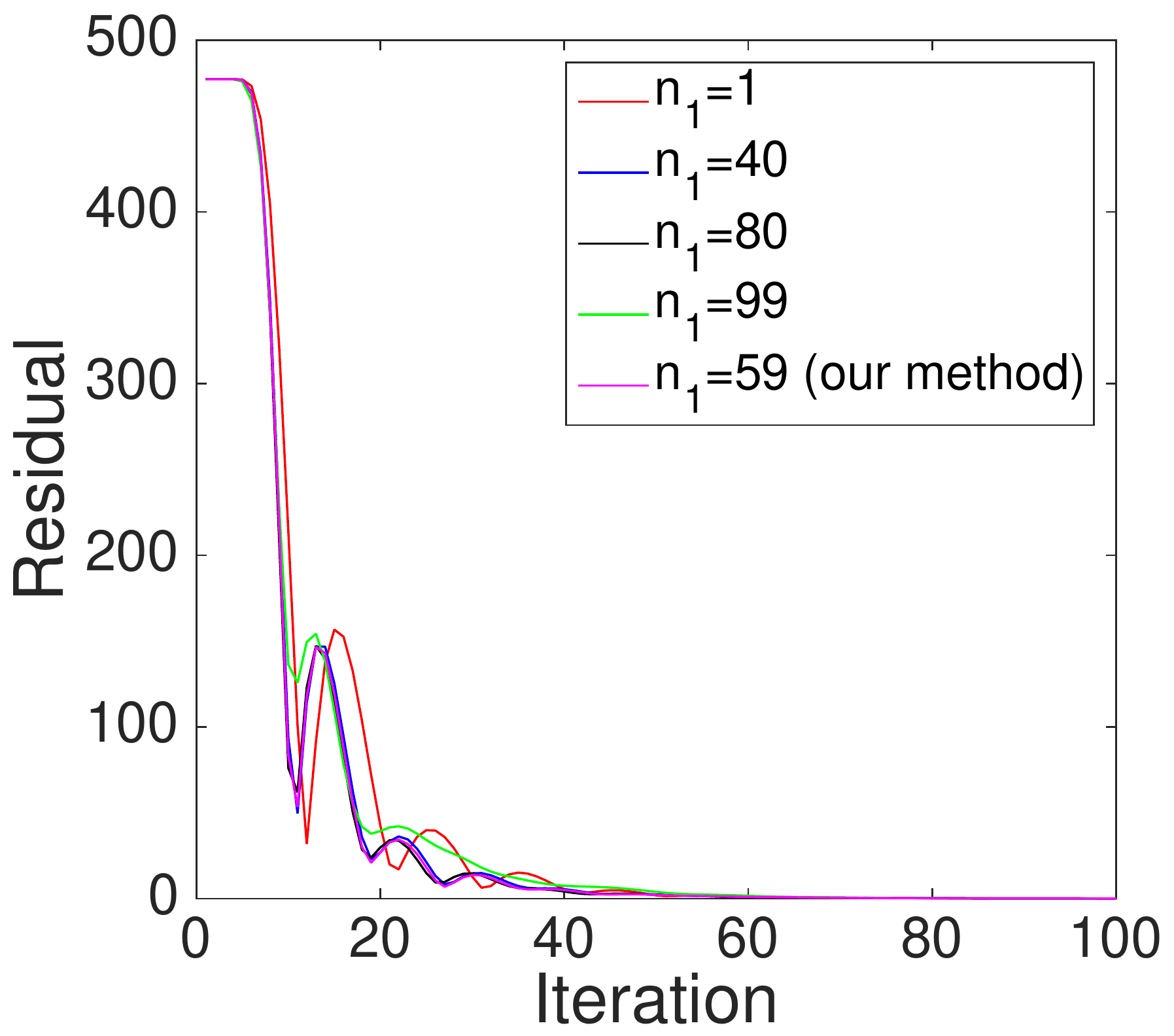} 
						\caption{\small }
					\end{subfigure} 
					\caption{\small Plots of (a) sorted $\{\|\A_i\|^2_2,i=1,\cdots,n\}$; (b) $L_{B_1}+L_{B_2}$ v.s. $n_1$; (c) $f(\xk)$  v.s. $k$ and (d) $\norm{\Axk-\b}$ v.s. $k$ based on different partitions (corresponding to different $n_1$).}\label{fig_blockpartition}
					\vspace{-0.4cm}
				\end{figure}

		\subsubsection{Analysis of the Proposed Partition Strategy}
		\label{sec_varpartstra}
		We conduct an experiment to compare the different convergence behaviors of M-ADMM with different variable partitions and demonstrate the effectiveness of the proposed partition method for Case I  in Section \ref{sec_varpat}. By   choosing $\G_i \succeq \eta_i\I-\A_i^\top\A_i$ in (\ref{rone}) and (\ref{rtwo}), M-ADMM solves (\ref{nonsparse11}) by the following rules		
		 	\begin{equation*}
		 	\left\{
		 	\begin{aligned}
		 	\x_i^{k+1}=&\arg\min_{\x_i\geq \bzero} \|\x_i\|_1+\frac{\betak\eta_i}{2}\left\|\x_i- \u_i^k\right\|^2, \ i\in B_1, \\
		 	\x_i^{k+1}=&\arg\min_{\x_i\geq \bzero} \|\x_i\|_1 	+\frac{\betak\eta_i}{2}\left\|\x_i-\v_i^k\right\|^2, \ i\in B_2, \\
		 	\lambdakk=&\lambdak+\betak(\A\xkk-\y),
		 	\end{aligned}
		 	\right.
		 	\end{equation*}
		 	where $\u_i^k = \x_i^k-\frac{\A_i^T(\blambda^k+\betak(\A\xk-\y))}{\betak\eta_i}$,
		 	$\v_i^k = \x_i^k-\frac{\A_i^T(\blambda^k+\betak(\A_{B_1}\xkk_{B_1}+\A_{B_2}\xk_{B_2}-\y))}{\betak\eta_i}$, $\eta_i\geq n_1\|\A_i\|_2^2$, $i\in B_1$, and $\eta_i>n_2\|\A_i\|_2^2$, $i\in B_2$. In M-ADMM, $\x_i$ and $\blambda$ are initialized as zeros. We set $\beta^{(0)}=10^{-4}$ and update $\betakk=\min(\rho\betak,10^6)$ with $\rho=1.1$. Let $\eta_i = n_1\|\A_i\|_2^2$ for $i\in B_1$, and  $\eta_i=1.02 n_2\|\A_i\|_2^2$ for $i\in B_2$. We test M-ADMM for (\ref{nonsparse11}) on the synthetic data generated as follows. We set $n=100$,  $d=50$, $m_i=10i$ and the elements of $\A_i\in\mathbb{R}^{d\times m_i}$ are independently sampled from an $N(0, 1)$ distribution. We generate $\x$ with $90\%$ elements being zeros  and others independently sampled from an $N(0, 1)$ distribution. Then $\y=[\A_1,\cdots,\A_n ]\x$.   The sizes of $\A_i$'s are different, and so are the Lipschitz constants $\|\A_i\|_2^2$'s. We plot the sorted $\|\A_i\|_2^2$'s in Figure~\ref{fig_blockpartition}~(a).
		M-ADMM requires   dividing these $n$ blocks of variables into two super blocks, i.e., $\xbone$ with $n_1$ blocks, and $\xbtwo$ with $n_2$ blocks. Our partition strategy finds $n_1$ by minimizing $L_{B_1}+L_{B_2}$, where $L_{B_1}=(n_1-1)\sum_{i\in B_1}\|\A_i\|^2_2-\norm{\Abone}_2^2$ and $L_{B_2}=(n_2-1)\sum_{i\in B_2}\|\A_i\|^2_2$. In this experiment, our method gives the best $n_1=59$. See the plot of $L_{B_1}+L_{B_2}$ v.s. $n_1$ in Figure~\ref{fig_blockpartition}~(b). Note that one may have many other choices of $n_1\in\{1,2\cdots,100\}$. Figure \ref{fig_blockpartition} (c) plots the objective function value $f(\xk)$ v.s. iteration $k\ (\leq 100)$ and Figure \ref{fig_blockpartition} (d) plots the residual $\norm{\A\xk-\b}$  v.s. iteration $k\ (\leq 100)$,  based on different choices of $n_1\in\{1,40,80,99,69\}$. Generally, the convergences of $\norm{\A\xk-\b}$ based on different partitions are quite similar since it heavily depends on the same updating rule of $\betakk$. However, different $n_1$ leads to quite different convergences of $f(\xk)$, and $n_1=59$, predicted by our method, performs well. The choice of $n_1=1$ is the worst case since $L_{B_1}+L_{B_2}$ is the largest. These results verify that M-ADMM converge faster when using our proposed variable partition strategy, which leads to tight majorant surrogates.

		
				\begin{figure}
					\centering
					\begin{subfigure}[b]{0.24\textwidth}
						\centering
						\includegraphics[width=\textwidth]{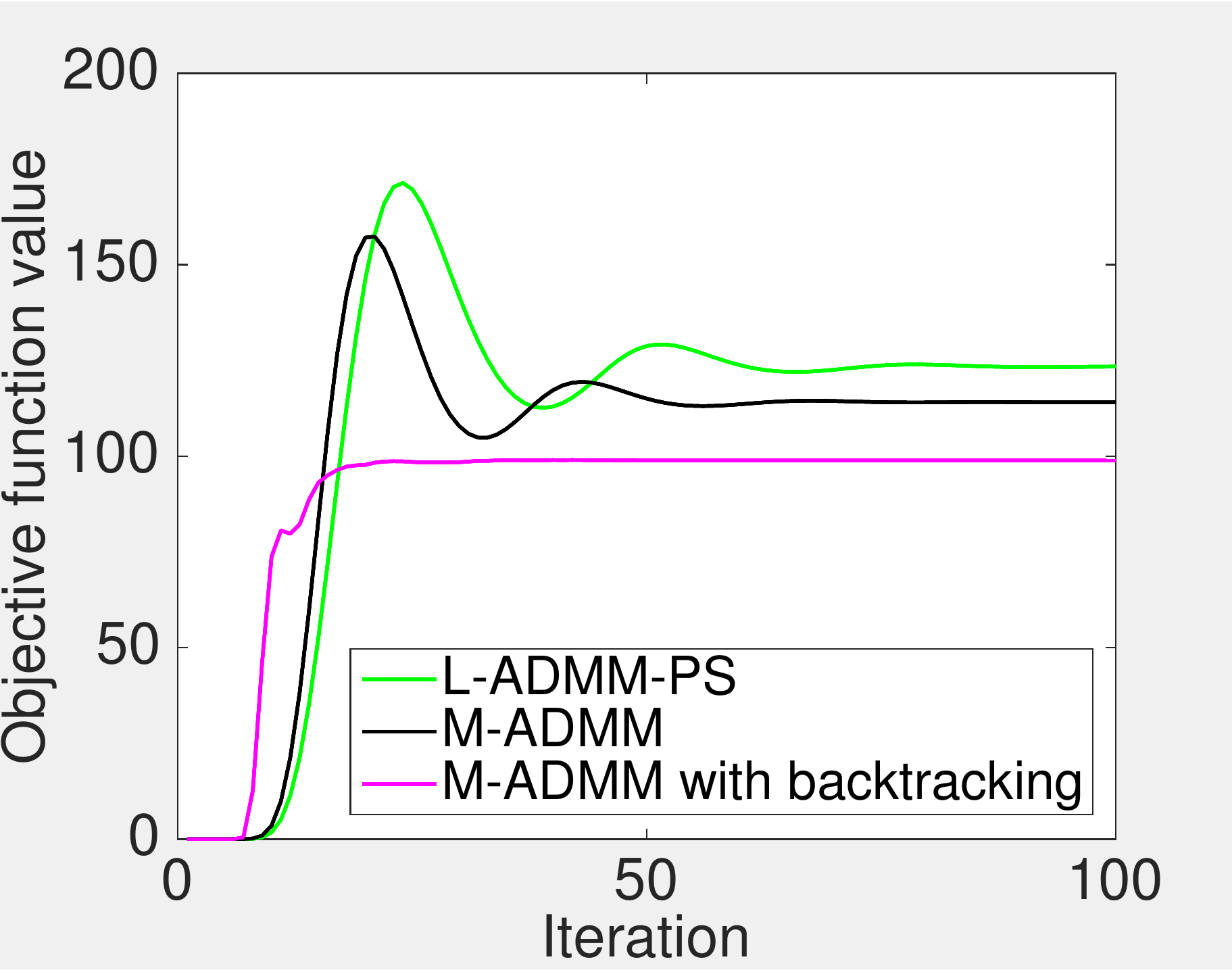} 
					\end{subfigure}
					\begin{subfigure}[b]{0.24\textwidth}
						\centering
						\includegraphics[width=\textwidth]{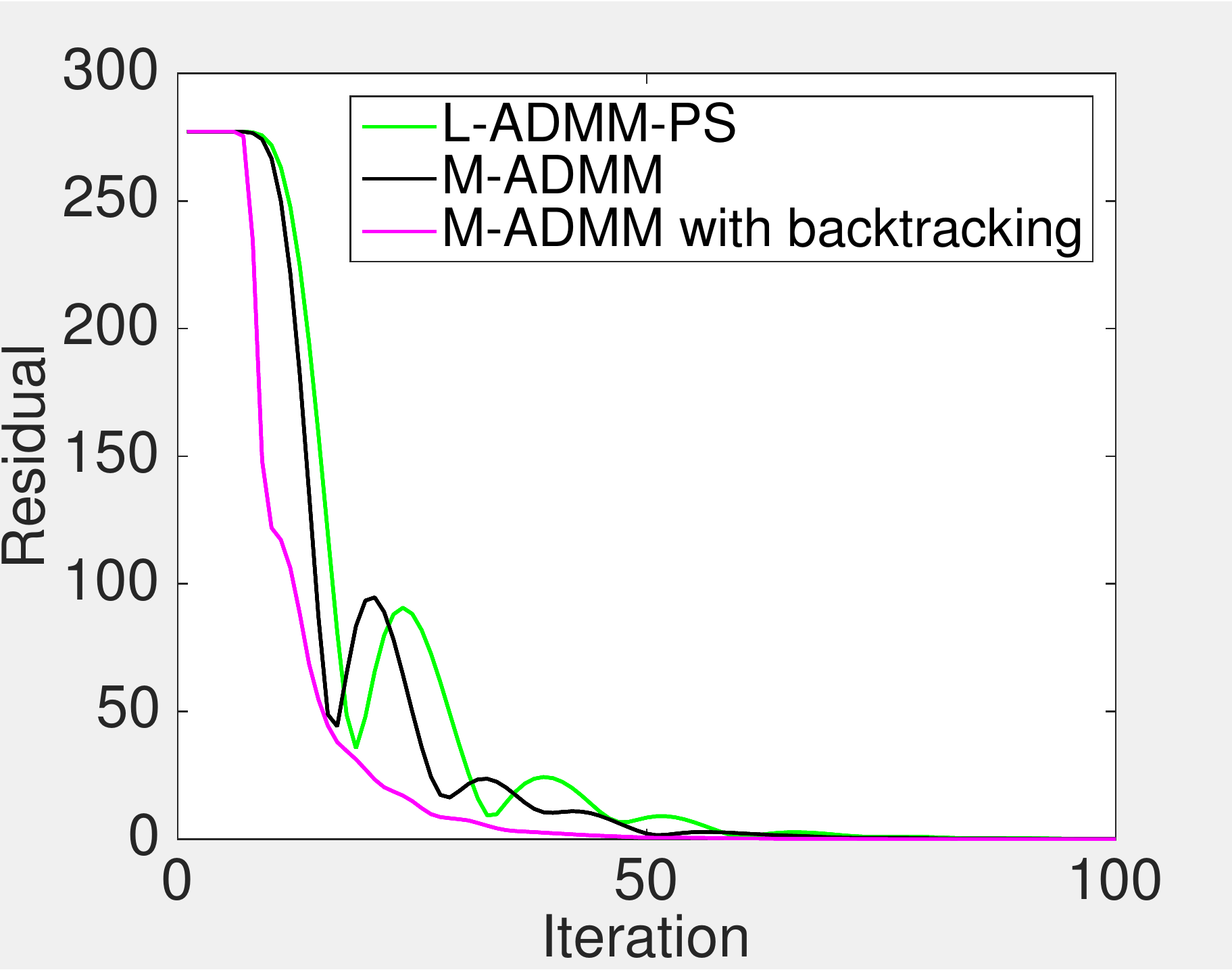} 
					\end{subfigure}
					\begin{subfigure}[b]{0.155\textwidth}
						\centering
						\includegraphics[width=\textwidth]{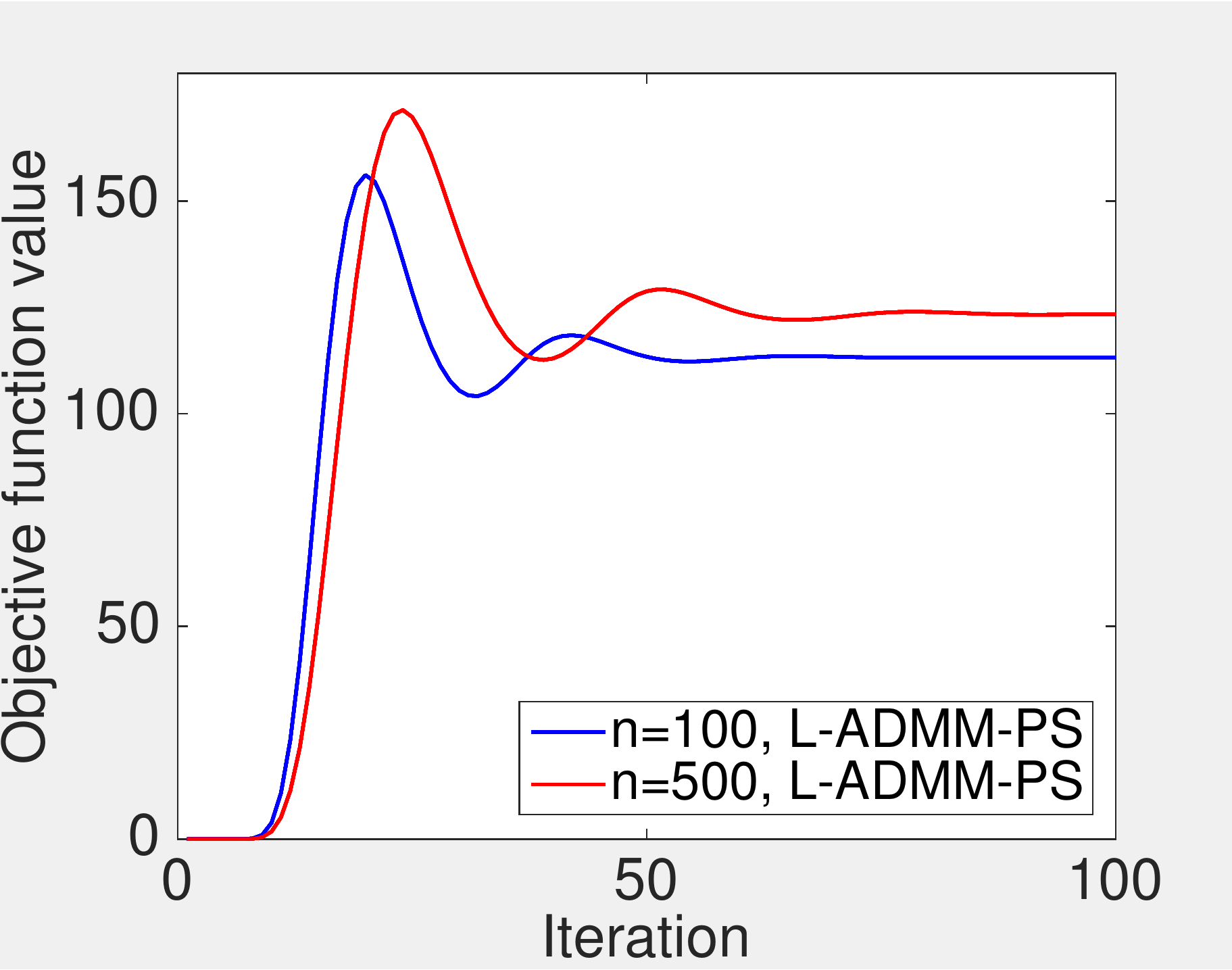} 
					\end{subfigure} 
					\begin{subfigure}[b]{0.155\textwidth}
						\centering
						\includegraphics[width=\textwidth]{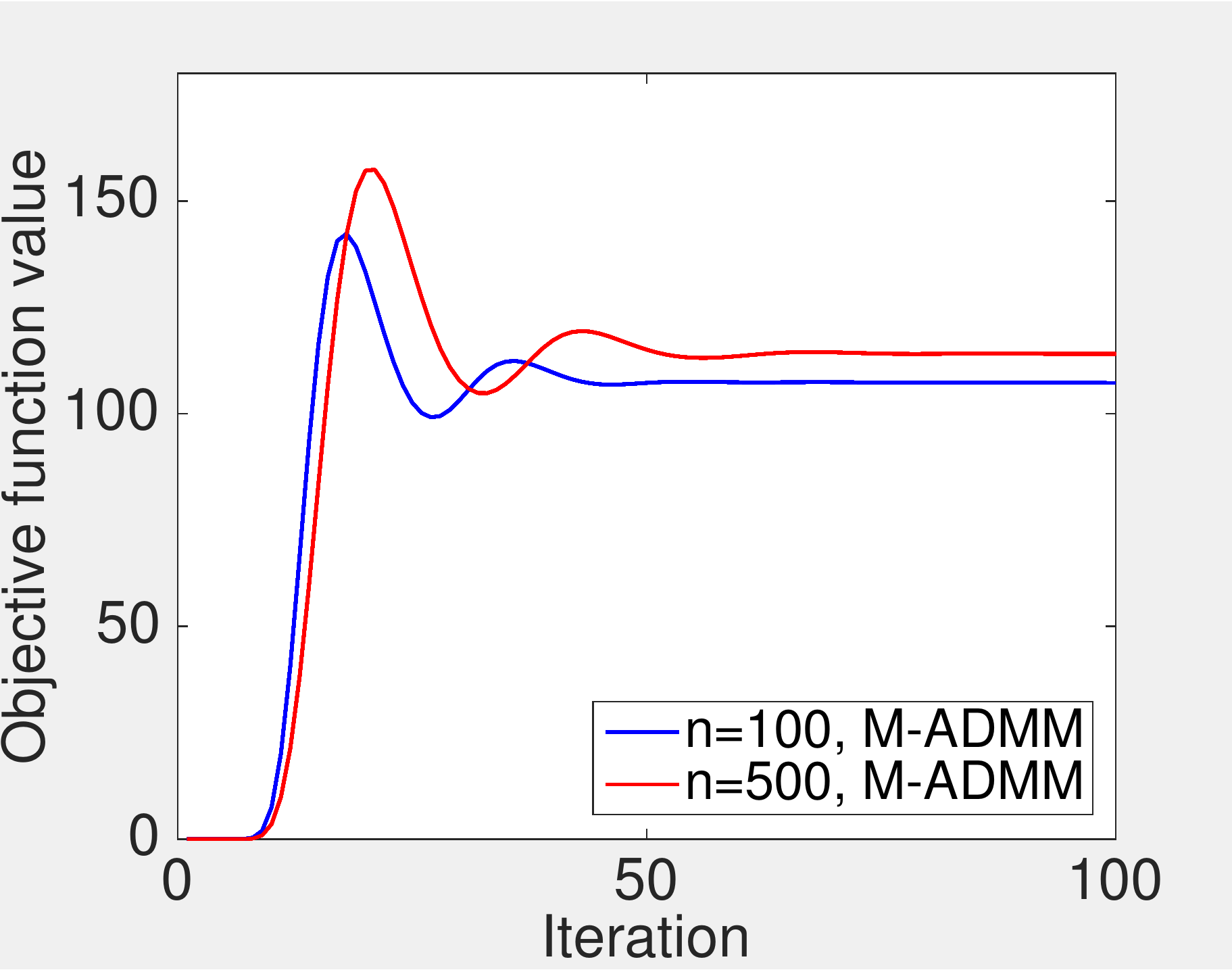} 
					\end{subfigure} 
					\begin{subfigure}[b]{0.155\textwidth}
						\centering
						\includegraphics[width=\textwidth]{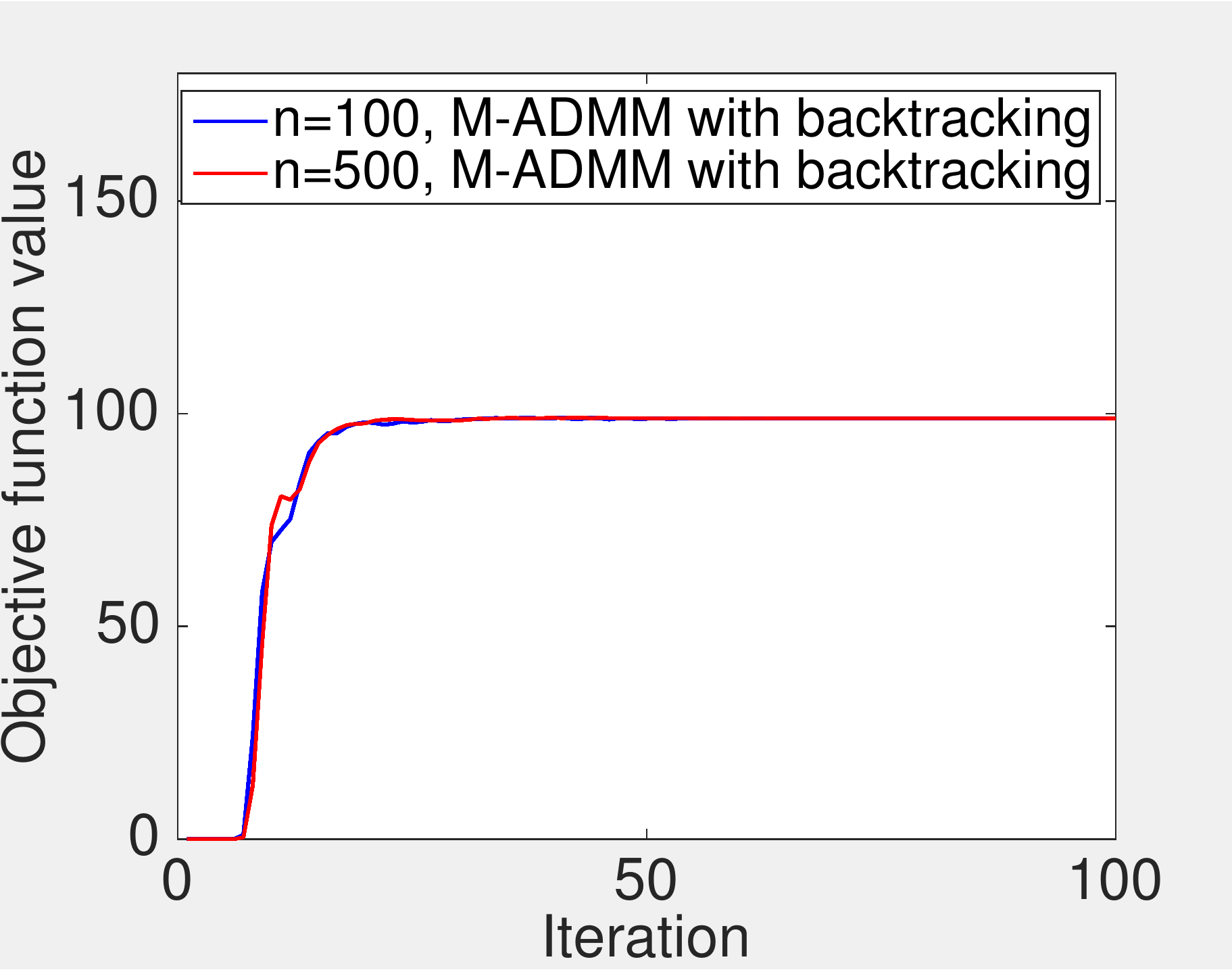} 
					\end{subfigure} 
					\caption{\small Top row: comparison of L-ADMM-PS, M-ADMM and M-ADMM with backtracking based on  $f(\xk)$  v.s. $k$ (left) and $\norm{\Axk-\b}$ v.s. $k$ (right) in  the case $n=500$. Bottom row: comparison of L-ADMM-PS (left), M-ADMM (middel), and M-ADMM with backtracking (right) in different cases of $n=100$ and $n=500$.  }\label{fig_backtrackinganalysis}
				\end{figure}


		\subsubsection{Analysis of M-ADMM with backtracking}
		We conduct three experiments to show the advantage of M-ADMM with backtracking, which uses tighter majorant surrogate, over L-ADMM-PS and M-ADMM. We still consider (\ref{nonsparse11}) on synthetic data. We generate $\A\in\mathbb{R}^{d\times m}$, where $d=50$ and $m=10,000$, with its elements independently sampled from an $N(0, 1)$ distribution. We generate $\x$ with $90\%$ elements being zeros  and others independently sampled from an $N(0, 1)$ distribution, and $\y=\A\x$. Then we uniformly split $\A$ and $\x$ into $n$ blocks, $[\A_1,\cdots,\A_n]$ and $\x=[\x_1;\cdots;\x_n]$, respectively. We   consider two cases:   $n=100$ and $n=500$. Though $n$ is different, the solved problems are equivalent. We are interested in the different convergence behaviors of the used solvers in both cases.  In M-ADMM with backtracking, we set $\tau=1.3$,  $\eta_i=0.01n_1\norm{\A_i}_2^2$, $i\in B_1$ and $\eta_i=0.01n_2\norm{\A_i}_2^2$, $i\in B_2$. The other settings and the initializations are the same as M-ADMM in Section \ref{sec_varpartstra}. Note that though the backtracking in Algorithm \ref{alg4}  requires some additional cost to estimate $\eta_i$'s, the cost can be ignored since the conditions in (\ref{eqbteq2}) and (\ref{eqbteq3}) fail only in a few iterations. Considering that the per-iteration complexity of the three solvers are the same, we simply compare their performance  based on  $f(\x^k)$ v.s. $k$  and $\norm{\A\x^k-\b}$ v.s.  $k$.

\begin{figure}
	\begin{subfigure}[b]{0.16\textwidth}
		\centering
		\includegraphics[width=\textwidth]{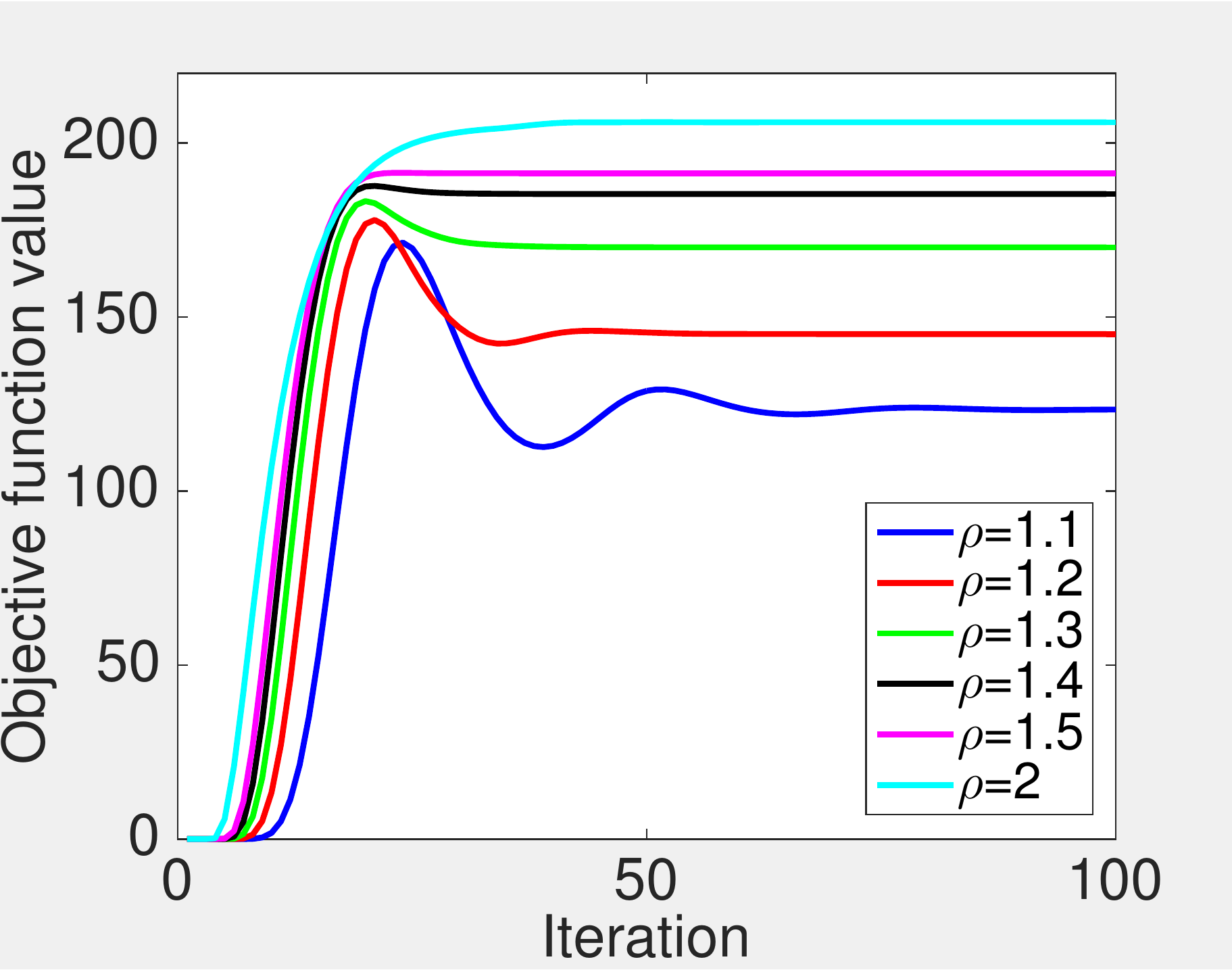} 
	\end{subfigure}
	\begin{subfigure}[b]{0.16\textwidth}
		\centering
		\includegraphics[width=\textwidth]{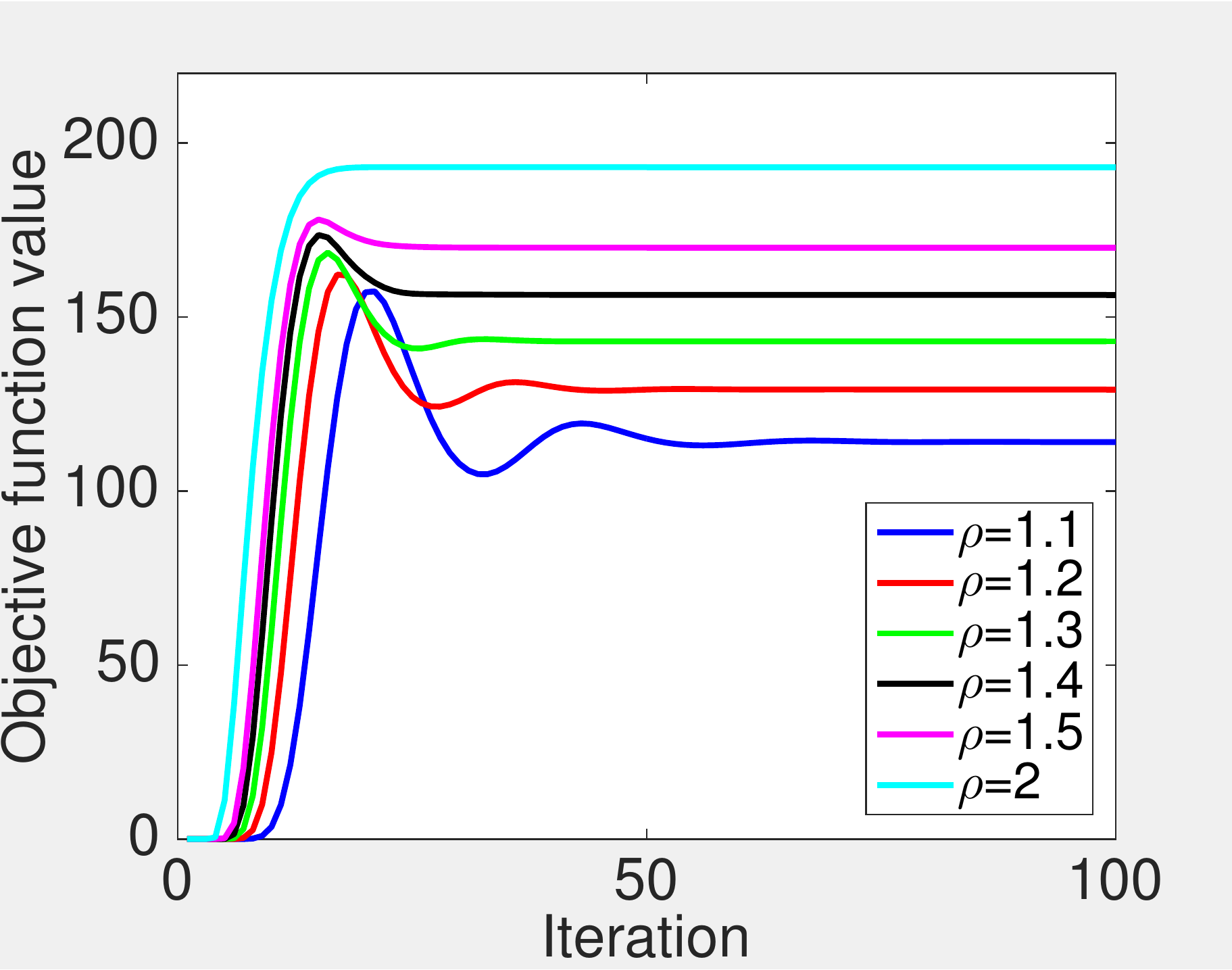} 
	\end{subfigure} 
	\begin{subfigure}[b]{0.16\textwidth}
		\centering
		\includegraphics[width=\textwidth]{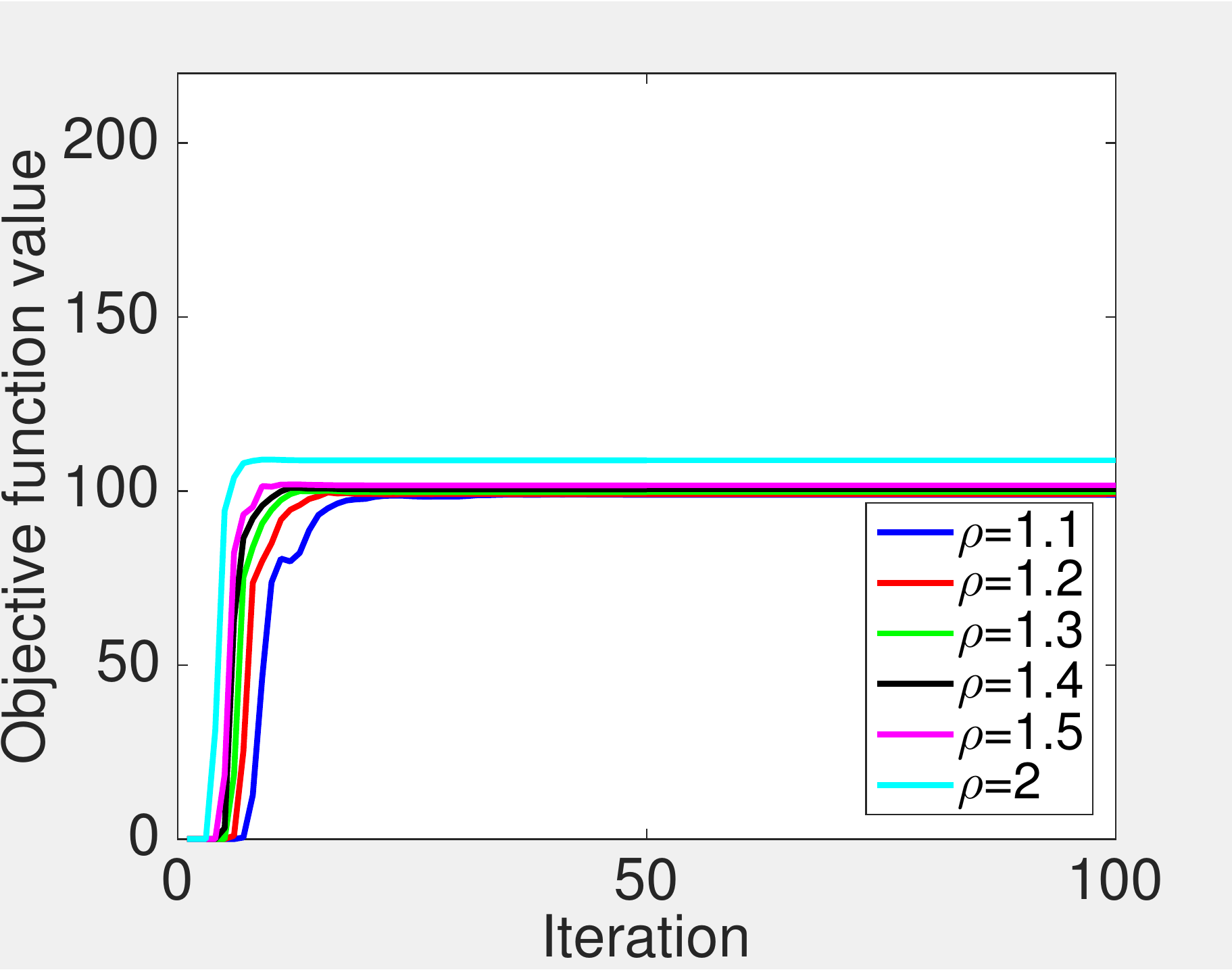} 
	\end{subfigure} 	
	
	\begin{subfigure}[b]{0.16\textwidth}
		\centering
		\includegraphics[width=\textwidth]{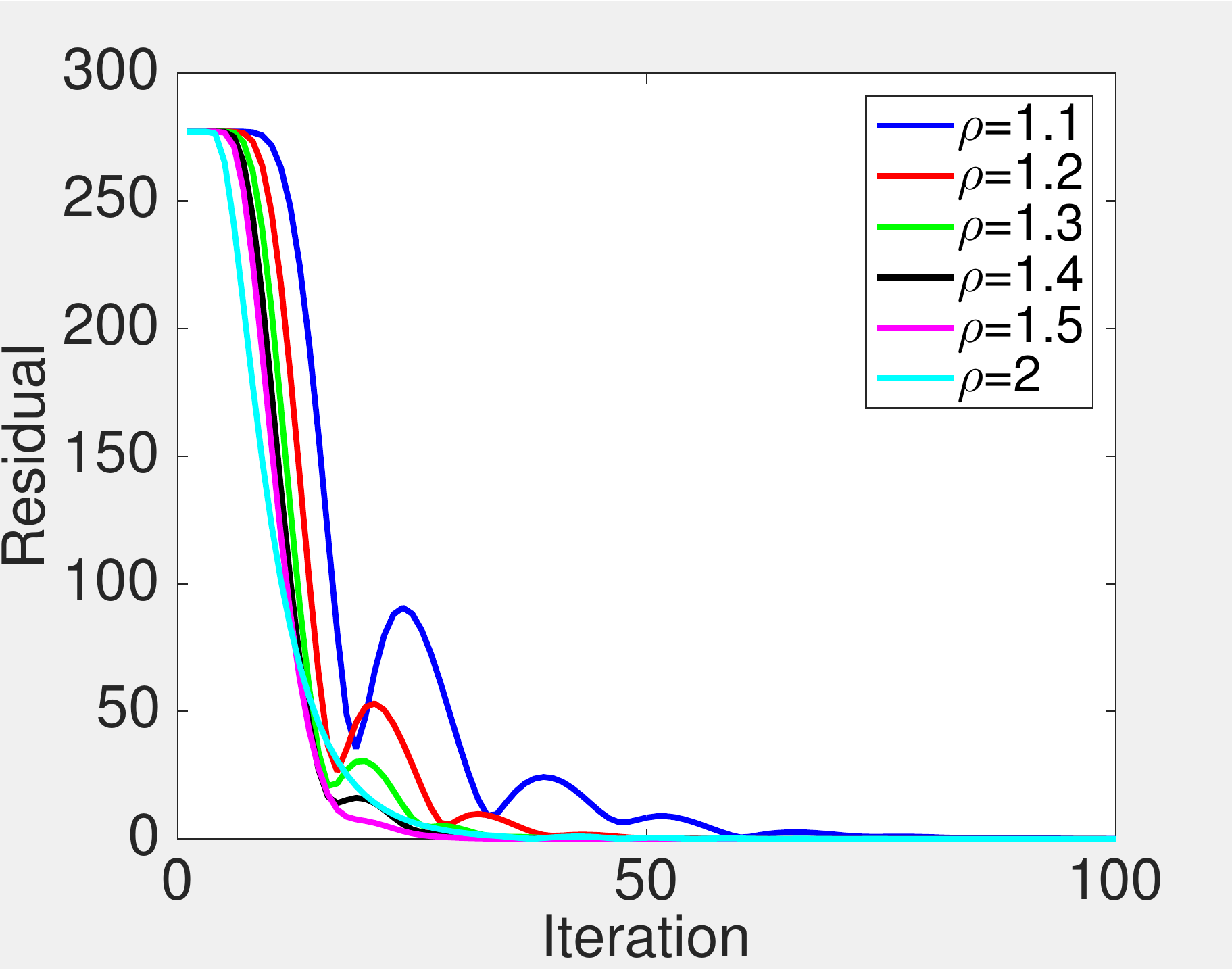} 
		\caption{}
	\end{subfigure} 
	\begin{subfigure}[b]{0.16\textwidth}
		\centering
		\includegraphics[width=\textwidth]{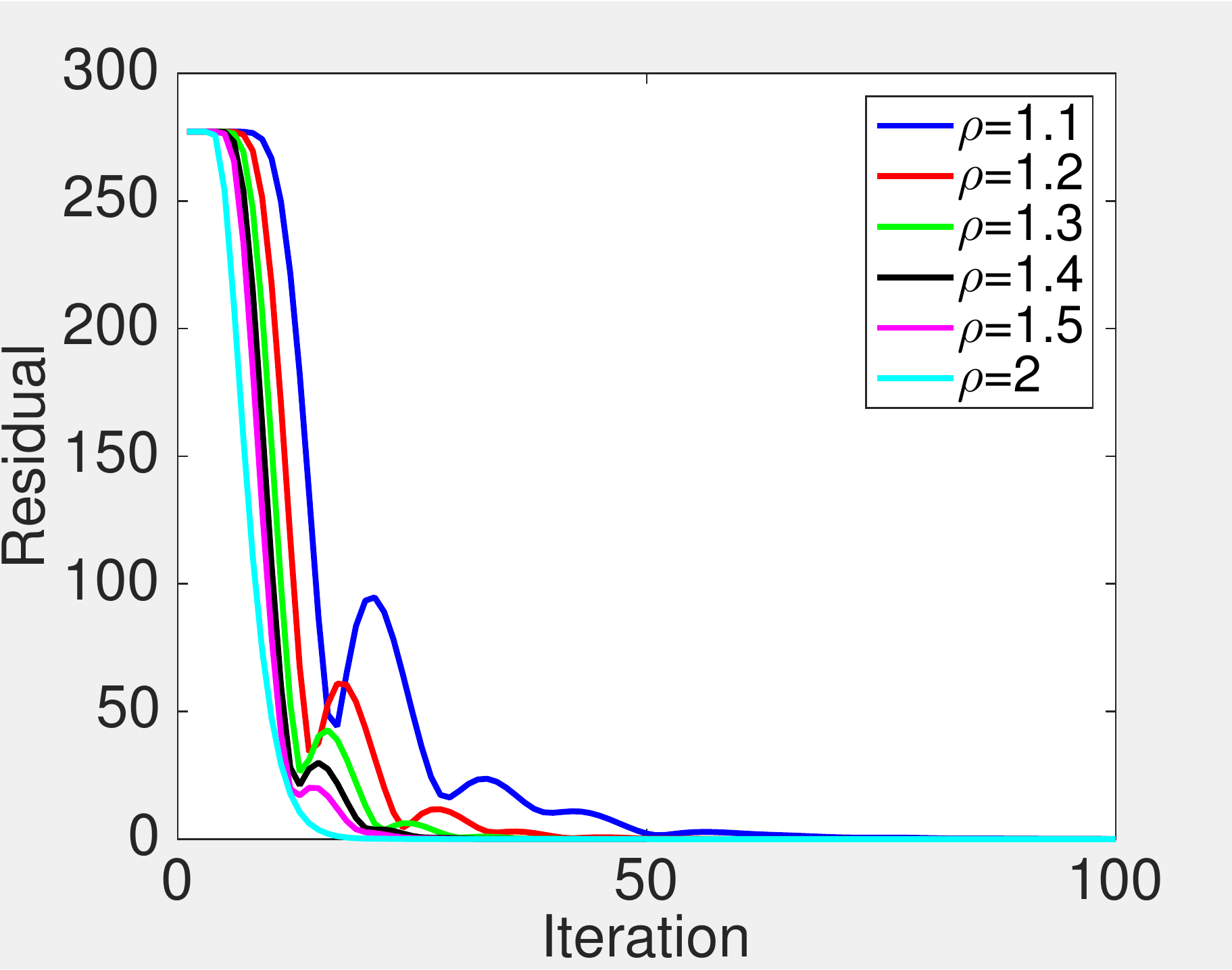} 
		\caption{}
	\end{subfigure} 
	\begin{subfigure}[b]{0.16\textwidth}
		\centering
		\includegraphics[width=\textwidth]{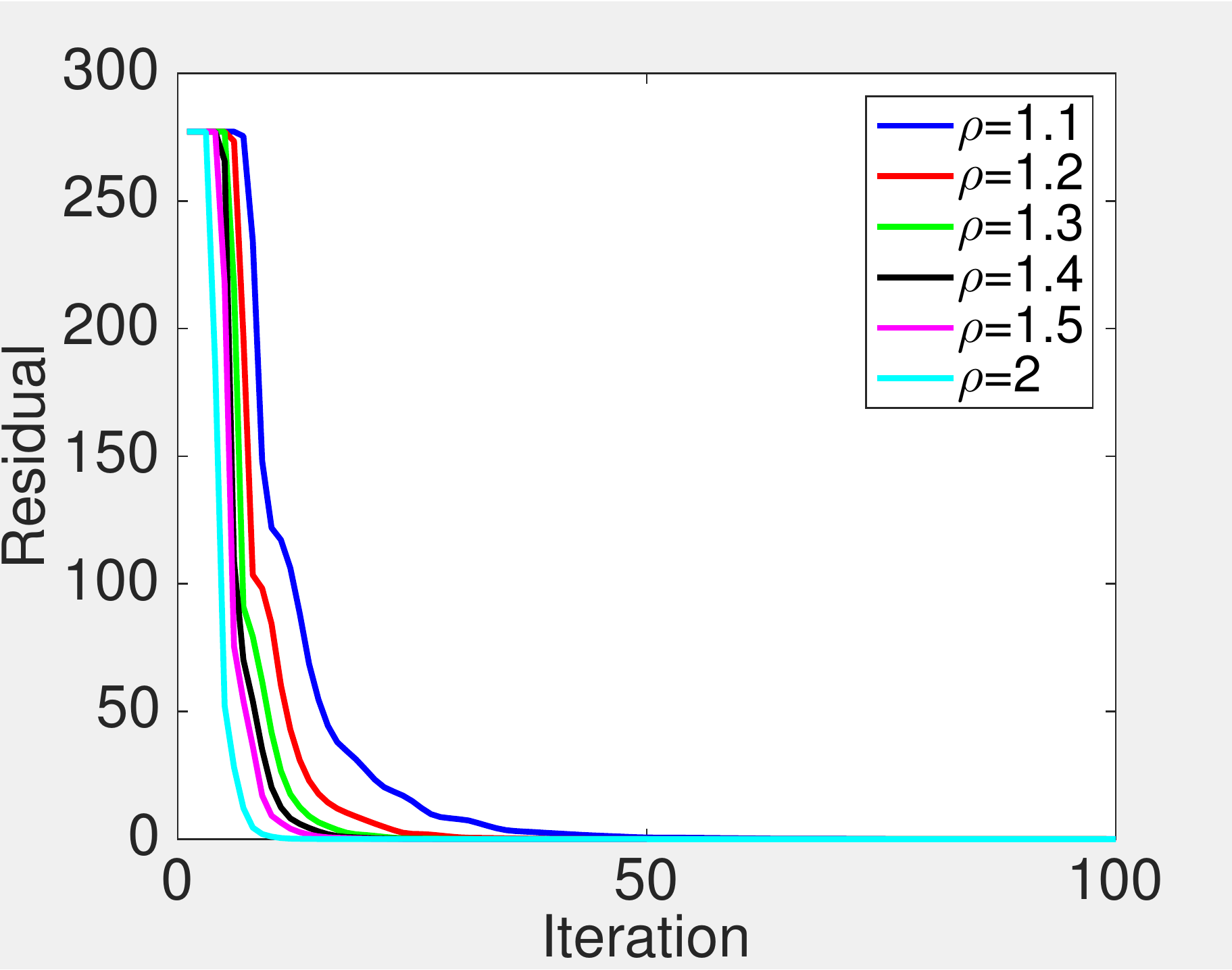} 
		\caption{}
	\end{subfigure} 
	\caption{\small Comparison of (a) L-ADMM-PS, (b) M-ADMM and (c) M-ADMM with backtracking on different choices of $\rho$, where $\betakk=\rho\betak$. Top row: plots of  $f(\xk)$  v.s.    $k$; bottom row: plots of $\norm{\Axk-\b}$ v.s. $k$.    }\label{fig_nl1_beta}
\end{figure}

		Figure \ref{fig_backtrackinganalysis} shows the comparison results. In Figure \ref{fig_backtrackinganalysis} (a)-(b), we consider the case $n=500$ and compare the three solvers based on $f(\x^k)$ v.s. $k$ ($k\leq 100$) and  $\norm{\A\x^k-\b}$ v.s. $k$. It can be seen that M-ADMM with backtracking achieves the smallest objective function value when the algorithm converges and it reduces the residual much faster than the other two methods. M-ADMM also outperforms L-ADMM-PS. These results well verifies the effectiveness of M-ADMM and the proposed backtracking technique, and are consistent with our theoretical analysis. Second,  we compare the convergence behaviors of the three solvers based on different splits of $\A$, i.e., $n=100$ and $n=500$. 
		From Figure \ref{fig_backtrackinganalysis} (c)-(d), it can be seen  that L-ADMM-PS and M-ADMM for the case $n=100$ perform much better than the case $n=500$, respectively. This is not a surprise since both two solvers use constant $\eta_i$'s which depend on the block number (see Theorem \ref{them3}). Intuitively, the smaller $n$ leads to a tighter majorant surrogate, e.g., (\ref{rone}), and thus it leads to a better approximated solution.
		However, M-ADMM with backtracking performs the best and it is not sensitive to   $n$, since it estimates $\eta_i$'s locally and this leads to a tight majorant surrogate.
		
		Furthermore, we compare the three solvers based on different choices of  $\rho\in\{1.1,1.2,1.3,1.4,1.5,2\}$, where $\betakk=\rho\betak$. We test on the same dataset as the above experiment with $n=500$, and plot $f(\xk)$ v.s. $k$ and $\norm{\A\xk-\b}$ v.s. $k$ in Figure \ref{fig_nl1_beta}. For all the three solvers, when $\rho$ is larger, the residual $\norm{\A\xk-\b}$ decreases faster. More importantly, the price is that the objective $f(\xk)$ decreases slower. Considering  the convergence of $f(\xk)$, both L-ADMM-PS and M-ADMM are sensitive to the choice of $\rho$, though the later one performs better. However,  M-ADMM with backtracking  performs very well even when $\rho$ increases. The reason is that the larger $\rho$ implies that $\betak$ increases much faster and this  makes the majorant surrogates in (\ref{rone})-(\ref{rtwo}) much looser. In contrast, the surrogates $\hatr^k_{B_1}$ and $\hatr^k_{B_2}$  in M-ADMM with backtracking are computed  locally based on (\ref{eqbteq2})  and (\ref{eqbteq3})   and thus the surrogates are much tighter.   
		This experiment verifies that the backtracking technique allows a relatively faster increasing sequence  $\{\betak\}$ and improves the convergence. 
		
\begin{figure}
	\begin{subfigure}[b]{0.16\textwidth}
		\centering
		\includegraphics[width=\textwidth]{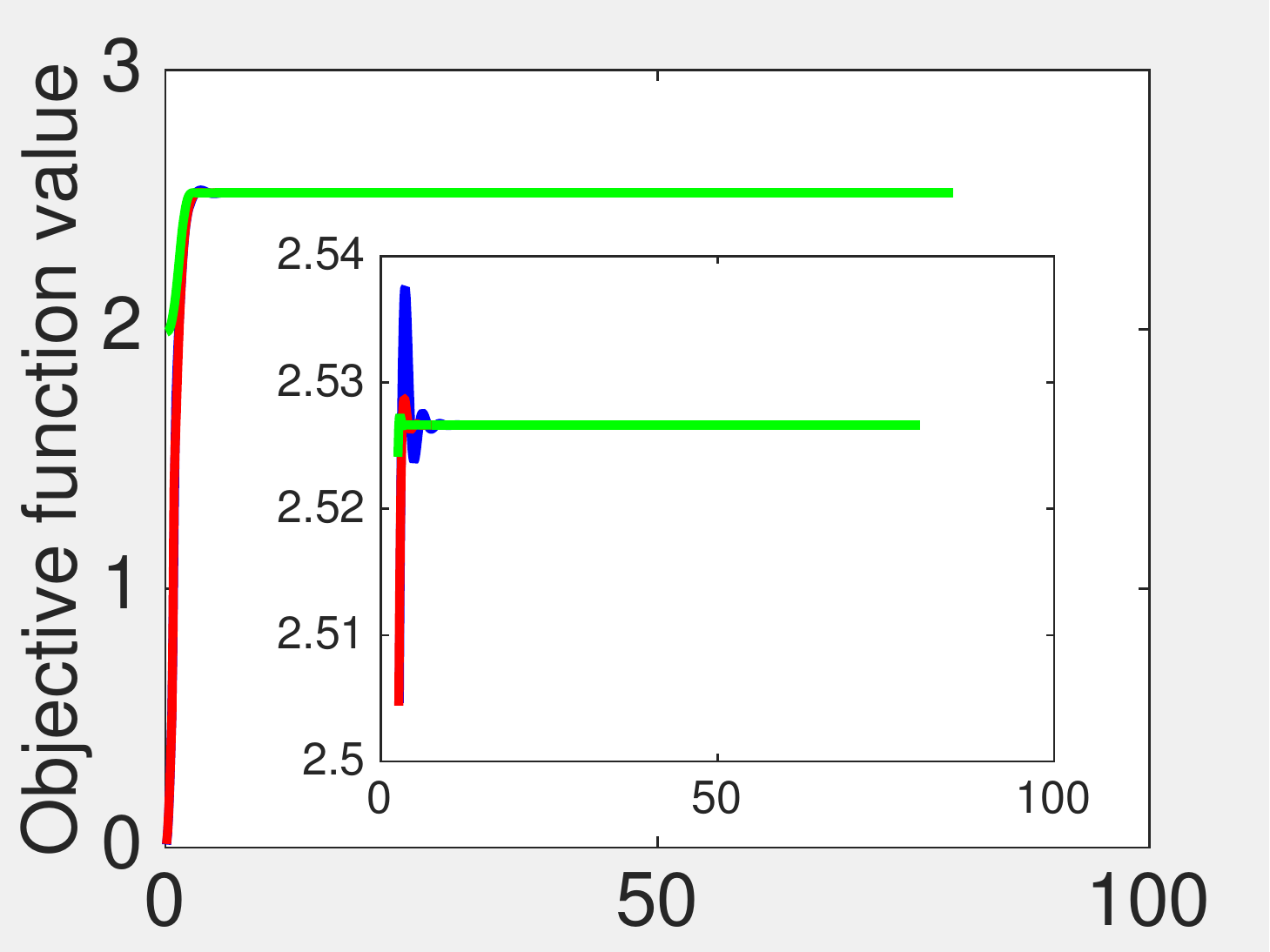} 
	\end{subfigure}
	\begin{subfigure}[b]{0.16\textwidth}
		\centering
		\includegraphics[width=\textwidth]{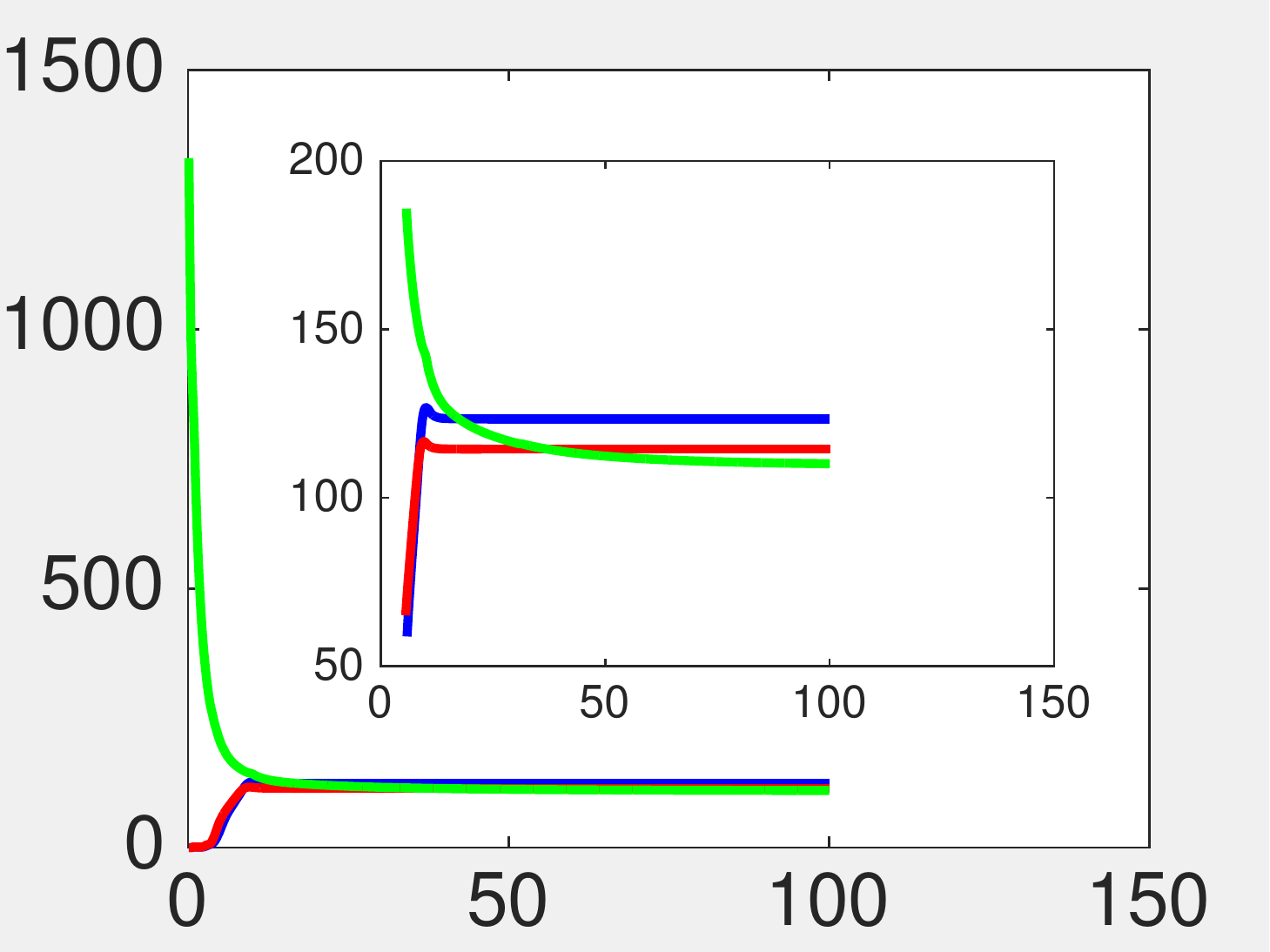} 
	\end{subfigure} 
	\begin{subfigure}[b]{0.16\textwidth}
		\centering
		\includegraphics[width=\textwidth]{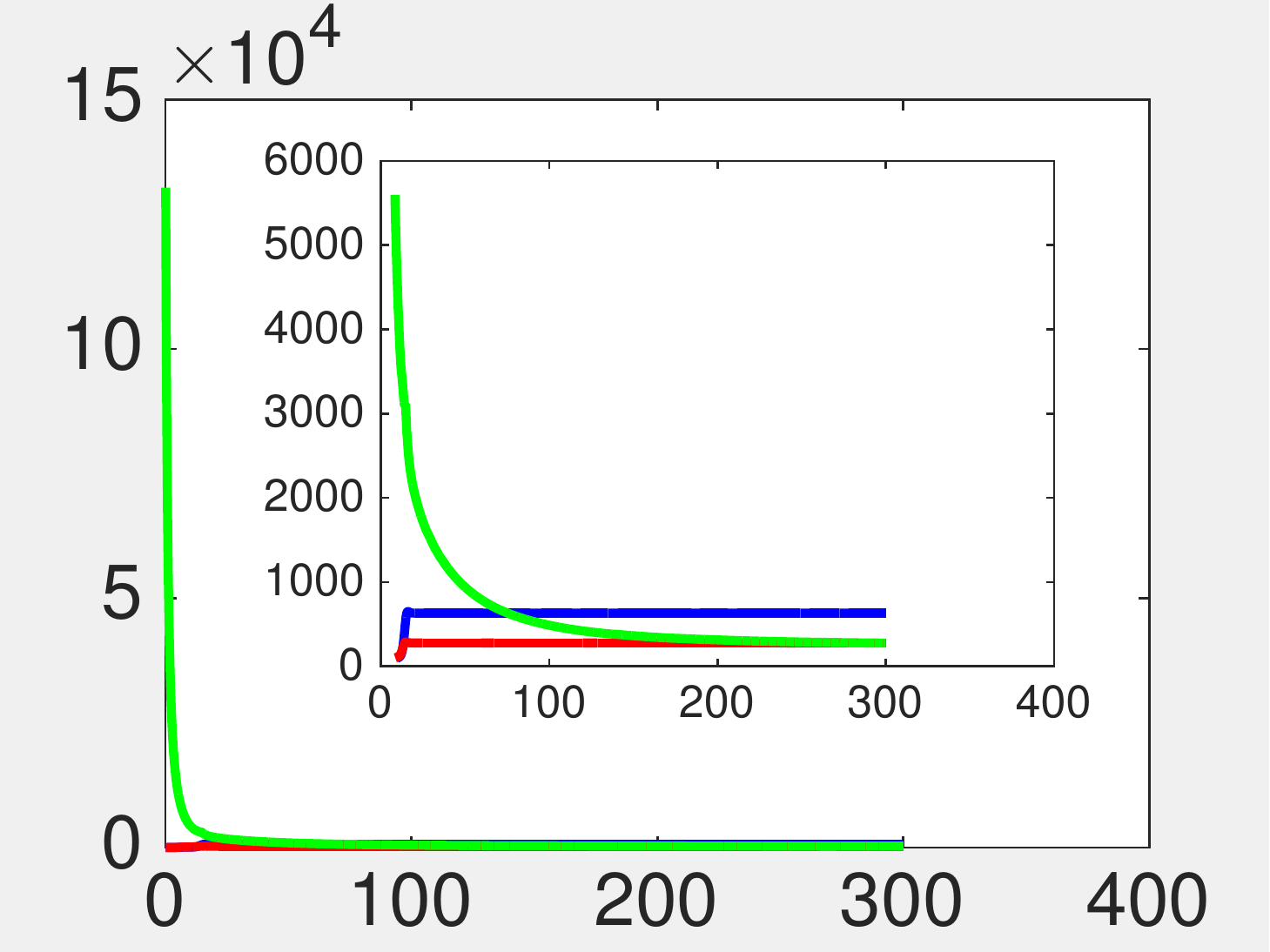} 
	\end{subfigure}

	\begin{subfigure}[b]{0.16\textwidth}
		\centering
		\includegraphics[width=\textwidth]{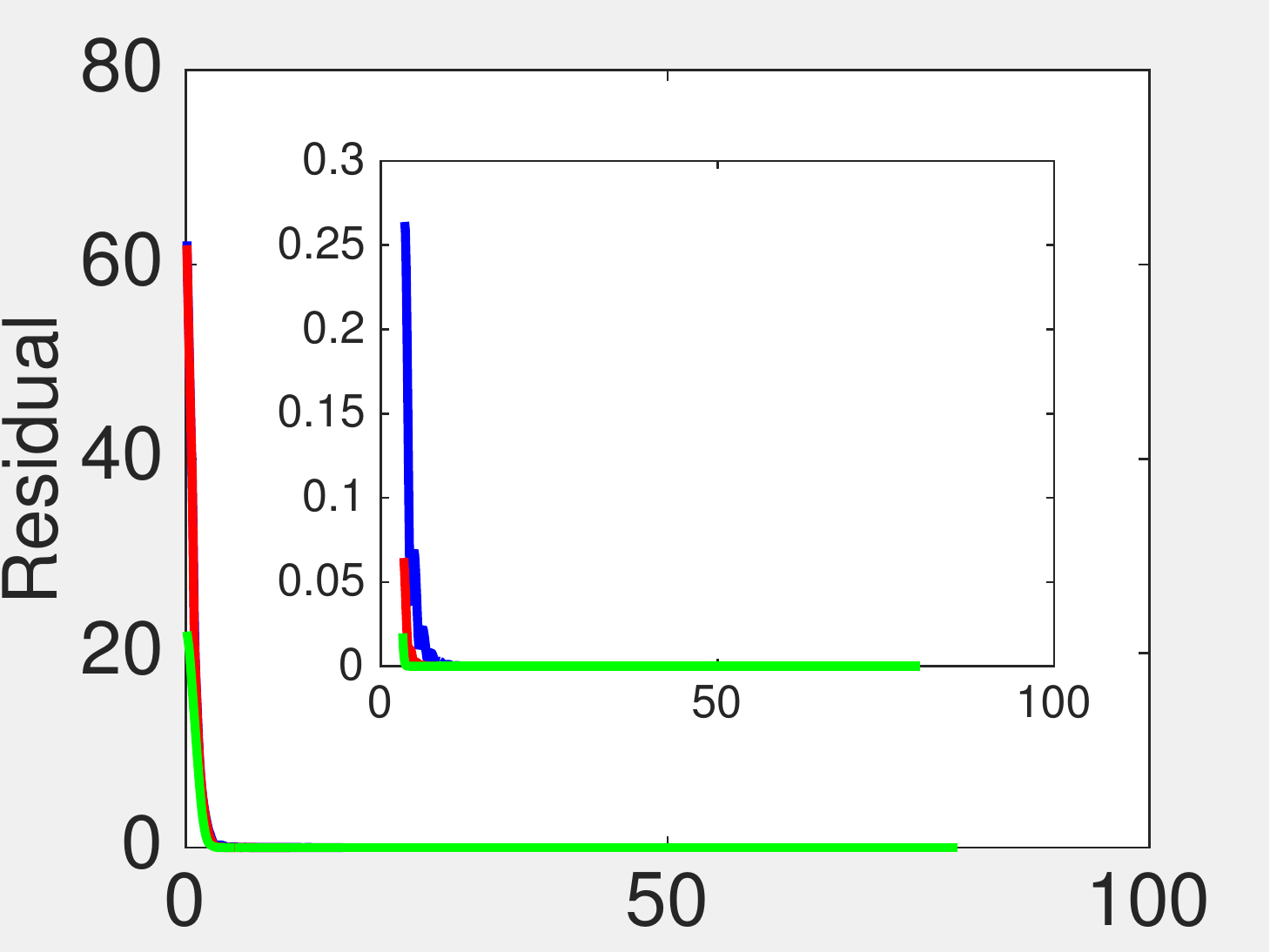} 
		\caption{$\lambda=0.001$}
	\end{subfigure} 
	\begin{subfigure}[b]{0.16\textwidth}
		\centering
		\includegraphics[width=\textwidth]{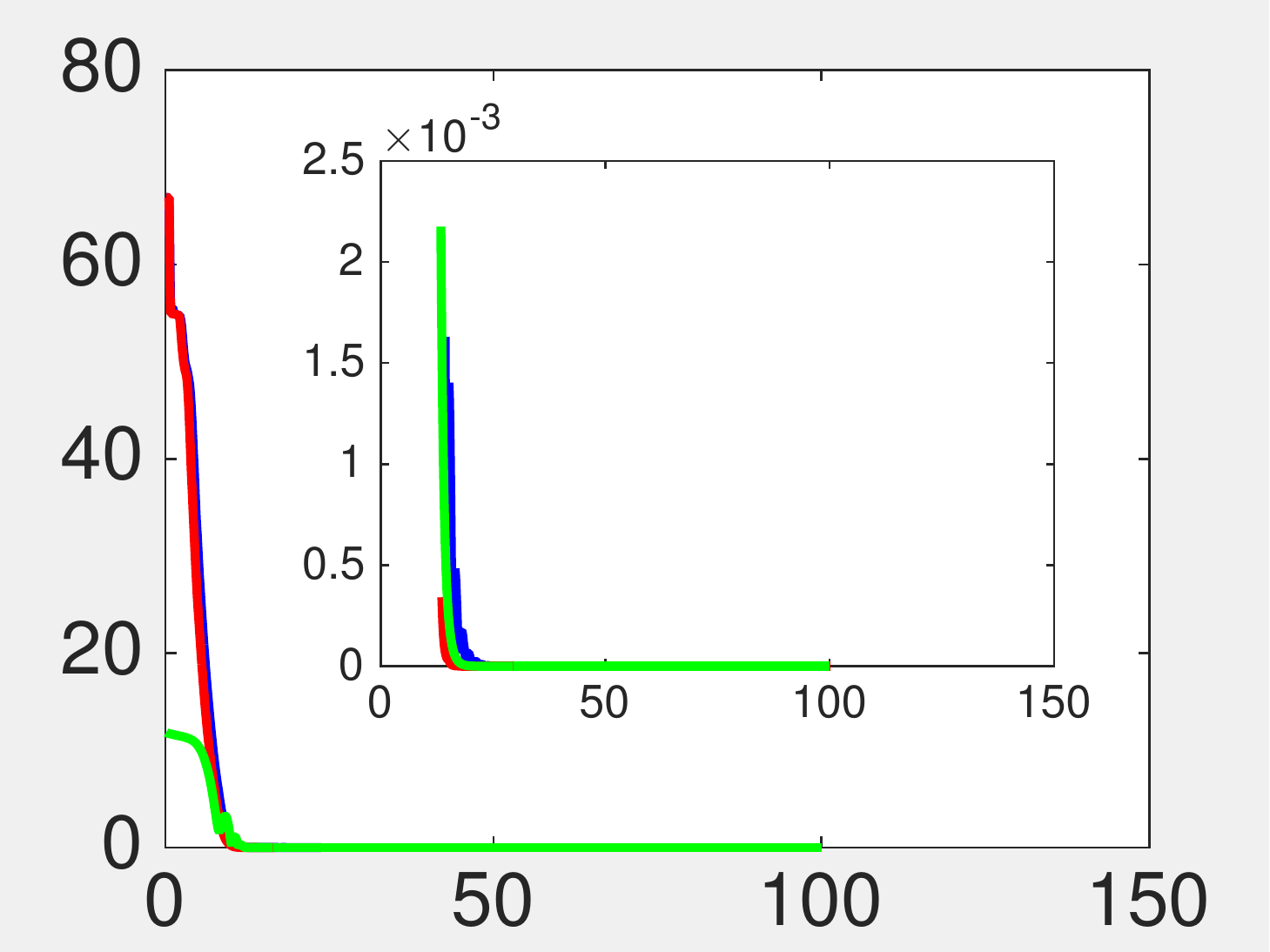} 
		\caption{$\lambda=0.1$}
	\end{subfigure} 
	\begin{subfigure}[b]{0.16\textwidth}
		\centering
		\includegraphics[width=\textwidth]{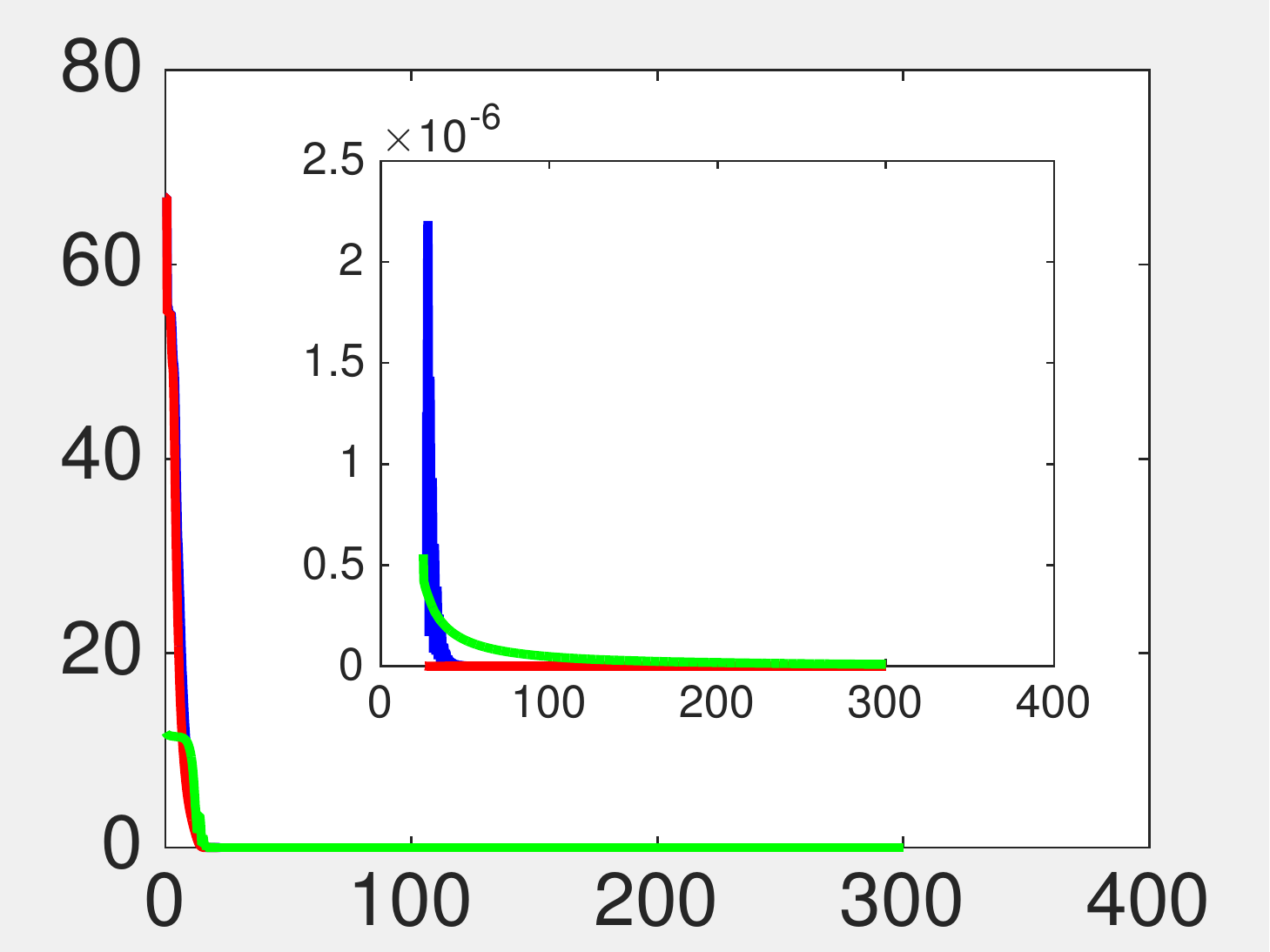} 
		\caption{$\lambda=10$}
	\end{subfigure} 
	\begin{subfigure}[b]{0.16\textwidth}
		\centering
		\includegraphics[width=\textwidth]{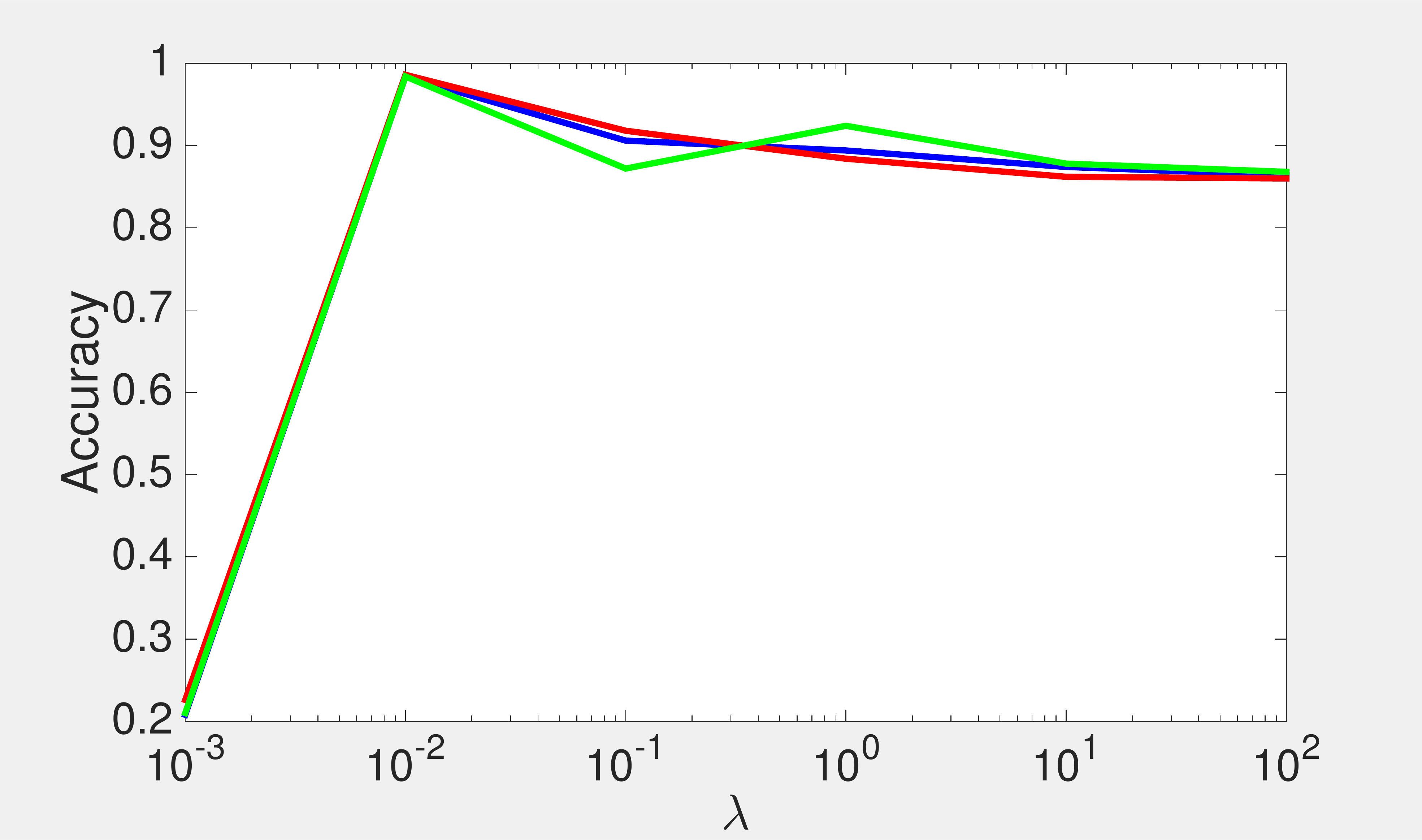} 
		\caption{}
	\end{subfigure}
	\begin{subfigure}[b]{0.16\textwidth}
		\centering
		\includegraphics[width=\textwidth]{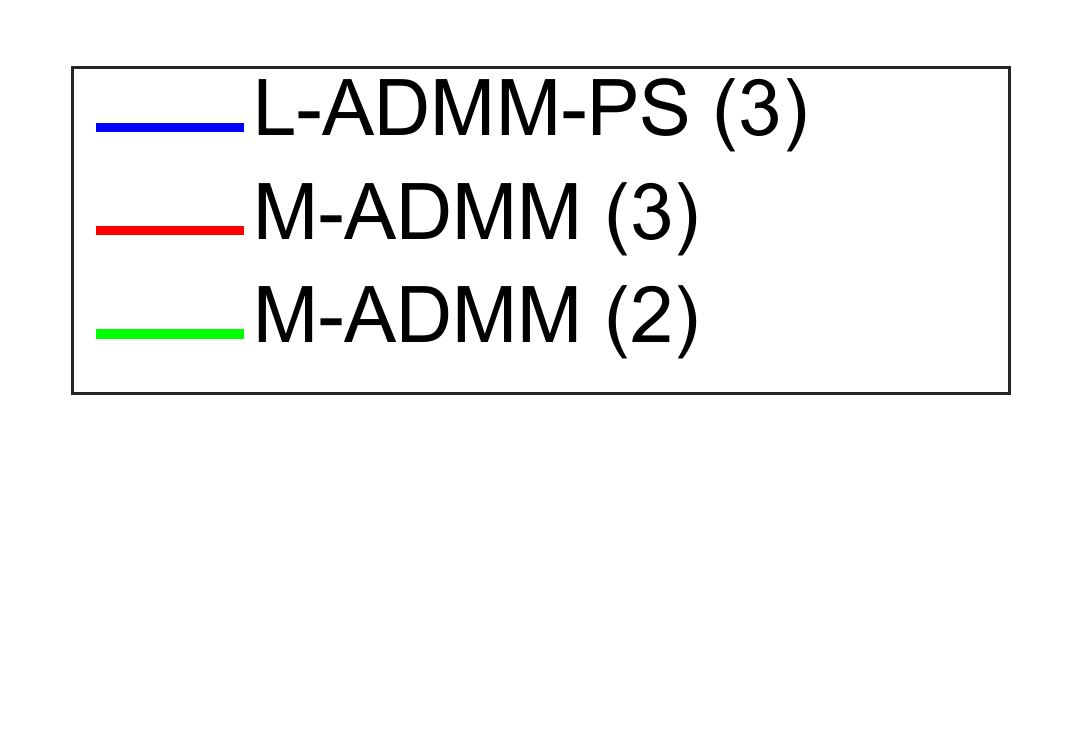} 
	\end{subfigure} 
	\caption{\small Comparison of L-ADMM-PS (3), M-ADMM (3) and M-ADMM (2) on different choices of $\lambda$: (a) $\lambda=0.001$; (b) $\lambda=0.1$ and (c) $\lambda=10$. Top row: plots of  $f(\xk)$  v.s.   CPU time; middle row: plots of $\norm{\Axk-\b}$ v.s. CPU  time.   (d) Subspace clustering accuracy v.s. $\lambda$. In (a)-(c), for better visualization, we plot the objective value and residual within a relatively smaller range of CPU time in sub-figures.  }\label{fig_latlrrtoy}
\end{figure}

\subsection{Solving Non-separable Objective Problem}		
		To show that M-ADMM can solve the problem with non-separable objective, we consider the Latent Low-Rank Representation (LatLRR) problem \cite{liu2011latent} for affine subspace clustering
		\begin{equation}\label{latlrr}
		\min_{\Z,\bL} \|\Z\|_*+\|\bL\|_*+\frac{\lambda}{2}\|\X\Z+\bL\X-\X\|_F^2, \ \text{s.t.} \ \bm{1}^\top\Z=\bm{1}^\top,
		\end{equation}
		where $\lambda>0$ and the constraint is due to the affine subspace structure of data $\X$ \cite{6482137}. The objective of (\ref{latlrr}) is non-separable and can  be rewritten as the following one with separable objective 
		\begin{equation}\label{latlrr2}
		\begin{split}
		\min_{\Z,\bL,\E} \|\Z\|_*+\|\bL\|_*+\frac{\lambda}{2}\|\E\|_F^2, \\
		\ \text{s.t.}  \ \bm{1}^\top\Z=\bm{1}^\top,\ \X\Z+\bL\X-\X=\E.
		\end{split}	
		\end{equation}
		We compare the following three solvers which own the convergence guarantee to solve  the latent LRR problem:
		\begin{itemize}
			\item {L-ADMM-PS (3)}:  use  (\ref{ladmmpsxi}) for   3 blocks problem (\ref{latlrr2}).
			\item {M-ADMM (3)}: use M-ADMM for   3  blocks problem (\ref{latlrr2}).
			\item {M-ADMM (2)}: use M-ADMM for 2 blocks problem (\ref{latlrr}). 
		\end{itemize}
		Note that $h(\Z,\bL)=\frac{1}{2}\|\X\Z+\bL\X-\X\|_F^2$ in (\ref{latlrr}) is $\{2\|\X\|_2^2\I,2\|\X\|_2^2\I\}$-smooth. M-ADMM (2) uses the Lipschitz gradient surrogate in  (\ref{sqrsurr1}) to make the subproblems   separable. For M-ADMM (3), we partition the three variables into two super blocks: $\{\Z\}$ and $\{ \bL,\E\}$, and update them in the Gauss-Seidel way. In contrast, L-ADMM-PS updates $\Z$, $\bL$ and $\E$ in parallel.  
			
		We apply latent LRR for subspace clustering by using the learned $\Z$ based on both the synthesized and real data. For the synthesized data, we  generate $\X=[\X_1,\X_2,\cdots]$ with its columns sampled from different subspaces. We construct $k=5$ independent subspaces $\{\mathcal{S}_i\}_{i=1}^5\subseteq \mathbb{R}^{200}$ whose bases $\{\U_i\}_{i=1}^5$ are computed by $\U_i=\T\U_i$, $1\leq i \leq 4$, where $\T$ is a random rotation and $\U_1\in\mathbb{R}^{200\times 5}$ is a random orthogonal matrix. We sample 100 vectors from each subspace by $\X_i = \U_i\Q+0.1$, $1\leq i\leq 5$ with $\Q\in\mathbb{R}^{5\times 100}$ being an  i.i.d. $N(0,1)$ matrix. Furthermore, $20\%$ of data vectors are chosen to be corrupted, e.g., for a data vector $\x$ chosen to be corrupted, its observed vector is computed by adding Gaussian noise with zero mean and variance $0.2\norm{\x}$. Given $\X\in\mathbb{R}^{200\times 500}$ by the above way, we can solve the latent LRR problem by the three solvers and obtain the solution $\Z^*$. Then the data vectors can be grouped into $k$ groups based   on the affinity matrix $(|\Z^*|+|(\Z^*)^\top|)/2$ by spectral clustering \cite{liu2011latent}. The clustering accuracy is used to evaluate the clustering performance \cite{liu2011latent}.  We test on different choices of $\lambda$ and compare the three solvers based on   $f(\x^k)$ v.s. CPU time (in seconds),   $\norm{\A\x-\b}$ v.s. CPU time and clustering accuracy. The results are  shown in Figure \ref{fig_latlrrtoy} and we have the following observations:
		\begin{itemize}
			\item M-ADMM (3) always outperforms    L-ADMM-PS (3) in the sense that the objective value is smaller when the algorithms converge and the residual decreases much faster. Both  solve the same problem (\ref{latlrr2}) with 3 blocks of variables. But M-ADMM (3) updates ${\Z}$ and $\{\bL,\E\}$ sequentially, and thus it is faster than L-ADMM-PS (3) which updates them in parallel.   This   is consistent with our analysis at the end of Section \ref{subsectionmadmm}. 
			\item When $\lambda$ is relatively small, M-ADMM (2)  converges faster than M-ADMM (3). When $\lambda$ is relatively large, M-ADMM (2) leads to a smaller objective value, but it requires much more running time (many more iterations). Both solvers have their advantages and disadvantages. In this experiment, the block number $n$ and the looseness of the surrogate are two crucial factors.
			M-ADMM (2) solves (\ref{latlrr}) with only 2 blocks, but it requires   constructing the Lipschitz gradient surrogate by (\ref{sqrsurr1}) for $h(\Z,\bL)=\frac{\lambda}{2}\|\X\Z+\bL\X-\X\|_F^2$.  This surrogate is looser when $\lambda$ is lager. This is why M-ADMM (2) is slower when $\lambda$ increases (the same phenomenon also appears in ISTA and FISTA \cite{beck2009fast}). On the other hand,  M-ADMM (3) for 3 blocks problem (\ref{latlrr2}) converges quickly regardless of the choice of  $\lambda$. The issue of M-ADMM (3) is that the surrogate $\hatr_i^k(\x_i)$ in (\ref{rone})-(\ref{rtwo}) also becomes looser when $\betak$ increases. So M-ADMM (3) may quickly get stuck and the final objective value is larger than M-ADMM (2).
			In practice, one has to balance the effects of both the block number $n$ and the looseness of the surrogate, by considering the specific problems.  
		\end{itemize}

 \begin{table}[t]
 	\caption{\small Comparison of L-ADMM-PS (3) and M-ADMM (3) and M-ADMM (2)  for latent LRR on the Hopkins 155 dataset.}
 	\footnotesize
 	\centering
 	\begin{tabular}{c|c|c|c}
 		\hline
 		Methods	        & L-ADMM-PS (3)& M-ADMM (3)    & M-ADMM (2)    \\ \hline\hline
 		Accuracy ($\%$)             & 90.9     & \textbf{92.7}&  87.1 \\ \hline
 		CPU Time (s)  & 756.2    & \textbf{738.5}  & 932.1      \\ \hline
 	\end{tabular} \label{tabhopkins155}   
 \end{table} 
 
		We further apply latent LRR for motion segmentation and test on the Hopkins 155 dataset \cite{tron2007benchmark}. This dataset contains 156 sequences,  each with 39$\sim$550 vectors drawn from two or three motions (one motion corresponds to one subspace). Each sequence is a sole segmentation (clustering) task and thus there are 156 clustering tasks in total. We follow the experimental settings in \cite{liu2011latent} but without the complex post-processing. We set $\lambda=500$ and compare the performance by using M-ADMM (2), L-ADMM-PS (3) and M-ADMM (3). We  stop the algorithms when
		\begin{equation}\label{eqstopcri}
		\norm{\A\xk-\b}/\norm{\b} \leq \epsilon, \text{ and } \norm{\xkk-\xk}/\norm{\b}\leq \epsilon,
		\end{equation}
		where $\epsilon=10^{-4}$. For each motion sequence, we record the clustering accuracy and the CPU time of  solvers. Then  the mean clustering accuracy and the total CPU time of all 156 sequences  are reported in Table \ref{tabhopkins155}. It can be seen that, due to the same stopping criteria in (\ref{eqstopcri}), the CPU time of L-ADMM-PS (3) and that of M-ADMM (3) are similar. But the solution to latent LRR obtained by  M-ADMM (3) achieves better clustering accuracy than L-ADMM-PS (3). The reason is that M-ADMM (3) obtains a better solution with much smaller objective value within similar running time (or similar number of iterations). In this experiment, M-ADMM (2) for (\ref{latlrr}) is inferior to the  other two solvers since the used $\lambda$ is relatively large and thus the used  majorant surrogate is loose.


\begin{figure}
	\begin{subfigure}[b]{0.057\textwidth}
		\centering
		\includegraphics[width=\textwidth]{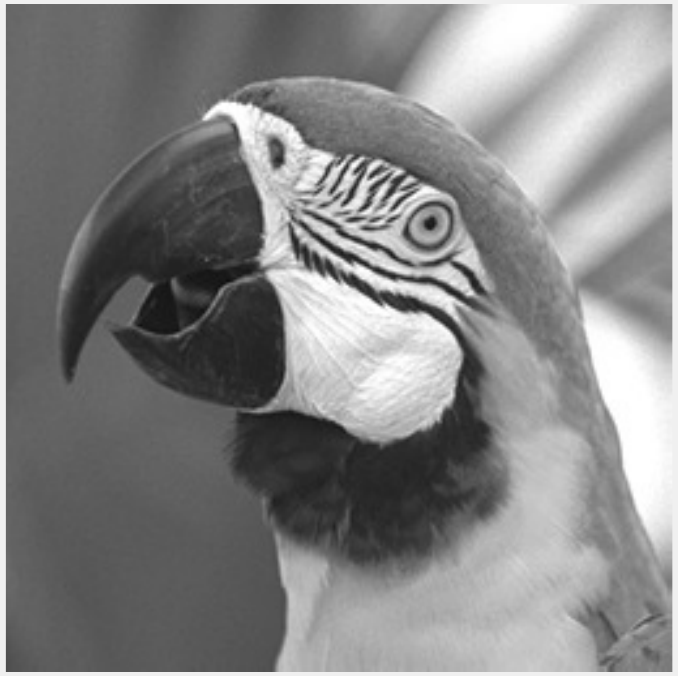} \vspace{-7mm}
		\caption*{\tiny parrot}
	\end{subfigure} 
	\begin{subfigure}[b]{0.057\textwidth}
		\centering
		\includegraphics[width=\textwidth]{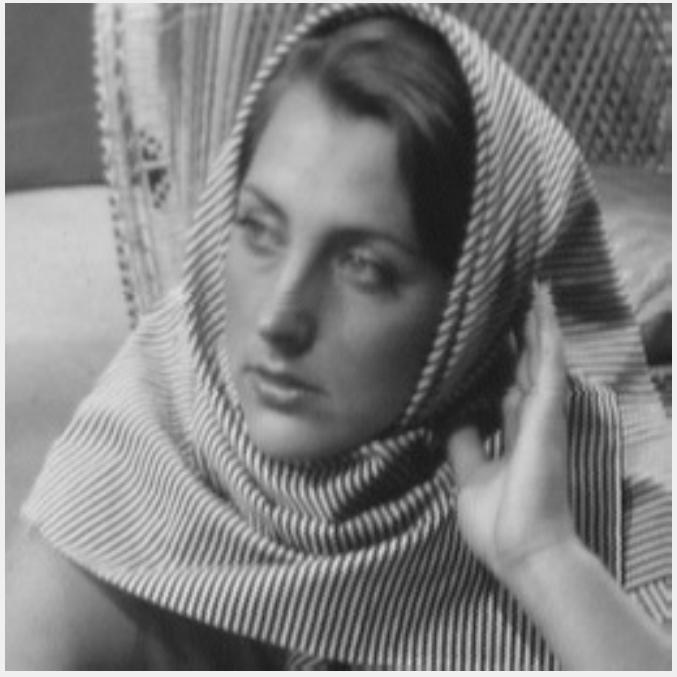} \vspace{-7mm}
		\caption*{\tiny barbara }
	\end{subfigure}
	\begin{subfigure}[b]{0.057\textwidth}
		\centering
		\includegraphics[width=\textwidth]{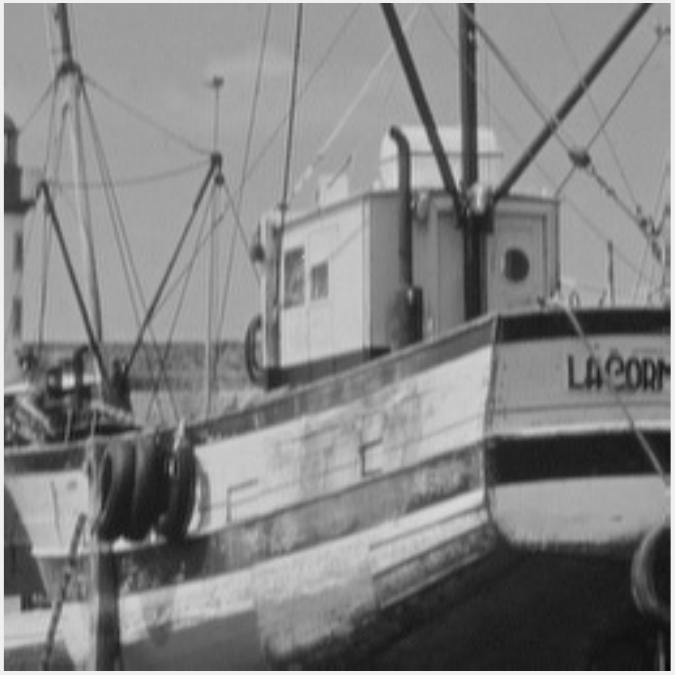} \vspace{-7mm}
		\caption*{\tiny boat }
	\end{subfigure} 
	\begin{subfigure}[b]{0.057\textwidth}
		\centering
		\includegraphics[width=\textwidth]{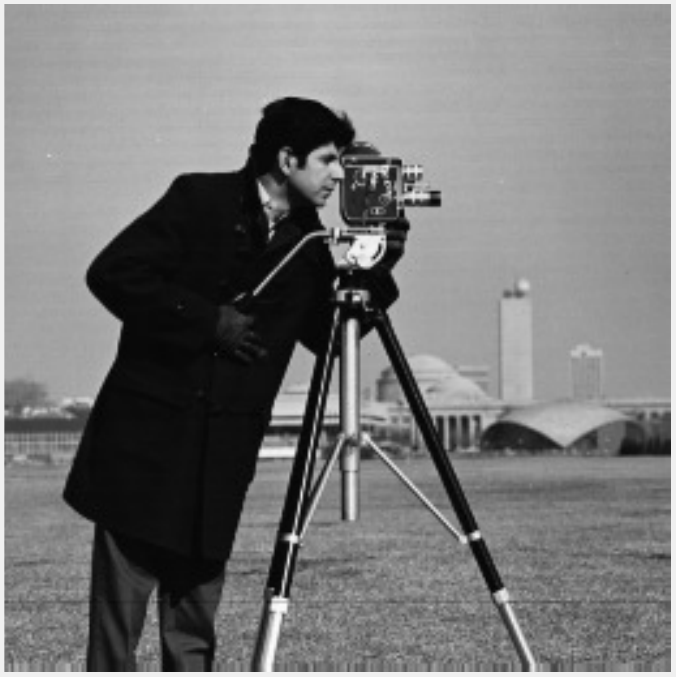} \vspace{-7mm}
		\caption*{\tiny  cameraman}
	\end{subfigure} 		
	\begin{subfigure}[b]{0.057\textwidth}
		\centering
		\includegraphics[width=\textwidth]{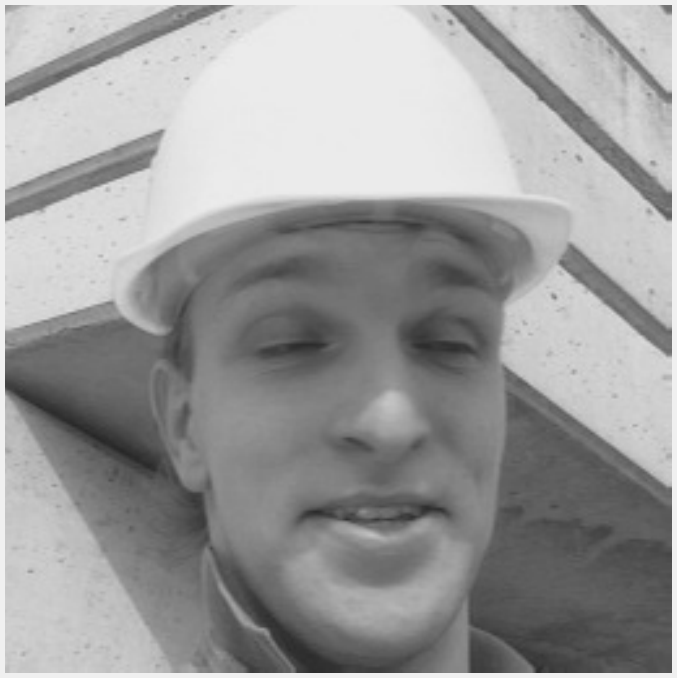} \vspace{-7mm}
		\caption*{\tiny foreman }
	\end{subfigure} 
	\begin{subfigure}[b]{0.057\textwidth}
		\centering
		\includegraphics[width=\textwidth]{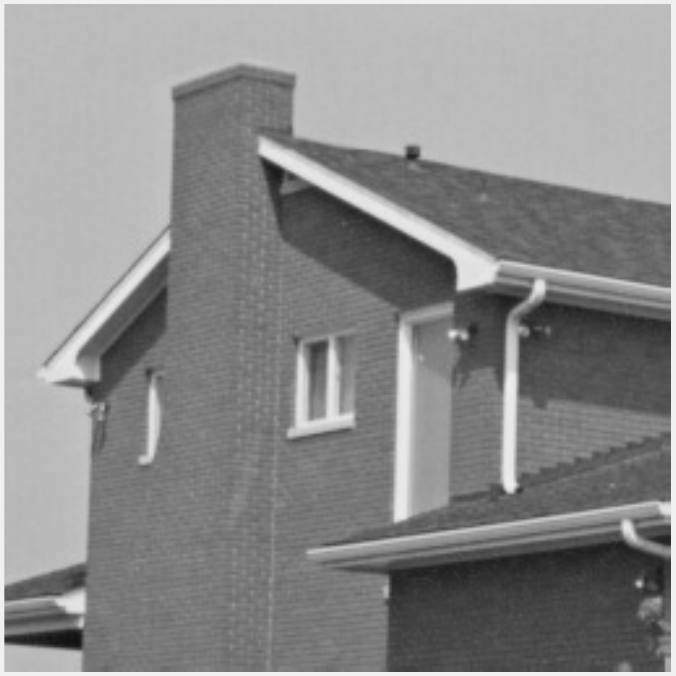} \vspace{-7mm}
		\caption*{\tiny house }
	\end{subfigure} 
	\begin{subfigure}[b]{0.057\textwidth}
		\centering
		\includegraphics[width=\textwidth]{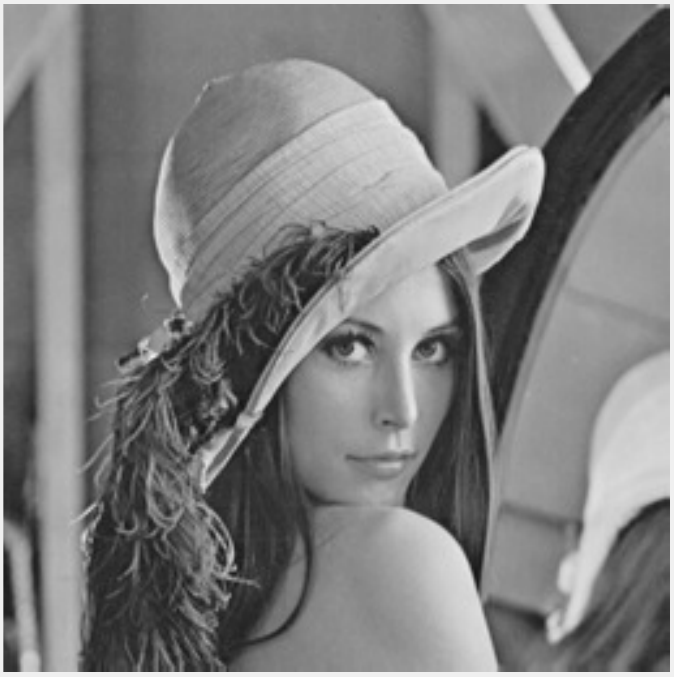} \vspace{-7mm}
		\caption*{\tiny lena}
	\end{subfigure} 
	\begin{subfigure}[b]{0.057\textwidth}
		\centering
		\includegraphics[width=\textwidth]{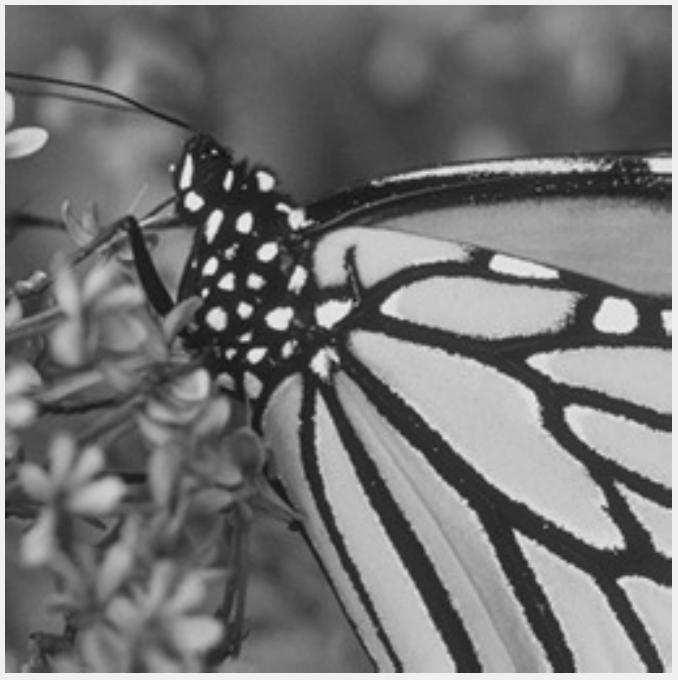} \vspace{-7mm}
		\caption*{\tiny  monarch}
	\end{subfigure} \vspace{-5mm}	
	\caption{\small  Images used for nonnegative matrix completion.  }\label{fig_nmcimgs}
\end{figure}

\begin{figure}	
	
	\begin{subfigure}[b]{0.16\textwidth}
		\centering
		\includegraphics[width=\textwidth]{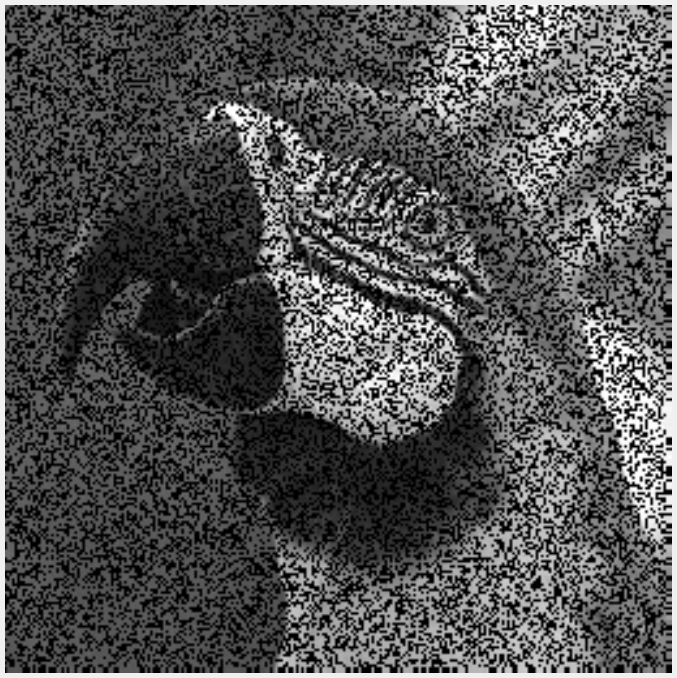} 
		\caption{}
	\end{subfigure} 
	\begin{subfigure}[b]{0.16\textwidth}
		\centering
		\includegraphics[width=\textwidth]{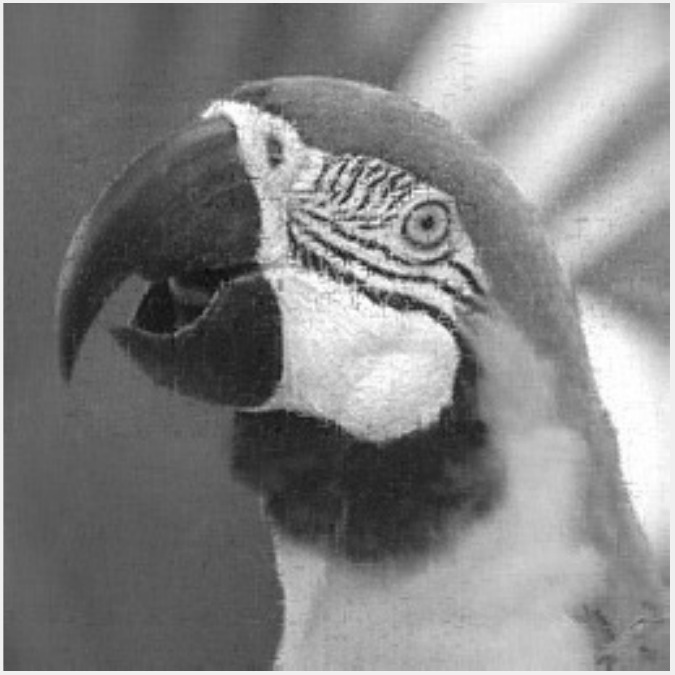} 
		\caption{}
	\end{subfigure} 
	\begin{subfigure}[b]{0.16\textwidth}
		\centering
		\includegraphics[width=\textwidth]{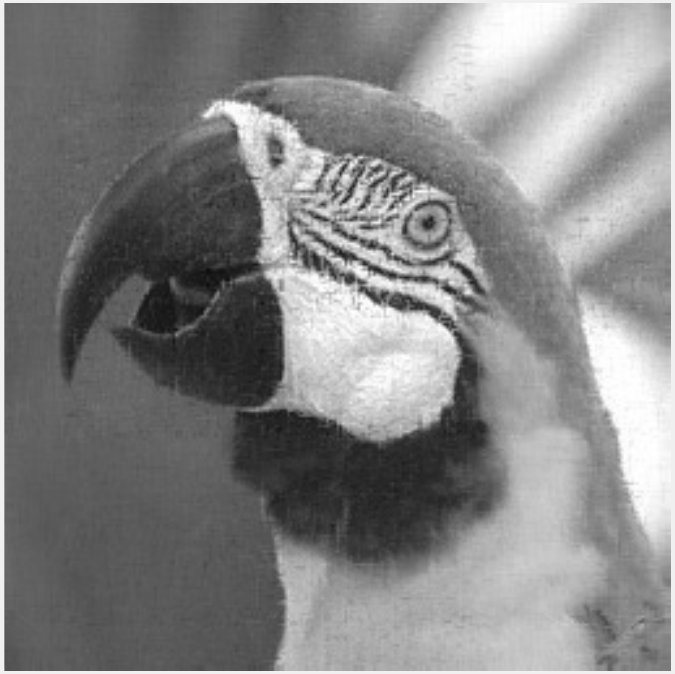} 
		\caption{}
	\end{subfigure} 	
	
	\begin{subfigure}[b]{0.16\textwidth}
		\centering
		\includegraphics[width=\textwidth]{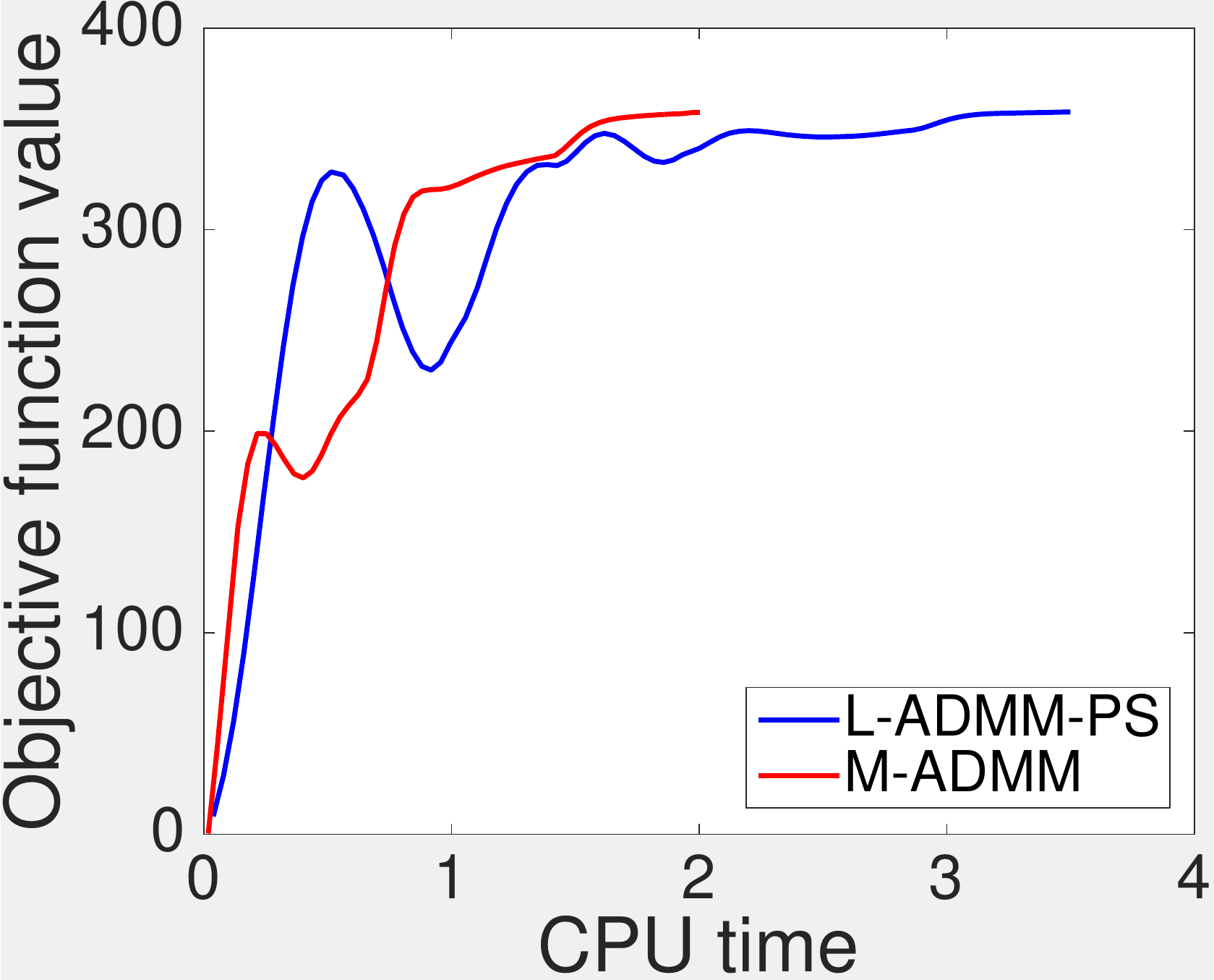} 
		\caption{}
	\end{subfigure} 
	\begin{subfigure}[b]{0.16\textwidth}
		\centering
		\includegraphics[width=\textwidth]{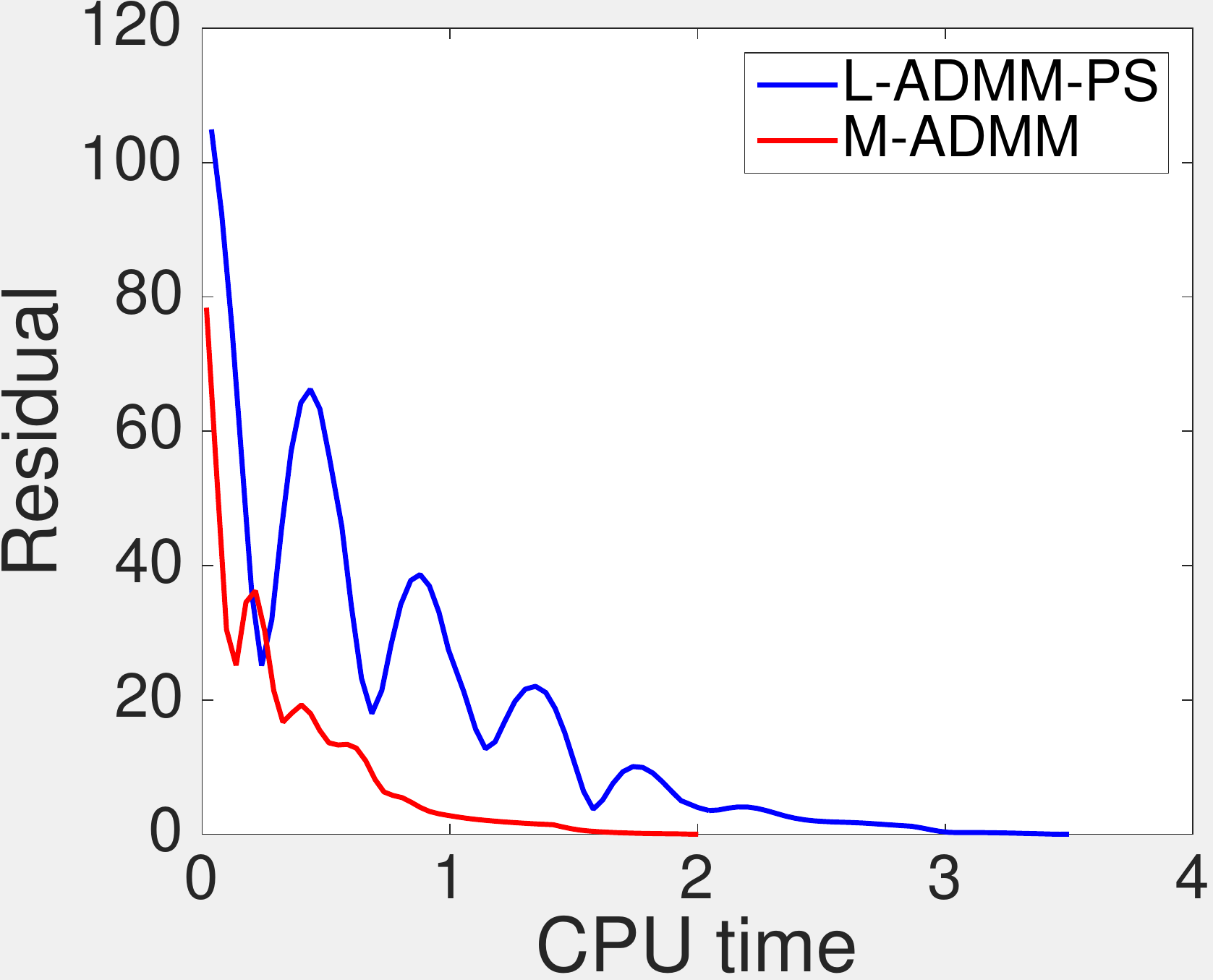} 
		\caption{}
	\end{subfigure} 
	\begin{subfigure}[b]{0.16\textwidth}
		\centering
		\includegraphics[width=\textwidth]{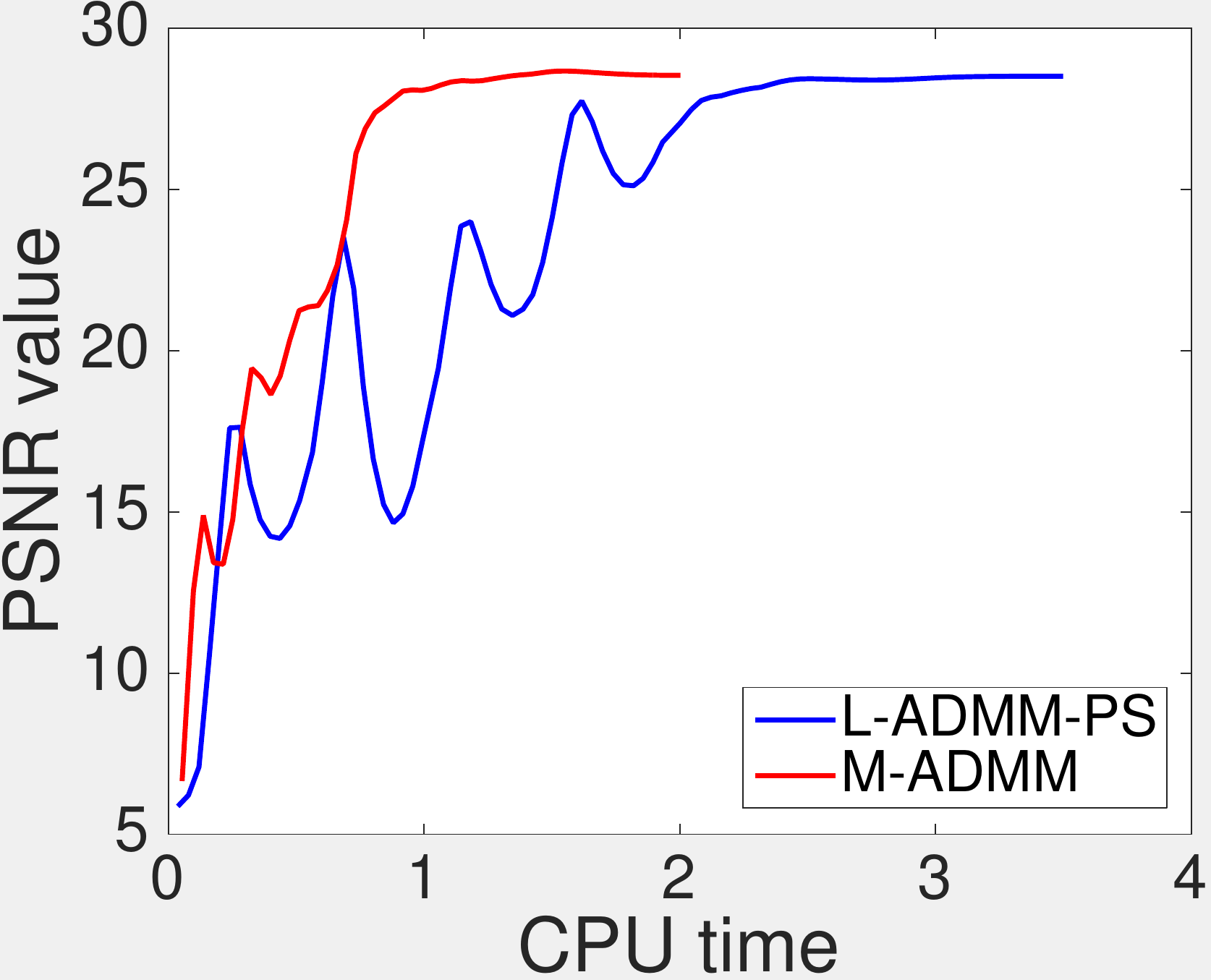} 
		\caption{}
	\end{subfigure} 
	\caption{\small Top row: the observed noisy image (left), recovered image by L-ADMM-PS (middle), and recovered image by M-ADMM (right). Bottom row: plots of  $f(\xk)$  v.s. CPU time (left), plots of $\norm{\Axk-\b}$ v.s. CPU time (middle), and PSNR values v.s. CPU time (right).  }\label{fig_nmcresexm}
\end{figure}

\begin{table}[t]
	\caption{\small Numerical comparison on the image inpainting.}
	\scriptsize
	\centering
	\begin{tabular}{c|ccc|ccc}
		\hline			
		& \multicolumn{3}{|c|}{L-ADMM-PS} & \multicolumn{3}{|c}{ M-ADMM }  	\\
		images	& PSNR & CPU & $\#$ Iter. & PSNR & CPU & $\#$Iter.	\\ \hline
		parrot & 28.51 & 3.50 & 87 & 28.54 & 2.00 & 55  \\    
		barbara & 27.69 & 3.36 & 85 & 27.72 & 2.27 & 60  \\    
		boat & 28.91 & 3.54 & 85 & 28.93 & 2.21 & 58  \\    
		cameraman & 26.06 & 3.33 & 84 & 26.08 & 2.15 & 58  \\    
		foreman & 31.83 & 3.80 & 86 & 31.84 & 2.06 & 54  \\    
		house & 31.26 & 3.48 & 87 & 31.26 & 2.29 & 56  \\    
		lena & 27.65 & 3.55 & 85 & 27.68 & 2.33 & 62  \\    
		monarch & 25.29 & 3.47 & 85 & 25.33 & 2.52 & 63  \\        \hline 
	\end{tabular} \label{tab_nmcres}   
\end{table}

\subsection{Solving Nonnegative Matrix Completion}
In this subsection, we show how to use Gauss-Seidel ADMM to solve a class of problems ($n>2$) with the condition (\ref{sepcondi}) being satisfied. We consider the following nonnegative noisy matrix completion problem \cite{liu2013linearized} 
\begin{equation}\label{eqnnmc}
\min_{\X,\E} \ \norm{\X}_*+\frac{\lambda}{2}\norm{\E}^2, \text{ s.t.} \ \Pomega(\X)+\E=\B, \ \X\geq \bm{0},
\end{equation}
where $\Omega$ is an index set and $\Pomega$ is a linear mapping that keeps the entries in $\Omega$ unchanged and those outside $\Omega$ zeros.
The above problem can be reformulated as a 3 blocks problem by (94) in \cite{liu2013linearized} and then solved by  L-ADMM-PS. We    instead reformulate (\ref{eqnnmc}) as
\begin{equation}\label{eqnnmc3}
\begin{split}
&\min_{\X,\E,\Z} \  \norm{\X}_*+\frac{\lambda}{2}\norm{\E}^2, \\
\text{ s.t.} \ & \Pomega(\Z)+\E=\B, \ \X=\Z,\ \Z\geq \bm{0}.
\end{split}
\end{equation}
Note that  (\ref{sepcondi}) holds for (\ref{eqnnmc3}) with the partition  $\{\X,\E\}$ and $\{\Z\}$. Thus (\ref{eqnnmc3})  can be solved  using (\ref{updatexbone})-(\ref{updatexbtwo}) with closed form solutions for each variable. We still refer to this method as M-ADMM in this experiment.

We consider the same image inpainting problem in \cite{liu2013linearized}  which is to fill in the missing pixel values of a corrupted image. As the pixel values are nonnegative, the image inpainting problem can be solved by (\ref{eqnnmc}).
The corrupted image is generated from the original image by sampling 60$\%$
of the pixels uniformly at random and adding Gaussian noise with mean zero
and standard deviation 0.1. We use the same adaptive penalty to update $\betak$ as \cite{liu2013linearized}. We set $\lambda=10$, $\epsilon_1  = 10^{-3}$, $\epsilon_2  = 10^{-4}$ and $\beta^{(0)}=\min{(d_1,d_2)}\epsilon_2$, where $d_1\times d_2$ is the size of $\X$. We update $\betakk=\max{(10\betak,10^6)}$ when $\max_i(\betak\norm{\xkk_i-\xk_i}/\norm{\b})\leq \epsilon_1$. The stopping criteria are $\max_i(\norm{\xkk_i-\xk_i}/\norm{\b}) \leq \epsilon_2$ and $\norm{\Axk-\b}/\norm{\b}\leq \epsilon_1$. We test on 8 images, all with size $256\times 256$,  in Figure \ref{fig_nmcimgs} and evaluate the recovery performance based on the PSNR value. The higher PSNR value indicates better recovery performance. The quantitative results are reported in Table \ref{tab_nmcres} and Figure \ref{fig_nmcresexm} gives more results test on the parrot image. It can be seen that, with  slightly better recovery performance, M-ADMM converges faster than L-ADMM-PS. The improvement benefits from the sequential updating of  $\{\X\}$ and $\{\Z,\E\}$ and avoids computing of the majorant surrogate as that in L-ADMM-PS.

		\section{Conclusions  }
		\label{sec7}
		
		This paper revisits   ADMM, an old but reborn method for convex problems with linear constraint. Many previous ADMMs can be categorized into the Gauss-Seidel ADMMs and Jacobian ADMMs according to different updating orders of the primal variables. We observed that many previous ADMMs update the primal variables by minimizing different majorant functions. Then we proposed the majorant first-order surrogate functions and presented the unified frameworks with unified convergence analysis. 	They not only draw the connections with existing ADMMs, but also can be used to solve new problems with non-separable objectives. The convergence bound show that the convergence speed depends on the tightness of the used majorant functions. We then analyzed how to improve the tightness to improve the efficiency. We improve Jacobian ADMMs by introducing the   Mixed Gauss-Seidel and Jacobian ADMM and the backtracking technique. We also discussed  how to perform variable partition for efficient implementations. Experiments on both synthesized and real-world data well demonstrated the effectiveness of our new ADMMs.	
		
		In the future, one may consider  extending our unified analysis based on MM to develop new ADMMs or solve other problems, e.g.,  strongly convex or nonconvex problems, and other ADMMs, e.g., stochastic ADMMs.


		{
			\bibliographystyle{ieee}
			\bibliography{supp}
		}
		
		
		\appendix

 \begin{lemma}\label{Lem_Pythagoras}
 	Given any $\bm{a}$, $\bm{b}$, $\bm{c}$, $\bm{d}$ and $\G\succeq0$ of compatible sizes, we have
 	\begin{align}
 	&\langle\bm{a}-\bm{b},\bm{c}-\bm{a}\rangle_{\G}\notag\\
 	=&\frac{1}{2}\left(\| \bm{b}-\bm{c} \|_{\G}^2-\|\bm{a}-\bm{c}\|_{\G}^2-\|\bm{a}-\bm{b}\|_{\G}^2\right),\label{ineq1} \\
 	&\langle\bm{a}-\bm{b},\bm{c}-\bm{d}\rangle \notag \\
 	=&\frac{1}{2}\left(\| \bm{a}-\bm{d} \|^2-\|\bm{a}-\bm{c}\|^2-\|\bm{b}-\bm{d}\|^2+\|\bm{b}-\bm{c}\|^2\right).\label{Lem_Pythagoras2}
 	\end{align}
 \end{lemma}

 \begin{lemma} \textit{\em{\textit{\textbf{(Combination Rules for Majorant First-Order Surrogates)}}}}\label{lemma5}
 	Let $\hat{f}\in\mathcal{S}_{\{\bL_i\}_{i=1}^n}(f,\bm{\kappa})$ and $\hat{f}'\in\mathcal{S}_{\{\bL'_i\}_{i=1}^n}(f',\bm{\kappa})$. Then the following combination rules hold:
 	\begin{itemize}
 		\item \textbf{Linear combination}: for any $\alpha, \alpha'>0$, $\alpha f+\alpha'f'$ is a majorant surrogate function in $\mathcal{S}_{\{\alpha \bL_i+\alpha'\bL'_i\}_{i=1}^n}(\alpha f+\alpha' f',\bm{\kappa})$;
 		\item \textbf{Transitivity}: let $F\in\mathcal{S}_{\{\bL''_i\}_{i=1}^n}(\hat{f},\bm{\kappa})$. Then $F$ is a majorant surrogate in $\mathcal{S}_{\{\bL_i+\bL''_i\}_{i=1}^n}(f,\bm{\kappa})$.
 	\end{itemize}	
 \end{lemma}

 \noindent\textbf{Proof of Lemma \ref{lemma44}.}   We deduce
 \begin{align}
 &f(\x)	\overset{\text{\ding{172}}}{\leq}
 \hat{f}(\x)\notag\\
 =&\left(\hat{f}(\x)-\frac{1}{2}\sumi\|\x_i-\bm{\kappa}_i\|^2_{\bP_i}\right)+\frac{1}{2}\sumi\|\x_i-\bm{\kappa}_i\|^2_{\bP_i}\notag\\
 \overset{\text{\ding{173}}}{\leq}&\left(\hat{f}(\y)-\frac{1}{2}\sumi\|\y_i-\bm{\kappa}_i\|^2_{\bP_i}\right)-\langle\mathbf{u},\y-\x\rangle\notag\\
 &+\sumi\langle\x_i-\bm{\kappa}_i,\y_i-\x_i\rangle_{\bP_i} 
 +\frac{1}{2}\sumi\|\x_i-\bm{\kappa}_i\|^2_{\bP_i}  \notag\\
 \overset{\text{\ding{174}}}{\leq}& \left(f(\y)+\frac{1}{2}\sumi\|\y_i-\bm\kappa_i\|^2_{\bL_i-\bP_i}\right)-\langle\mathbf{u},\y-\x\rangle\notag\\
 &+\sumi\langle\x_i-\bm{\kappa}_i,\y_i-\x_i\rangle_{\bP_i}+\frac{1}{2}\sumi\|\x_i-\bm{\kappa}_i\|^2_{\bP_i}    \notag\\
 \overset{\text{\ding{175}}}{=}&\left(f(\y)+\frac{1}{2}\sumi\|\y_i-\bm\kappa_i\|^2_{\bL_i-\bP_i}\right)-\langle\mathbf{u},\y-\x\rangle\notag\\
 &-\frac{1}{2}\sumi\left(\|\x_i-\bm{\kappa}_i\|^2_{\bP_i}+\|\y_i-\x_i\|^2_{\bP_i}-\|\y_i-\bm{\kappa}_i\|^2_{\bP_i}\right)\notag\\
 &+\frac{1}{2}\sumi\|\x_i-\bm{\kappa}_i\|_{\bP_i}^2 \notag\\
 =&f(\y)-\langle\mathbf{u},\y-\x\rangle+\frac{1}{2}\sum_{i=1}^{n}\left(\|\y_i-\bm\kappa_i\|^2_{\bL_i}-\|\y_i-\x_i\|^2_{\bP_i}\right),\notag
 \end{align}
 where   \ding{172} is from the fact that $\hatf$ is a majorant function of $f$, \ding{173} is from the convexity of $\hat{f}(\x)-\frac{1}{2}\sumi\|\x_i-\bm{\kappa}_i\|^2_{\bP_i}$ (or $\hat{f}$ is $\{\bP_i\}_{i=1}^n$-strongly convex),  \ding{174} uses (\ref{keylm11}), and \ding{175}  is from (\ref{ineq1}).
 $\hfill\blacksquare$

 \noindent
 \textbf{Proof of Lemma \ref{lem5}.} By using (\ref{ineq1}), for any $\x$ and $\y$, we have 
 \begin{align}
 &\frac{1}{2}\|\A\x-\b\|^2-\frac{1}{2}\|\A\y-\b\|^2\label{eqn5454545}\\
 =&\frac{1}{2}\|\A(\x-\y)\|^2+\langle\A(\x-\y),\A\y-\b\rangle\notag\\
 \leq&\frac{1}{2} \sumi  \|\x_i-\y_i\|_{\bL'_i}^2+\langle\A(\x-\y),\A\y-\b\rangle \label{eqn555552} \\
 \leq&\frac{1}{2}\sumi \left( \|\x_i-\y_i\|_{\G_i}^2+\|\A_i\left(\x_i-\y_i\right)\|^2 \right)\label{proofLem61}\\
 &+\sumi\langle\A_i(\x_i-\y_i),\A\y-\b\rangle\notag\\
 =&\frac{1}{2}\sumi\|\x_i-\y_i\|_{\G_i}^2\notag \\&+\frac{1}{2}\sumi\left(\|\A_i\left(\x_i-\y_i\right)+\A\y-\b\|^2-\|\A\y-\b\|^2\right)\notag\\
 =&\frac{1}{2}\sumi\left(\|\x_i-\y_i\|_{\G_i}^2+\left\|\A_i\x_i+\sum_{j\neq i}\A_j\y_j-\b\right\|^2\right)\notag \\
 &-\frac{n}{2}\|\A\y-\b\|^2\notag,
 \end{align}
 where (\ref{eqn555552}) holds for some $\bL'_i$'s; e.g., we can choose $\bL'_i\succeq n\A_i^\top\A_i$, and (\ref{proofLem61}) uses $\G_i\succeq \bL'_i-\A_i^\top\A_i$. 	Note that $r(\x)$ is convex and 
 (\ref{eqn5454545})-(\ref{eqn555552}) imply that (\ref{lipmuv}) holds. Thus, $r$ is  $\{\bL'_i\}_{i=1}^n$-smooth.
 By the definition of $\hatr$ in (\ref{sqrsurr2}), the above inequality implies that $r(\x)\leq\hatr(\x)$. Furthermore, it is easy to obtain (\ref{sqrsurr1})  by substituting  $\G_i=\eta_i\I-\A_i^\top\A_i$ into (\ref{proofLem61}).		
 $\hfill\blacksquare$

		We give the proof of Theorem \ref{them4}. 
		In the following, we define
		\begin{equation}\label{hatlambdatb}
		\hatlambdakk=\lambdak+\betak(\Abone\xbonekk+\Abtwo\xbtwok-\b).
		\end{equation}
		
		\begin{proposition}\label{pro1madmm}
			In Algorithm \ref{alg4}, under the assumptions of Theorem \ref{them4}, for any $\x$, we have
			\begin{align}
			 &f(\xkk)-f(\x)-\langle\A^\top\hatlambdakk,\x-\xkk\rangle \notag \\
			\leq& \frac{\betak}{2}\sum_{j=1}^2\left(\|\x_{B_j}-\xk_{B_j}\|_{\Hk_j}^2-\|\x_{B_j}-\xkk_{B_j}\|_{\Hkk_j}^2\right)\notag \\
			&-\frac{\betak}{2}\|\xbtwokk-\xbtwok\|^2_{ \Kk_2},\label{lm5}
			\end{align}
			where $\Hk_1=\text{Diag}\left\{\frac{1}{\betak}\bL_i+\A_i^\top\A_i+\Gk_i,i\in B_1\right\}-\Abone^\top\Abone$,  $\Hk_2=\text{Diag}\left\{\frac{1}{\betak}\bL_{i}+\A_i^\top\A_i+\Gk_i,i\in B_2\right\}$ and $\Kk_2=\text{Diag}\{\A_i^\top\A_i+\Gk_i,i\in B_2\}$.
		\end{proposition}

			\begin{proposition}\label{pro2madmm}
				In Algorithm \ref{alg4}, for any $\blambda$, we have  
				\begin{align}
				&\langle\A\xkk-\b,\blambda-\hatlambdakk\rangle+\frac{\beta^{(0)}\alpha}{2}\|\Axkk-\b\|^2\notag\\
				\leq&\frac{\betak}{2}\left(\|\blambda-\lambdak\|^2_{\Hk_3}-\|\blambda-\lambdakk\|^2_{\Hkk_3}\right)\notag\\
				&+\frac{\betak}{2}\|\xbtwokk-\xbtwok\|_{\Kk_2}^2\label{lm7},
				\end{align}
				where
				$\Hk_3=\left({1/\betak}\right)^2\I$ and $\alpha=\min\left\{\frac{1}{2},\frac{\tau}{2\norm{\Abtwo}_2^2}\right\}$.
			\end{proposition}


		\noindent\textbf{Proof of Theorem \ref{them4}.} Let $\x=\x^*$ and $\blambda=\blambda^*$ in (\ref{lm5}) and (\ref{lm7}). We have
		\begin{align*}
			&f(\xkk)-f(\x^*)+\langle \A^\top\blambda^*,\xkk-\x^*\rangle\notag\\
			&+\frac{\beta^{(0)}\alpha}{2}\|\Axkk-\b\|^2\notag\\
			\leq&\langle \A^\top(\blambda^*-\hatlambdakk),\xkk-\x^*\rangle-\langle\Axkk-\b,\blambda^*-\hatlambdakk\rangle\notag\\
			&+\frac{\betak}{2}\left(\sum_{i=j}^2\|\x_{B_j}^*-\xk_{B_j}\|^2_{\Hk_j}+\|\blambda^*-\lambdak\|^2_{\Hk_{3}}\right) \\
			&-\frac{\betak}{2}\left(\sum_{i=j}^2\|\x_{B_j}^*-\xkk_{B_j}\|^2_{\Hkk_j}+\|\blambda^*-\lambdakk\|^2_{\Hkk_{3}}\right)\notag\\
			=&\frac{\betak}{2}\left(\sum_{i=j}^2\|\x_{B_j}^*-\xk_{B_j}\|^2_{\Hk_j}+\|\blambda^*-\lambdak\|^2_{\Hk_{3}}\right) \\
			&-\frac{\betak}{2}\left(\sum_{i=j}^2\|\x^*_{B_j}-\xkk_{B_j}\|^2_{\Hkk_j}+\|\blambda^*-\lambdakk\|^2_{\Hkk_{3}}\right),
		\end{align*}
		where the last equation uses the fact $\A\x^*=\b$. Note that $\sum_{k=0}^K\gammak=1$. Multiplying $\gammak$ on both sides of the above inequalities and summing them from $0$ to $K$, we have
		\begin{align*}
			&\sum_{k=0}^{K}\gammak f(\xkk)-f(\x^*)+\left\langle \A^\top\blambda^*,\sum_{k=0}^{K}\gammak\xkk-\x^*\right\rangle\notag\\
			&+\frac{\beta^{(0)}\alpha}{2}\sum_{k=0}^{K}\gammak\|\Axkk-\b\|^2\notag\\ \leq&\frac{\sum_{j=1}^2\|\x_{B_j}^*-\x_{B_j}^0\|^2_{\HH^0_j}+\|\blambda^*-\blambda^0\|^2_{\HH^0_{3}}}{2\sum_{k=0}^{K}\left(\betak\right)^{-1}}.
		\end{align*}
		By the definition of $\bar{\x}^K$ and the convexity of $f$ and $\|\cdot\|^2$, we have
		\begin{align*}
			&f(\bar{\x}^K)-f(\x^*)+\langle \A^\top\blambda^*,\bar{\x}^K-\x^*\rangle+\frac{\beta^{(0)}\alpha}{2}\|\A\bar{\x}^K-\b\|^2\notag\\
			\leq&\sum_{k=0}^{K}\gammak f(\xkk)-f(\x^*)+\left\langle \A^\top\blambda^*,\sum_{k=0}^{K}\gammak\xkk-\x^*\right\rangle\notag\\
			&+\frac{\beta^{(0)}\alpha}{2}\sum_{k=0}^{K}\gammak\|\Axkk-\b\|^2\notag\\ \leq&\frac{\sum_{j=1}^2\|\x_{B_j}^*-\x_{B_j}^0\|^2_{\HH^0_j}+\|\blambda^*-\blambda^0\|^2_{\HH^0_{3}}}{2\sum_{k=0}^{K}\left(\betak\right)^{-1}}.
		\end{align*}
		The proof is completed.
		$\hfill\blacksquare$

		\vspace{-40pt}
		\begin{IEEEbiography}[{\includegraphics[width=1in,height=1.25in,clip,keepaspectratio]{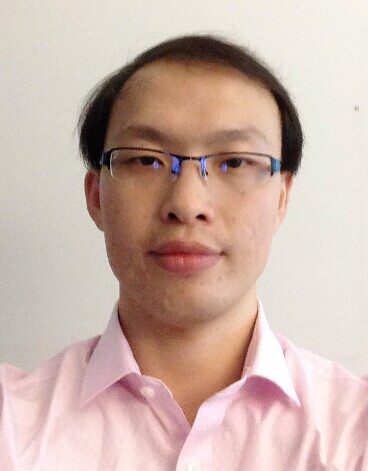}}]{Canyi Lu} received the bachelor degree in mathematics from the Fuzhou University in 2009, and the master degree in the pattern recognition and intelligent system from the University of Science and Technology of China in 2012. He is currently a Ph.D. student with the Department of Electrical and Computer Engineering at the National University of Singapore. His current research interests include computer vision, machine learning, pattern recognition and optimization. He was the winner of the Microsoft Research Asia Fellowship 2014.
		\end{IEEEbiography}
		\vspace{-40pt}
		
		\begin{IEEEbiography}[{\includegraphics[width=1in,height=1.25in,clip,keepaspectratio]{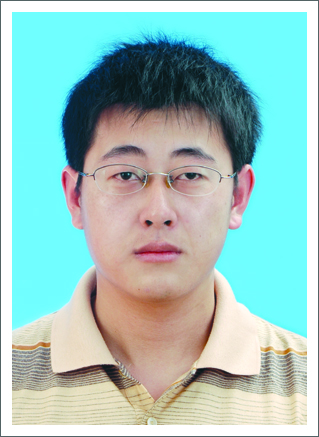}}]{Jiashi Feng} is currently an assistant Professor in the department of electrical and computer engineering in the National University of Singapore. He got his B.E. degree from University of Science and Technology, China in 2007 and Ph.D. degree from National University of  Singapore in 2014. He was a postdoc researcher at University of California from 2014 to 2015. His current research interest focus on machine learning and computer vision techniques for large-scale data analysis. Specifically, he has done work in object recognition, deep learning, machine learning, high-dimensional statistics and big data analysis.
		\end{IEEEbiography}
		\vspace{-40pt}
		
		\begin{IEEEbiography}[{\includegraphics[width=1in,height=1.25in,clip,keepaspectratio]{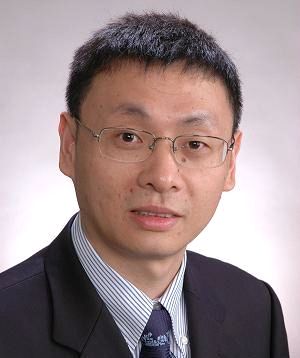}}]{Shuicheng Yan}
			is currently an Associate Professor at the Department of Electrical and Computer Engineering at National University of Singapore, and the founding lead of the Learning and Vision Research Group (http://www.lv-nus.org). Dr. Yan's research areas include machine learning, computer vision and multimedia, and he has authored/co-authored hundreds of technical papers over a wide range of research topics, with Google Scholar citation $>$30,000 times and H-index 64. He is ISI Highly-cited Researcher, 2014 and IAPR Fellow 2014. He has been serving as an associate editor of IEEE TKDE, TCSVT and ACM Transactions on Intelligent Systems and Technology (ACM TIST). He received the Best Paper Awards from ACM MM'13 (Best Paper and Best Student Paper), ACM MM’12 (Best Demo), PCM'11, ACM MM’10, ICME’10 and ICIMCS'09, the runner-up prize of ILSVRC'13, the winner prize of ILSVRC’14 detection task, the winner prizes of the classification task in PASCAL VOC 2010-2012, the winner prize of the segmentation task in PASCAL VOC 2012, the honourable mention prize of the detection task in PASCAL VOC'10, 2010 TCSVT Best Associate Editor (BAE) Award, 2010 Young Faculty Research Award, 2011 Singapore Young Scientist Award, and 2012 NUS Young Researcher Award.
		\end{IEEEbiography}
		\vspace{-40pt}
		\begin{IEEEbiography}[{\includegraphics[width=1in,height=1.25in,clip,keepaspectratio]{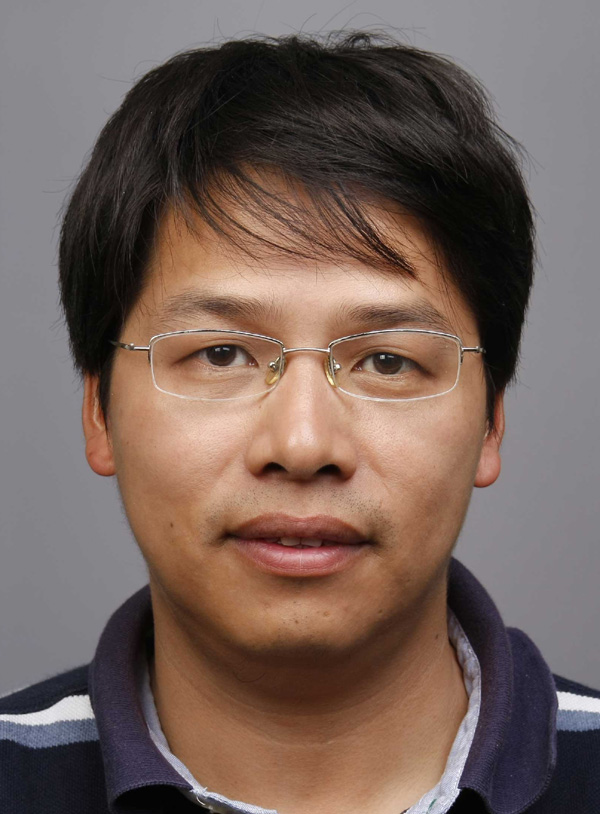}}]{Zhouchen Lin}
			received the Ph.D. degree in Applied Mathematics from Peking University, in 2000. He is currently a Professor at Key Laboratory of Machine Perception (MOE), School of Electronics Engineering and Computer Science, Peking University. He is also a Chair Professor at Northeast Normal University and a Guest Professor at Beijing Jiaotong University. Before March 2012, he was a Lead Researcher at Visual Computing Group, Microsoft Research Asia. He was a Guest Professor at Shanghai Jiaotong University and Southeast University, and a Guest Researcher at Institute of Computing Technology, Chinese Academy of Sciences. His research interests include computer vision, image processing, computer graphics, machine learning, pattern recognition, and numerical computation and optimization. He is an Associate Editor of IEEE Trans. Pattern Analysis and Machine Intelligence and International J. Computer Vision, an area chair of CVPR 2014, ICCV 2015, NIPS 2015 and AAAI 2016, and a Senior Member of the IEEE.
		\end{IEEEbiography}
		
	\newpage
	\onecolumn
	~\
	\begin{center}
		\LARGE\textbf{Supplementary Material}
	\end{center}
	~\
	
	This document contains two parts. First, we give the  proofs of some lemmas and propositions which are used to prove Theorem \ref{them4}.  Second, we give the implementation details of some problems in  the experiments.
	
		~\\\
	\noindent\textbf{\large 1. \ Proofs}
		~\\\
		
	   \noindent\textbf{Proof of Lemma \ref{lemma5}.}  Lemma \ref{lemma5} is obvious by using the  definition of the majorant first order surrogate function and the following lemma.
	 
	  \begin{lemma}\label{lem33}
	  	Let $f, f': \mathbb{R}^{p_1}\times\cdots\times\mathbb{R}^{p_n}\rightarrow\mathbb{R}$ be convex, and $\{\bL_i\}_{i=1}^n$-smooth and $\{\bL'_i\}_{i=1}^n$-smooth, respectively. We only consider two cases: (1) if $\bL_i\succeq\bL'_i$, define $\max\{\bL_i,\bL'_i\}=\bL_i$; (2) if $\bL'_i\succeq\bL_i$, define $\max\{\bL_i,\bL'_i\}=\bL'_i$. Then $f-f'$ is $\{\max\{\bL_i,\bL'_i\}\}_{i=1}^n$-smooth, and $f+f'$ is $\{\bL_i+\bL'_i\}_{i=1}^n$-smooth.
	  \end{lemma}

	 \noindent\textbf{Proof of Lemma \ref{lem33}.} Let $h=f-f'$. By using (\ref{lipmuv}) and the convexity of $f$ and $f'$, for any $\x=[\x_1;\cdots;\x_n]$ and $\y=[\y_1;\cdots;\y_n]$ with $\x_i, \y_i\in\mathbb{R}^{p_i}$, $i=1,\cdots,n$.  we have
	 \begin{eqnarray*}
	 	0 \leq  f(\x)-f(\y)-\langle\nabla f(\y),\x-\y\rangle  \leq\frac{1}{2}\sum_{i=1}^{n}\|\x_i-\y_i\|_{\bL_i}^2, \\
	 	-\frac{1}{2}\sum_{i=1}^{n}\|\x_i-\y_i\|_{\bL'_i}^2\leq -f'(\x)+f'(\y)+\langle\nabla f'(\y),\x-\y\rangle\leq 0.
	 \end{eqnarray*}
	 Summing the above two inequalities we have
	 \begin{equation*}
	 |h(\x)-h(\y)-\langle\nabla h(\y),\x-\y\rangle|\leq\frac{1}{2}\sum_{i=1}^{n}\|\x_i-\y_i\|_{\max\{\bL_i,\bL'_i\}}^2.
	 \end{equation*}
	 Thus $h$ is  $\{\max\{\bL_i,\bL'_i\}\}_{i=1}^n$-smooth. It is easy to see that $f+f'$ is $\{\bL_i+\bL'_i\}_{i=1}^n$-smooth  by applying (\ref{lipmuv}) for $f$ and $f'$.
	 $\hfill\blacksquare$

		~\
		
		\noindent
		\textbf{Proof of Proposition \ref{pro1madmm}.} First,  for $i\in B_1$, by the optimality of $\xkk_i$ to   problem (\ref{updatexbone2}) in Algorithm \ref{alg4},  there exists $\u_i^{k+1}\in\partial \hat{f}_i^k(\xkk_i)$ such that
		\begin{align}
		&-\u^{k+1}_i	=\nabla \hatr_i^k(\xkk_i)\notag\\
		\overset{\text{\ding{172}}}{=}&\A_i^\top\left(\betak\A_i(\xkk_i-\xk_i)+\betak(\A\xk-\b)+\lambdak\right)+\betak\Gk_i(\xkk_i-\xk_i)\notag\\
		=&\A_i^\top\left(\betak(\Abone\xbonekk+\Abtwo\xbtwok-\b)+\lambdak\right)-\betak\A_i^\top\Abone(\xbonekk-\xbonek)  +\betak(\A_i^\top\A_i+\Gk_i)(\xkk_i-\xk_i)\notag\\
		\overset{\text{\ding{173}}}{=}&\A_i^\top\hatlambdakk-\betak\A_i^\top\Abone(\xbonekk-\xbonek)+\betak(\A_i^\top\A_i+\Gk_i)(\xkk_i-\xk_i)\notag,
		\end{align}
		where \ding{172}  uses the definition of $\hat{r}^k_i$ in (\ref{rone}), and   \ding{173} uses the definition of $\hatlambdakk$ in  (\ref{hatlambdatb}). 	
		A dot-product with $\xkk_i-\x_i$ on both sides of the above equation gives
		\begin{align}
		&-\langle\u_{B_1}^{k+1},\xbonekk-\xbone\rangle	=-\sumione\langle\u^{k+1}_i,\xkk_i-\x_i\rangle\notag\\
		=&\sumione\left\langle\A_i^\top\hatlambdakk-\betak\A_i^\top\Abone(\xbonekk-\xbonek),\xkk_i-\x_i \right\rangle +\sumione\left\langle \betak(\A_i^\top\A_i+\Gk_i)(\xkk_i-\xk_i), \xkk_i-\x_i\right\rangle\notag\\
		=&\langle\Abone^\top\hatlambdakk,\xbonekk-\xbone\rangle +\betak\langle\xbonekk-\xbonek,\xbonekk-\xbone\rangle_{\Kk_1-\Abone^\top\Abone}\notag \\
		=&\langle\Abone^\top\hatlambdakk,\xbonekk-\xbone\rangle+\frac{\betak}{2}\|\xbonekk-\xbonek\|^2_{\Kk_1-\Abone^\top\Abone} +\frac{\betak}{2}\|\xbone-\xbonekk\|^2_{\Kk_1-\Abone^\top\Abone} -\frac{\betak}{2}\|\xbone-\xbonek\|^2_{\Kk_1-\Abone^\top\Abone}\notag\\
		\overset{\text{\ding{172}}}{\geq}&\langle\Abone^\top\hatlambdakk,\xbonekk-\xbone\rangle+	\frac{\betak}{2}\|\xbone-\xbonek\|^2_{\Kk_1-\Abone^\top\Abone}-\frac{\betak}{2}\|\xbone-\xbonekk\|^2_{\Kk_1-\Abone^\top\Abone},\label{prop5main1}
		\end{align}
		where $\Kk_1=\text{Diag}\{\A_i^\top\A_i+\Gk_i,i\in B_1\}$  and \ding{172} uses $\|\xbonekk-\xbonek\|^2_{\Kk_1-\Abone^\top\Abone}\geq0$ due to (\ref{eqbteq2}).  
		
		Second, for $i\in B_2$, by the optimality of $\xkk_i$ to problem (\ref{updatexbtwo2}) in Algorithm \ref{alg4}, there exists $\u_i^{k+1}\in\partial \hat{f}_i^k(\xkk_i)$ such that
		\begin{align*}
		&-\u^{k+1}_i			=\nabla \hatr_i^k(\xkk_i)\\
		=&\A_i^\top\left(\betak\A_i(\xkk_i-\xk_i)+\betak(\Abone\xbonekk+\Abtwo\xbtwok-\b)\right)\\
		&+\A_i^\top\lambdak+\betak\Gk_i(\xkk_i-\xk_i)\\
		=&\A_i^\top\hatlambdakk+\betak(\A_i^\top\A_i+\Gk_i)(\xkk_i-\xk_i),
		\end{align*}
		where we use the definitions of $\hat{r}^k_i$ in (\ref{rtwo}) and $\hatlambdakk$ in (\ref{hatlambdatb}).
		A dot-product with $\xkk_i-\x_i$ on both sides of the above equation gives
		\begin{align}
		&-\langle\u_{B_2}^{k+1},\xbtwokk-\xbtwo\rangle	=-\sumitwo\langle\u^{k+1}_i,\xkk_i-\x_i\notag\rangle\\
		=&\sumitwo\left\langle\A_i^\top\hatlambdakk+\betak(\A_i^\top\A_i+\Gk_i)(\xkk_i-\xk_i),\xkk_i-\x_i \right\rangle\notag\\
		=&\langle\Abtwo^\top\hatlambdakk,\xbtwokk-\xbtwo\rangle+\betak\langle\xbtwokk-\xbtwok,\xbonekk-\xbone\rangle_{\Kk_2} \notag\\
		=&\langle\Abtwo^\top\hatlambdakk,\xbtwokk-\xbtwo\rangle-\frac{\betak}{2} \|\xbtwo-\xbtwok\|^2_{\Kk_2}  +\frac{\betak}{2} \|\xbtwo-\xbtwokk\|^2_{\Kk_2}+\frac{\betak}{2}\|\xbtwokk-\xbtwok\|^2_{\Kk_2},\label{prop5main2}
		\end{align}
		where $\Kk_2=\text{Diag}\{\A_i^\top\A_i+\Gk_i,i\in B_2\}$. 
		
		Third, note that $\hatf^k\in\mathcal{S}_{\{\bL_i,\bP_i\}_{i=1}^n}(f,\xk)$. By using (\ref{keylm22}), we have
		\begin{align}
		& f(\xkk)-f(\x)\notag\\
		\leq&\langle \u^{k+1},\xkk-\x\rangle+\frac{1}{2}\sumi\left(\|\x_i-\xk_i\|^2_{\bL_i}-\|\x_i-\xkk_i\|^2_{\bP_i}\right)\notag\\
		\leq&\langle \u^{k+1},\xkk-\x\rangle+\frac{1}{2}\sumi\left(\|\x_i-\xk_i\|^2_{\bL_i}-\|\x_i-\xkk_i\|^2_{\bL_i}\right)\notag\\
		=&\langle \u_{B_1}^{k+1},\xbonekk-\xbone\rangle+\langle \u_{B_2}^{k+1},\xbtwokk-\xbtwo\rangle +\frac{1}{2}\sum_{j=1}^2\left(\|\x_{B_j}-\x^k_{B_j}\|_{\bL_{B_j}}^2-\|\x_{B_j}-\x^{k+1}_{B_j}\|^2_{\bL_{B_j}}\right)\notag\\
		\overset{\text{\ding{172}}}{\leq}&-\langle\A^\top\hatlambdakk,\xkk-\x\rangle-\frac{\betak}{2}\|\xbtwokk-\xbtwok\|^2_{\K_2}+\frac{\betak}{2}\sum_{j=1}^2\left(\|\x_{B_j}-\xk_{B_j}\|_{\Hk_j}^2-\|\x_{B_j}-\xkk_{B_j}\|_{\Hkk_j}^2\right) \notag
		\end{align}
		where   $\bL_{B_j}=\text{Diag}\left\{\bL_i, i\in B_j\right\}$  and \ding{172} uses (\ref{prop5main1})-(\ref{prop5main2}), the definitions of $\Hk_j$ in Proposition \ref{pro1madmm} and the fact $\betakk\geq\betak$.
		The proof is completed.  $\hfill\blacksquare$

		~\
	
		\noindent
		\textbf{Proof of Proposition \ref{pro2madmm}.} By using line 4 of Algorithm \ref{alg3}, (\ref{Lem_Pythagoras2}) and the fact that $\betakk\geq\betak$, we have
		\begin{align}
		&\langle\A\xkk-\b,\blambda-\hatlambdakk\rangle		=\frac{1}{\betak}\langle\lambdakk-\lambdak,\blambda-\hatlambdakk\rangle\notag\\
		=&\frac{1}{2\betak}\left(\|\blambda-\lambdak\|^2-\|\blambda-\lambdakk\|^2\right) -\frac{1}{2\betak}\left(\|\hatlambdakk-\lambdak\|^2-\|\lambdakk-\hatlambdakk\|^2\right) \notag\\ &-\frac{1}{2\betak}\left(\|\hatlambdakk-\lambdak\|^2-\|\lambdakk-\hatlambdakk\|^2\right)\label{lm71}.
		\end{align}
		Consider the last two terms in (\ref{lm71}). We deduce
		\begin{align}
		&\frac{1}{2\betak}\left(\|\hatlambdakk-\lambdak\|^2-\|\lambdakk-\hatlambdakk\|^2\right)\notag\\
		\overset{\text{\ding{172}}}{=}&\frac{\betak}{2}\|\Abone\xbonekk+\Abtwo\xbtwok-\b\|^2-\|\Abtwo(\xbtwokk-\xbtwok)\|^2\notag\\
		=&\frac{\betak}{2}\left(\|\Abone\xbonekk+\Abtwo\xbtwok-\b\|^2-\|\xbtwokk-\xbtwok\|_{\Kk_2}^2\right)+\frac{\betak}{2}\|\xbtwokk-\xbtwok\|_{\Kk_2-\Abtwo^\top\Abtwo}^2	\notag\\
		\overset{\text{\ding{173}}}{\geq}&\frac{\betak}{2}\left(\|\Abone\xbonekk+\Abtwo\xbtwok-\b\|^2- \|\xbtwokk-\xbtwok\|_{\Kk_2}^2\right)+\frac{\tau}{\|\Abtwo\|^2_2}\|\Abtwo(\xbtwokk-\xbtwok)\|^2 \notag\\
		\overset{\text{\ding{174}}}{\geq}&{\betak\alpha}\|\Abone\xbonekk+\Abtwo\xbtwok-\b\|^2-\frac{\betak}{2}\|\xbtwokk-\xbtwok\|_{\Kk_2}^2+\betak\alpha\|\Abtwo(\xbtwokk-\xbtwok)\|^2\notag\\
		\geq&\frac{\betak\alpha}{2}\|\Abone\xbonekk+\Abtwo\xbtwokk-\b\|^2-\frac{\betak}{2}\|\xbtwokk-\xbtwok\|_{\Kk_2}^2\notag\\
		\overset{\text{\ding{175}}}{\geq}&\frac{\beta^{(0)}\alpha}{2}\|\Axkk-\b\|^2-\frac{\betak}{2}\|\xbtwokk-\xbtwok\|_{\Kk_2}^2\label{lm72},
		\end{align}
		where \ding{172} uses  (\ref{updatelambdaall}) and (\ref{hatlambdatb}), \ding{173} uses (\ref{eqbteq3}), \ding{174}  uses $\alpha=\min\left\{\frac{1}{2},\frac{\tau}{2\norm{\Abtwo}_2^2}\right\}$,
		and \ding{175} uses $\betak\geq\beta^{(k-1)}\geq\cdots\geq \beta^{(0)}$. 		
		The proof is completed by substituting (\ref{lm72}) into (\ref{lm71}). 
		$\hfill\blacksquare$
		
		~\

		~\\\\\\\\
		\noindent\textbf{\large 2. \ Implementation Details}
		~\
				
		~\\
		\noindent\textbf{\large 2.1 \ Latent Low-Rank Representation}
		~\

		Consider the following Latent Low-Rank Representation (LRR) problem
		\begin{equation}\label{latlrr111}
		\min_{\Z,\bL} \|\Z\|_*+\|\bL\|_*+\frac{\lambda}{2}\|\X\Z+\bL\X-\X\|_F^2, \ \text{s.t.} \ \bm{1}^T\Z=\bm{1}^T,
		\end{equation}
		Problem (\ref{latlrr111}) is equivalent to
		\begin{align}
		\min_{\Z,\bL,\E} \|\Z\|_*+\|\bL\|_*+\frac{\lambda}{2}\|\E\|_F^2, 
		\ \text{s.t.} \ \bm{1}^T\Z=\bm{1}^T,\ \X\Z+\bL\X-\X=\E.\label{latlrr222}
		\end{align}
		
		
		\begin{enumerate}[(a)]
			\item Solve (\ref{latlrr111}) by M-ADMM (2)			
			
			The augmented Lagrangian function of (\ref{latlrr111}) is
			\begin{equation*}
			\mathcal{L}(\Z,\bL,\blambda,\beta)=\|\Z\|_*+\|\bL\|_*+\frac{\lambda}{2}\|\X\Z+\bL\X-\X\|_F^2+\langle\blambda,\bm{1}^T\Z-\bm{1}^T\rangle+\frac{\beta}{2}\|\bm{1}^T\Z-\bm{1}^T\|^2.
			\end{equation*}
			It is easy to verify that $\frac{1}{2}\|\X\Z+\bL\X-\X\|_F^2$ is $\{L_1\I,L_2\I\}$-smooth, where $L_1=L_2=2\|\X\|_2^2$, and $\frac{1}{2}\|\bm{1}^T\Z-\bm{1}^T\|^2$  is $\eta$-smooth, where $\eta>\|\bm{1}\|^2$. By using these properties, M-ADMM (2) solves (\ref{latlrr111}) by the following updating rules
			
			\begin{equation*}
			\left\{
			\begin{aligned}
			\Z^{k+1}=&\arg\min_{\Z} \|\Z\|_*+\frac{\lambda L_1+\betak\eta}{2}\left\|\Z-\Z^k+\frac{\lambda\X^T(\X\Z^k+\bL^k\X-\X)+\bm{1}\left(\betak(\bm{1}^T\Z^k-\bm{1}^T)+\lambdak\right)}{\lambda L_1+\betak\eta} \right\|^2_F,  \\
			\bL^{k+1}=&\arg\min_{\bL}\|\bL\|_*+\frac{\lambda L_2}{2}\left\|\bL-\bL^k+\frac{(\X\Z^k+\bL^k\X-\X)\X^T}{L_2}\right\|_F^2,\\
			\lambdakk=&\lambdak+\betak(\bm{1}^T\Z^{k+1}-\bm{1}^T).
			\end{aligned}
			\right.
			\end{equation*}

			\item Solve (\ref{latlrr222}) by L-ADMM-PS (3)
			
			The augmented Lagrangian function of (\ref{latlrr222}) is
			\begin{align*}
			\mathcal{L}(\Z,\bL,\blambda,\beta)=&\|\Z\|_*+\|\bL\|_*+\frac{\lambda}{2}\|\E\|_F^2+\langle\blambda_1,\bm{1}^T\Z-\bm{1}^T\rangle+\frac{\beta}{2}\|\bm{1}^T\Z-\bm{1}^T\|^2\\
			&+\langle\blambda_2,\X\Z+\bL\X-\X-\E\rangle+\frac{\beta}{2}\|\X\Z+\bL\X-\X-\E\|_F^2.
			\end{align*}
			
			Note that $h(\Z,\bL,\E)=\frac{1}{2}\|\bm{1}^T\Z-\bm{1}^T\|^2+\frac{1}{2}\|\X\Z+\bL\X-\X-\E\|_F^2$ is $\{\eta_1\I,\eta_2\I,\eta_3\I\}$-smooth, where $\eta_1>\|\bm{1}\|^2+3\|\X\|_2^2$, $\eta_2>3\|\X\|_2^2$ and $\eta_3>3$.
			By using such a property, 
			L-ADMM-PS (3) solves (\ref{latlrr222}) by the following updating rules 
			
			\begin{equation*}
			\left\{
			\begin{aligned}
			\Z^{k+1}=&\arg\min_{\Z} \|\Z\|_*+\frac{\betak\eta_1}{2}\left\|\Z-\Z^k+\frac{\bm{1}(\blambda^k_1+\betak(\bm{1}^T\Z^k-\bm{1}^T))+\X^T(\blambda_2^k+\betak(\X\Z^k+\bL^k\X-\X-\E^k))}{\betak\eta_1} \right\|_F^2,  \\
			\bL^{k+1}=&\arg\min_{\bL} \|\bL\|_*+\frac{\betak\eta_2}{2}\left\|\bL-\bL^k+\frac{(\blambda^k_2+\betak(\X\Z^k+\bL^k\X-\X-\E^k))\X^T}{\betak\eta_2} \right\|_F^2,\\
			\E^{k+1}=&\arg\min_{\E}\frac{\lambda}{2}\|\E\|_F^2-\langle\blambda^k_2+\betak(\X\Z^{k}+\bL^k\X-\X-\E^k),\E\rangle+\frac{\betak\eta_3}{2}\|\E-\E^k\|_F^2,\\
			\lambdakk_1=&\lambdak_1+\betak(\bm{1}^T\Z^{k+1}-\bm{1}^T),\\
			\lambdakk_2=&\lambdak_2+\betak(\X\Z^{k+1}+\bL^{k+1}\X-\X-\E^{k+1}),
			\end{aligned}
			\right.
			\end{equation*}

			\item Solve (\ref{latlrr222}) by M-ADMM (3)
			
			M-ADMM (3) divides the variables $\{\Z,\bL,\E\}$ into two super blocks, i.e., $\{\Z\}$ and $\{\bL,\E\}$. Then it solves (\ref{latlrr222}) by the following updating rules
			\begin{equation*}
			\left\{
			\begin{aligned}
			\Z^{k+1}=&\arg\min_{\Z} \|\Z\|_*+\frac{\betak\eta_1}{2}\left\|\Z-\Z^k+\frac{\bm{1}(\blambda^k_1+\betak(\bm{1}^T\Z^k-\bm{1}^T))+\X^T(\blambda_2^k+\betak(\X\Z^k+\bL^k\X-\X-\E^k))}{\betak\eta_1} \right\|_F^2,  \\
			\bL^{k+1}=&\arg\min_{\bL} \|\bL\|_*+\frac{\betak\eta_2}{2}\left\|\bL-\bL^k+\frac{(\blambda^k_2+\betak(\X\Z^{k+1}+\bL^k\X-\X-\E^k))\X^T}{\betak\eta_2} \right\|_F^2,\\
			\E^{k+1}=&\arg\min_{\E}\frac{\lambda}{2}\|\E\|_F^2+\frac{\betak}{2}\left\|\X\Z^{k+1}+\bL^k\X-\X-\E+\frac{\lambdak_2}{\betak}\right\|^2_F+\frac{\betak\eta_3}{2}\|\E-\E^k\|_F^2,\\
			\lambdakk_1=&\lambdak_1+\betak(\bm{1}^T\Z^{k+1}-\bm{1}^T),\\
			\lambdakk_2=&\lambdak_2+\betak(\X\Z^{k+1}+\bL^{k+1}\X-\X-\E^{k+1}),
			\end{aligned}
			\right.
			\end{equation*}
			where $\eta_1=\|\bm{1}\|^2+\|\X\|_2^2$, $\eta_2>2\|\X\|_2^2$ and $\eta_3>1$.

\end{enumerate}

		~\\
		\noindent\textbf{\large 2.2 \ Nonnegative Matirx Completion}
		~\

		\begin{equation}\label{eqnnmc222222}
		\min_{\X,\E} \ \norm{\X}_*+\frac{\lambda}{2}\norm{\E}^2, \text{ s.t.} \ \Pomega(\X)+\E=\B, \ \X\geq \bm{0},
		\end{equation}
		\begin{enumerate}[(a)]
		\item  L-ADMM-PS 
			
			Problem (\ref{eqnnmc222222}) is equivalent to (see (94) in \cite{LADMPS})
			\begin{equation}\label{nmceqf1}
			\begin{split}
			&\min_{\X,\E,\Z} \  \norm{\X}_*+\frac{\lambda}{2}\norm{\E}^2, \\
			\text{ s.t.} \ & \Pomega(\X)+\E=\B, \ \X=\Z,\ \Z\geq \bm{0}.
			\end{split}
			\end{equation}
			The partial augmented Lagrangian function is
			\begin{align*}
			\mathcal{L}(\X,\E,\Z,\beta) = \norm{\X}_*+\frac{\lambda}{2}\norm{\E}^2  &+  \langle\blambda_1,  \Pomega(\X)+\E-\B  \rangle+\frac{\beta}{2}\|\Pomega(\X)+\E-\B\|^2 \\ 
			 & + \langle\blambda_2,  \X-\Z  \rangle+\frac{\beta}{2}\|\X-\Z\|^2.
			\end{align*}
			Then L-ADMM-PS solves (\ref{nmceqf1}) by the following updating rules 
			\begin{equation*}
			\left\{
			\begin{aligned}
			\X^{k+1}=&\arg\min_{\X} \|\X\|_*+\frac{\betak\eta_1}{2}\left\|\X-\X^k+\frac{\Pomega(\lambdak_1) +\lambdak_2 +\betak \Pomega(\X^k+\E^k-\B) + \betak (\X^k-\Z^k) }{\betak\eta_1} \right\|_F^2,  \\
			\E^{k+1}=&\arg\min_{\E}  \frac{\lambda}{2}\norm{\E}^2+ \langle\E,\lambdak_1+ \betak(\Pomega(\X^k)+\E^k-\B)\rangle + \frac{\betak\eta_2}{2}\norm{\E-\E^k}^2,  \\
			\Z^{k+1}=&\arg\min_{\Z\geq \mathbf{0}} \langle\lambdak_2+\betak(\X^k-\Z^k),-\Z\rangle + \frac{\betak}{2} \norm{\X^k-\Z}^2+  \frac{\betak\eta_3}{2} \norm{\X^k-\Z}^2, \\
			\lambdakk_1=&\lambdak_1+\betak(\Pomega(\X^{k+1})+\E^{k+1}-\B),\\
			\lambdakk_2=&\lambdak_2+\betak(\X^{k+1}-\Z^{k+1}),
			\end{aligned}
			\right.
			\end{equation*}
			where $\eta_1 > 3+2$, $\eta_2>3+2$ and $\eta_3 > 2$.
			
			\item  M-ADMM
			
			Problem (\ref{eqnnmc222222}) is equivalent to 
			\begin{equation}\label{eqnmc33}
			\begin{split}
			&\min_{\X,\E,\Z} \  \norm{\X}_*+\frac{\lambda}{2}\norm{\E}^2, \\
			\text{ s.t.} \ & \Pomega(\Z)+\E=\B, \ \X=\Z,\ \Z\geq \bm{0}.
			\end{split}
			\end{equation}
			The partial augmented Lagrangian function is
			\begin{align*}
			\mathcal{L}(\X,\E,\Z,\beta) = \norm{\X}_*+\frac{\lambda}{2}\norm{\E}^2  &+  \langle\blambda_1,  \Pomega(\Z)+\E-\B  \rangle+\frac{\beta}{2}\|\Pomega(\Z)+\E-\B\|^2 \\ 
			& + \langle\blambda_2,  \X-\Z  \rangle+\frac{\beta}{2}\|\X-\Z\|^2.
			\end{align*}
			Partition the three blocks into two super blocks $\{\X,\E\}$ and $\{\Z\}$.
			Then M-ADMM solves (\ref{eqnmc33}) by the following updating rules 
			\begin{equation*}
			\left\{
			\begin{aligned}
			\X^{k+1}=&\arg\min_{\X} \|\X\|_*+\frac{\betak }{2}\left\|\X-\Z^k + \frac{\lambdak_2}{\betak} \right\|_F^2,  \\
			\E^{k+1}=&\arg\min_{\E}  \frac{\lambda}{2}\norm{\E}^2+  \langle\lambdak_1,  \E  \rangle+\frac{\betak}{2}\|\Pomega(\Z^k)+\E-\B\|^2,  \\
			\Z^{k+1}=&\arg\min_{\Z\geq \mathbf{0}} \langle\lambdak_1,  \Pomega(\Z)\rangle+\frac{\betak}{2}\|\Pomega(\Z)+\E^{k+1}-\B\|^2 + \langle\lambdak_2,   -\Z  \rangle+\frac{\betak}{2}\|\X^{k+1}-\Z\|^2, \\
			\lambdakk_1=&\lambdak_1+\betak(\Pomega(\X^{k+1})+\E^{k+1}-\B),\\
			\lambdakk_2=&\lambdak_2+\betak(\X^{k+1}-\Z^{k+1}).
			\end{aligned}
			\right.
			\end{equation*}
			Note that the $\Z^{k+1}$ updating has a closed form solution.
		
		\end{enumerate}
		
		~\\
		\noindent\textbf{\large 3 \ A List of Problems Involved in Our Released Toolbox}
		~\
		
		Table \ref{tabLlistofproblems} gives a list of convex problems in compressed sensing solved by M-ADMM in our released LibADMM package. For each problem, we consider its specific structure to implement efficient M-ADMM by using several techniques proposed in this work.

\renewcommand{\arraystretch}{1.3}
\begin{table}[!h]
	\centering
	\caption{Applicability of the LibADMM package}
	\label{tabLlistofproblems}
	\begin{tabular}{c|cl|l|l}
		\hline
		Model               &  \multicolumn{2}{c|}{Problem}                                                                                                             & Function                              & \multicolumn{1}{l}{Description and Reference}                                                    \\ \hline\hline
		\multirow{12}{*}{\begin{tabular}[c]{@{}l@{}}Sparse\\ \\ models\end{tabular}} & \multicolumn{1}{l|}{\multirow{6}{*}{\begin{tabular}[c]{@{}l@{}}$\min_{\x}\ r(\x)$\\ \\ $\text{s.t.} \ \A\x=\b$\end{tabular}}}                       & $r(\x)=\norm{\x}_1$                                                                                                  & \mcode{l1}                          & \multicolumn{1}{l}{$\ell_1$  }                                       \\ \cline{3-5} 
		& \multicolumn{1}{l|}{}                                                                                                                               & $r(\x)=\sum_{g\in\mathcal{G}}\norm{\x_g}_2$                                                                          & \mcode{groupl1}                     & \multicolumn{1}{l}{Group Lasso  }                                            \\ \cline{3-5} 
		& \multicolumn{1}{l|}{}                                                                                                                               & $r(\x)=\norm{\x}_1+\lambda_2\norm{\x}_2^2$                                                                           & \mcode{elasticnet}                  & \multicolumn{1}{l}{Elastic net  }                                    \\ \cline{3-5} 
		& \multicolumn{1}{l|}{}                                                                                                                               & $r(\x)=\norm{\x}_1+\lambda_2\sum_{i=2}^p|x_i-x_{i-1}|$                                                               & \mcode{fusedl1}                     & \multicolumn{1}{l}{Fused Lasso  }                                   \\ \cline{3-5} 
		& \multicolumn{1}{l|}{}                                                                                                                               & $r(\x)=\norm{\A\text{Diag}(\x)}_*$                                                                                   & \mcode{tracelasso}                  & \multicolumn{1}{l}{Trace Lasso  }                                        \\ \cline{3-5} 
		& \multicolumn{1}{l|}{}                                                                                                                               & $r(\x)=\frac{1}{2}\norm{\x}^2_{\text{ksp}}$                                                                                   & \mcode{ksupport}                  & \multicolumn{1}{l}{$k$ support norm  }                                        \\ \cline{2-5} 
		& \multicolumn{1}{l|}{\multirow{6}{*}{\begin{tabular}[c]{@{}l@{}}$\min_{\x,\e}\ l(\e) + \lambda r(\x)$\\ \\ $\text{s.t.} \ \A\x+\e=\b$\end{tabular}}} & \multirow{6}{*}{\begin{tabular}[c]{@{}l@{}}$\l(\e)=\norm{\e}_1$\\ \\ $\l(\e)=\frac{1}{2}\norm{\e}_2^2$\end{tabular}} & \mcode{l1R}                         & \multicolumn{1}{l}{Reg. $\ell_1$}                                                                   \\ \cline{4-5} 
		& \multicolumn{1}{l|}{}                                                                                                                               &                                                                                                                      & \mcode{groupl1R}                    & \multicolumn{1}{l}{Reg. Group Lasso}                                                             \\ \cline{4-5} 
		& \multicolumn{1}{l|}{}                                                                                                                               &                                                                                                                      & \mcode{elasticnetR}                 & \multicolumn{1}{l}{Reg. Elastic net}                                                             \\ \cline{4-5} 
		& \multicolumn{1}{l|}{}                                                                                                                               &                                                                                                                      & \mcode{fusedl1R}                    & \multicolumn{1}{l}{Reg. Fused Lasso}                                                             \\ \cline{4-5} 
		& \multicolumn{1}{l|}{}                                                                                                                               &                                                                                                                      & \mcode{tracelassoR}                 & \multicolumn{1}{l}{Reg. Trace Lasso}                                                             \\   \cline{4-5} 
		& \multicolumn{1}{l|}{}                                                                                                                               &                                                                                                                      & \mcode{ksupportR}                 & \multicolumn{1}{l}{Reg. $k$ support norm}                                                             \\ \hline\hline
		\multirow{12}{*}{\begin{tabular}[c]{@{}l@{}}Low-rank\\ \\ matrix\\ \\ models\end{tabular}}  
		& \multicolumn{2}{l|}{$\min_{\X,\E} \ \norm{\X}_*+\lambda l(\E), \ \st \ \Pomega(\X)+\E=\M$}                                                                                                                                                                                  & \mcode{lrmcR}                       & \multicolumn{1}{l}{Reg. Low-rank matrix completion}                                              \\ \cline{2-5} 
		& \multicolumn{2}{l|}{$\min_{\X,\E}\ \norm{\X}_*+\lambda l(\E),\ \st \ \A=\B\X+\E$}                                                                                                                                                                                          & \mcode{lrr}                         & \multicolumn{1}{l}{Low-rank representation  }                                    \\ \cline{2-5} 
		& \multicolumn{2}{l|}{$\min_{\Z,\bL,\E}\ \norm{\Z}_*+\norm{\bL}_*+\lambda l(\E)$}                                                                                                                                                                                             & \multirow{2}{*}{\mcode{latlrr}}     & \multicolumn{1}{l}{\multirow{2}{*}{Latent low-rank representation  }}        \\
		& \multicolumn{2}{l|}{$\st \ \X\Z+\bL\X-\X=\E$}                                                                                                                                                                                                                               &                                       & \multicolumn{1}{l}{}                                                                             \\ \cline{2-5} 
		& \multicolumn{2}{l|}{$\min_{\X,\E}\ \norm{\X}_*+\lambda_1\norm{\X}_1+ \lambda_2 l(\E)$}                                                                                                                                                                                      & \multirow{2}{*}{\mcode{lrsr}}       & \multicolumn{1}{l}{\multirow{2}{*}{Low-rank and sparse representation  }}    \\
		& \multicolumn{2}{l|}{$\st \ \A=\B\X+\E$}                                                                                                                                                                                                                                     &                                       & \multicolumn{1}{l}{}                                                                             \\ \cline{2-5}
		& \multicolumn{2}{l|}{$\min_{\bL_i,\bS_i} \ \norm{\bL}_*+\lambda\sum_{i=1}^m \norm{\bS_i}_1$,}                                                                                                                                                                                & \multirow{2}{*}{\mcode{rmsc}}       & \multicolumn{1}{l}{\multirow{2}{*}{Robust multi-view spectral clustering  }} \\
		& \multicolumn{2}{l|}{$\st \ \X_i=\bL+\bS_i$, $i=1,\cdots,m$, $\bL\geq 0$, $\bL\bm{1}=\bm{1}$}                                                                                                                                                                                &                                       & \multicolumn{1}{l}{}                                                                             \\ \cline{2-5} 
		& \multicolumn{2}{l|}{$\min_{\Z_i,\E_i} \ \sum_{i=1}^K (\norm{\Z_i}_*+\lambda l(\E_i))+\alpha\norm{\Z}_{2,1}$}                                                                                                                                                                & \multirow{2}{*}{\mcode{mlap}}       & \multicolumn{1}{l}{\multirow{2}{*}{Multi-task low-rank affinity pursuit  }} \\
		& \multicolumn{2}{l|}{$\st \ \X_i=\X_i\Z_i+\E_i$, $i=1,\cdots,K$}                                                                                                                                                                                                             &                                       & \multicolumn{1}{l}{}                                                                             \\ \cline{2-5} 
		& \multicolumn{2}{l|}{$\min_{\bL,\bS} \ \norm{\bL}_*+\lambda\norm{\C \circ\bS}_1, \ \st \ \A=\bL+\bS, 0\leq \bL\leq 1$}                                                                                                                                                      & \mcode{igc}                                   & \multicolumn{1}{l}{Improved graph clustering  }                           \\ \cline{2-5}
		\multirow{7}{*}{}  & \multicolumn{2}{l|}{\multirow{1}{*}{$\min_{\bP} \ \langle\bP,\bL\rangle + \lambda \norm{\bP}_1, \ \st \ 0\preceq \bP \preceq \I, \text{Tr}(\bP)=k$}}                                                                                                                         & \multirow{1}{*}{\mcode{sparsesc}}   & \multicolumn{1}{l }{\multirow{1}{*}{Sparse spectral clustering  }}             \\ \hline\hline
		\multirow{12}{*}{\begin{tabular}[c]{@{}l@{}}Low-rank\\ \\ tensor\\ \\ models\end{tabular}} & \multicolumn{2}{l|}{\multirow{2}{*}{$\min_{\tL,\tS} \ \sum_{i=1}^k\alpha_i\norm{\tL_{i(i)}}_*+ \norm{\tS}_1, \ \st \ \tX=\tL+\tS$}}                                                                                                                                        & \multirow{2}{*}{\mcode{trpca_snn}} & Tensor robust PCA based on                                                                        \\
		& \multicolumn{2}{l|}{}                                                                                                                                                                                                                                                      &                                       & sum of nuclear norm                                                      \\ \cline{2-5}
		& \multicolumn{2}{l|}{\multirow{2}{*}{$\min_{\tX} \ \sum_{i=1}^k\alpha_i\norm{\tX_{i(i)}}_*, \ \st \ \Pomega(\tX)=\Pomega(\tM)$}}                                                                                                                                            & \multirow{2}{*}{\mcode{lrtc_snn}}  & Low-rank tensor completion based on        \\
		& \multicolumn{2}{l|}{}                                                                                                                                                                                                                                                      &                                       & sum of nuclear norm                                             \\ \cline{2-5}
		& \multicolumn{2}{l|}{$\min_{\tX,\tE} \ \sum_{i=1}^k\alpha_i\norm{\tX_{i(i)}}_*+\lambda l(\tE)$}                                                                                                                                                                              & \multirow{2}{*}{\mcode{lrtcR_snn}} & Reg. low-tank tensor completion based on                                                                  \\
		& \multicolumn{2}{l|}{$\st \ \Pomega(\tX)+\tE=\tM$}                                                                                                                                                                                                                           &                                       &  sum of nuclear norm                                                                      \\ \cline{2-5}
		& \multicolumn{2}{l|}{\multirow{2}{*}{$\min_{\tL,\tS} \ \norm{\tL}_*+\lambda \norm{\tS}_1, \ \st \ \tX=\tL+\tS$}}                                                                                                                                                            & \multirow{2}{*}{\mcode{trpca_tnn}} & Tensor Robust PCA based on                                                                        \\
		& \multicolumn{2}{l|}{}                                                                                                                                                                                                                                                      &                                       & tensor nuclear norm                                                       \\ \cline{2-5}
		& \multicolumn{2}{l|}{\multirow{2}{*}{$\min_{\tX} \ \norm{\tX}_*, \ \st \ \Pomega(\tX)=\Pomega(\tM)$}}                                                                                                                                                                       & \multirow{2}{*}{\mcode{lrtc\_tnn}}  & Low-rank tensor completion based on                                                                       \\
		& \multicolumn{2}{l|}{}                                                                                                                                                                                                                                                      &                                       & tensor nuclear norm                                                    \\ \cline{2-5}
		& \multicolumn{2}{l|}{\multirow{2}{*}{$\min_{\tX,\tE} \ \norm{\tX}_*+\lambda l(\tE), \ \st \ \Pomega(\tX)+\tE=\tM$}}                                                                                                                                                                       & \multirow{2}{*}{\mcode{lrtcR\_tnn}}  & Reg. low-rank tensor completion based on                                                                       \\
		& \multicolumn{2}{l|}{}                                                                                                                                                                                                                                                      &                                       & tensor nuclear norm                                                     \\ \hline
		
	\end{tabular}\\ 
       \leftline{ \footnotesize{$^*$In this table,  the loss function $l(\cdot)$ can be $\norm{\cdot}_1$, $\frac{1}{2}\norm{\cdot}_F^2$ and  $\norm{\cdot}_{2,1}$. The  $\norm{\cdot}_{2,1}$ norm is only applicable to the matrix. } }
\end{table}



	\end{document}


\title{Supplementary Material of Unified Alternating Direction Method of Multipliers by Majorization Minimization}
\maketitle
This document mainly contains two parts. In the first section, we give the implementation details of the experiments and some more explanations of the results. In the second section, we present an experiment to show the effectiveness of our proposed partition strategy for M-ADMM. In the third section, we give some mathematical background knowledge and useful results which are useful for the convergence analysis of M-ADMM. Finally, we give the convergence analysis of the Mixed Gauss-Seidel and Jacobian ADMM (M-ADMM) in the last section. 


\section{Proofs}

 \noindent
 \textbf{Proof of Lemma \ref{lem33}.} Let $h=f-f'$. By using (\ref{lipmuv}) and the convexity of $f$ and $f'$, for any $\x=[\x_1;\cdots;\x_n]$ and $\y=[\y_1;\cdots;\y_n]$ with $\x_i, \y_i\in\mathbb{R}^{p_i}$, $i=1,\cdots,n$.  we have
 \begin{eqnarray*}
 	0 \leq  f(\x)-f(\y)-\langle\nabla f(\y),\x-\y\rangle  \leq\frac{1}{2}\sum_{i=1}^{n}\|\x_i-\y_i\|_{\bL_i}^2, \\
 	-\frac{1}{2}\sum_{i=1}^{n}\|\x_i-\y_i\|_{\bL'_i}^2\leq -f'(\x)+f'(\y)+\langle\nabla f'(\y),\x-\y\rangle\leq 0.
 \end{eqnarray*}
 Summing the above two inequalities we have
 \begin{equation*}
 |h(\x)-h(\y)-\langle\nabla h(\y),\x-\y\rangle|\leq\frac{1}{2}\sum_{i=1}^{n}\|\x_i-\y_i\|_{\max\{\bL_i,\bL'_i\}}^2.
 \end{equation*}
 Thus $h$ is  $\{\max\{\bL_i,\bL'_i\}\}_{i=1}^n$-smooth. It is easy to see that $f+f'$ is $\{\bL_i+\bL'_i\}_{i=1}^n$-smooth  by applying (\ref{lipmuv}) for $f$ and $f'$.
 $\hfill\blacksquare$

 ~\
 
 \noindent
 \textbf{Proof of Proposition \ref{pro1madmm}.} First,  for $i\in B_1$, by the optimality of $\xkk_i$ to   problem (\ref{updatexbone2}) in Algorithm \ref{alg4},  there exists $\u_i^{k+1}\in\partial \hat{f}_i^k(\xkk_i)$ such that
 \begin{align}
 &-\u^{k+1}_i	=\nabla \hatr_i^k(\xkk_i)\notag\\
 \overset{\text{\ding{172}}}{=}&\A_i^\top\left(\betak\A_i(\xkk_i-\xk_i)+\betak(\A\xk-\b)+\lambdak\right)\notag\\
 &+\betak\Gk_i(\xkk_i-\xk_i)\notag\\
 =&\A_i^\top\left(\betak(\Abone\xbonekk+\Abtwo\xbtwok-\b)+\lambdak\right)\notag\\
 &-\betak\A_i^\top\Abone(\xbonekk-\xbonek) \notag \\
 & +\betak(\A_i^\top\A_i+\Gk_i)(\xkk_i-\xk_i)\notag\\
 \overset{\text{\ding{173}}}{=}&\A_i^\top\hatlambdakk-\betak\A_i^\top\Abone(\xbonekk-\xbonek)\notag\\
 &+\betak(\A_i^\top\A_i+\Gk_i)(\xkk_i-\xk_i)\notag,
 \end{align}
 where \ding{172}  uses the definition of $\hat{r}^k_i$ in (\ref{rone}), and   \ding{173} uses the definition of $\hatlambdakk$ in  (\ref{hatlambdatb}). 	
 A dot-product with $\xkk_i-\x_i$ on both sides of the above equation gives
 \begin{align}
 &-\langle\u_{B_1}^{k+1},\xbonekk-\xbone\rangle	=-\sumione\langle\u^{k+1}_i,\xkk_i-\x_i\rangle\notag\\
 =&\sumione\left\langle\A_i^\top\hatlambdakk-\betak\A_i^\top\Abone(\xbonekk-\xbonek),\xkk_i-\x_i \right\rangle\notag\\
 &+\sumione\left\langle \betak(\A_i^\top\A_i+\Gk_i)(\xkk_i-\xk_i), \xkk_i-\x_i\right\rangle\notag\\
 =&\langle\Abone^\top\hatlambdakk,\xbonekk-\xbone\rangle\notag\\
 &+\betak\langle\xbonekk-\xbonek,\xbonekk-\xbone\rangle_{\Kk_1-\Abone^\top\Abone}\notag \\
 =&\langle\Abone^\top\hatlambdakk,\xbonekk-\xbone\rangle+\frac{\betak}{2}\|\xbonekk-\xbonek\|^2_{\Kk_1-\Abone^\top\Abone} \notag\\
 &+\frac{\betak}{2}\|\xbone-\xbonekk\|^2_{\Kk_1-\Abone^\top\Abone} \notag\\
 &-\frac{\betak}{2}\|\xbone-\xbonek\|^2_{\Kk_1-\Abone^\top\Abone}\notag\\
 \overset{\text{\ding{172}}}{\geq}&\langle\Abone^\top\hatlambdakk,\xbonekk-\xbone\rangle+	\frac{\betak}{2}\|\xbone-\xbonek\|^2_{\Kk_1-\Abone^\top\Abone}\notag\\
 &-\frac{\betak}{2}\|\xbone-\xbonekk\|^2_{\Kk_1-\Abone^\top\Abone},\label{prop5main1}
 \end{align}
 where $\Kk_1=\text{Diag}\{\A_i^\top\A_i+\Gk_i,i\in B_1\}$  and \ding{172} uses $\|\xbonekk-\xbonek\|^2_{\Kk_1-\Abone^\top\Abone}\geq0$ due to (\ref{eqbteq2}).  
 
 Second, for $i\in B_2$, by the optimality of $\xkk_i$ to problem (\ref{updatexbtwo2}) in Algorithm \ref{alg4}, there exists $\u_i^{k+1}\in\partial \hat{f}_i^k(\xkk_i)$ such that
 \begin{align*}
 &-\u^{k+1}_i			=\nabla \hatr_i^k(\xkk_i)\\
 =&\A_i^\top\left(\betak\A_i(\xkk_i-\xk_i)+\betak(\Abone\xbonekk+\Abtwo\xbtwok-\b)\right)\\
 &+\A_i^\top\lambdak+\betak\Gk_i(\xkk_i-\xk_i)\\
 =&\A_i^\top\hatlambdakk+\betak(\A_i^\top\A_i+\Gk_i)(\xkk_i-\xk_i),
 \end{align*}
 where we use the definitions of $\hat{r}^k_i$ in (\ref{rtwo}) and $\hatlambdakk$ in (\ref{hatlambdatb}).
 A dot-product with $\xkk_i-\x_i$ on both sides of the above equation gives
 \begin{align}
 &-\langle\u_{B_2}^{k+1},\xbtwokk-\xbtwo\rangle	=-\sumitwo\langle\u^{k+1}_i,\xkk_i-\x_i\notag\rangle\\
 =&\sumitwo\left\langle\A_i^\top\hatlambdakk+\betak(\A_i^\top\A_i+\Gk_i)(\xkk_i-\xk_i),\xkk_i-\x_i \right\rangle\notag\\
 =&\langle\Abtwo^\top\hatlambdakk,\xbtwokk-\xbtwo\rangle+\betak\langle\xbtwokk-\xbtwok,\xbonekk-\xbone\rangle_{\Kk_2} \notag\\
 =&\langle\Abtwo^\top\hatlambdakk,\xbtwokk-\xbtwo\rangle-\frac{\betak}{2} \|\xbtwo-\xbtwok\|^2_{\Kk_2} \notag \\
 &+\frac{\betak}{2} \|\xbtwo-\xbtwokk\|^2_{\Kk_2}+\frac{\betak}{2}\|\xbtwokk-\xbtwok\|^2_{\Kk_2},\label{prop5main2}
 \end{align}
 where $\Kk_2=\text{Diag}\{\A_i^\top\A_i+\Gk_i,i\in B_2\}$. 
 
 Third, note that $\hatf^k\in\mathcal{S}_{\{\bL_i,\bP_i\}_{i=1}^n}(f,\xk)$, by using (\ref{keylm22}), we have
 \begin{align}
 & f(\xkk)-f(\x)\notag\\
 \leq&\langle \u^{k+1},\xkk-\x\rangle+\frac{1}{2}\sumi\left(\|\x_i-\xk_i\|^2_{\bL_i}-\|\x_i-\xkk_i\|^2_{\bP_i}\right)\notag\\
 \leq&\langle \u^{k+1},\xkk-\x\rangle+\frac{1}{2}\sumi\left(\|\x_i-\xk_i\|^2_{\bL_i}-\|\x_i-\xkk_i\|^2_{\bL_i}\right)\notag\\
 =&\langle \u_{B_1}^{k+1},\xbonekk-\xbone\rangle+\langle \u_{B_2}^{k+1},\xbtwokk-\xbtwo\rangle\notag \\
 &+\frac{1}{2}\sum_{j=1}^2\left(\|\x_{B_j}-\x^k_{B_j}\|_{\bL_{B_j}}^2-\|\x_{B_j}-\x^{k+1}_{B_j}\|^2_{\bL_{B_j}}\right)\notag\\
 \overset{\text{\ding{172}}}{\leq}&-\langle\A^\top\hatlambdakk,\xkk-\x\rangle-\frac{\betak}{2}\|\xbtwokk-\xbtwok\|^2_{\K_2}\notag\\
 &+\frac{\betak}{2}\sum_{j=1}^2\left(\|\x_{B_j}-\xk_{B_j}\|_{\Hk_j}^2-\|\x_{B_j}-\xkk_{B_j}\|_{\Hkk_j}^2\right) \notag
 \end{align}
 where   $\bL_{B_j}=\text{Diag}\left\{\bL_i, i\in B_j\right\}$  and \ding{172} uses (\ref{prop5main1})-(\ref{prop5main2}), the definitions of $\Hk_j$ in Proposition \ref{pro1madmm} and the fact $\betakk\geq\betak$.
 The proof is completed.  $\hfill\blacksquare$

 ~\
 
 \noindent
 \textbf{Proof of Proposition \ref{pro2madmm}.} By using line 4 of Algorithm \ref{alg3}, (\ref{Lem_Pythagoras2}) and the fact that $\betakk\geq\betak$, we have
 \begin{align}
 &\langle\A\xkk-\b,\blambda-\hatlambdakk\rangle		=\frac{1}{\betak}\langle\lambdakk-\lambdak,\blambda-\hatlambdakk\rangle\notag\\
 =&\frac{1}{2\betak}\left(\|\blambda-\lambdak\|^2-\|\blambda-\lambdakk\|^2\right)\notag\\
 &-\frac{1}{2\betak}\left(\|\hatlambdakk-\lambdak\|^2-\|\lambdakk-\hatlambdakk\|^2\right)\notag\\
 \leq&\frac{\betak}{2}\left(\|\blambda-\lambdak\|^2_{\Hk_3}-\|\blambda-\lambdakk\|^2_{\Hkk_3}\right)\notag\\
 &-\frac{1}{2\betak}\left(\|\hatlambdakk-\lambdak\|^2-\|\lambdakk-\hatlambdakk\|^2\right)\label{lm71}.
 \end{align}
 Consider the last two terms in (\ref{lm71}). We deduce
 \begin{align}
 &\frac{1}{2\betak}\left(\|\hatlambdakk-\lambdak\|^2-\|\lambdakk-\hatlambdakk\|^2\right)\notag\\
 \overset{\text{\ding{172}}}{=}&\frac{\betak}{2}\|\Abone\xbonekk+\Abtwo\xbtwok-\b\|^2-\|\Abtwo(\xbtwokk-\xbtwok)\|^2\notag\\
 =&\frac{\betak}{2}\left(\|\Abone\xbonekk+\Abtwo\xbtwok-\b\|^2-\|\xbtwokk-\xbtwok\|_{\Kk_2}^2\right)\notag\\
 &+\frac{\betak}{2}\|\xbtwokk-\xbtwok\|_{\Kk_2-\Abtwo^\top\Abtwo}^2	\notag\\
 \overset{\text{\ding{173}}}{\geq}&\frac{\betak}{2}\left(\|\Abone\xbonekk+\Abtwo\xbtwok-\b\|^2- \|\xbtwokk-\xbtwok\|_{\Kk_2}^2\right)\notag\\
 &+\frac{\tau}{\|\Abtwo\|^2_2}\|\Abtwo(\xbtwokk-\xbtwok)\|^2 \notag\\
 \overset{\text{\ding{174}}}{\geq}&{\betak\alpha}\|\Abone\xbonekk+\Abtwo\xbtwok-\b\|^2-\frac{\betak}{2}\|\xbtwokk-\xbtwok\|_{\Kk_2}^2\notag\\
 &+\betak\alpha\|\Abtwo(\xbtwokk-\xbtwok)\|^2\notag\\
 \geq&\frac{\betak\alpha}{2}\|\Abone\xbonekk+\Abtwo\xbtwokk-\b\|^2-\frac{\betak}{2}\|\xbtwokk-\xbtwok\|_{\Kk_2}^2\notag\\
 \overset{\text{\ding{175}}}{\geq}&\frac{\beta^{(0)}\alpha}{2}\|\Axkk-\b\|^2-\frac{\betak}{2}\|\xbtwokk-\xbtwok\|_{\Kk_2}^2\label{lm72},
 \end{align}
 where \ding{172} uses  (\ref{updatelambdaall}) and (\ref{hatlambdatb}), \ding{173} uses (\ref{eqbteq3}), \ding{174}  uses $\alpha=\min\left\{\frac{1}{2},\frac{\tau}{2\norm{\Abtwo}_2^2}\right\}$,
 and \ding{175} uses $\betak\geq\beta^{(k-1)}\geq\cdots\geq \beta^{(0)}$. 		
 The proof is completed by substituting (\ref{lm72}) into (\ref{lm71}). 
 $\hfill\blacksquare$
 
 ~\

\section{Implementation Details}
\label{secresults}
In this section, we give the implementation details of the experiments and some more explanations of the results.

\subsection{Latent Low-Rank Representation}
Consider the following Latent Low-Rank Representation (LRR) problem
\begin{equation}\label{latlrr0}
\min_{\Z,\bL} \|\Z\|_*+\|\bL\|_*+\frac{\lambda}{2}\|\A\Z+\bl\B-\X\|_F^2, \ \text{s.t.} \ \bm{1}^T\Z=\bm{1}^T,
\end{equation}
where $\lambda>0$ and $\|\cdot\|_F$ denotes the Frobenius norm of a matrix.
Problem (\ref{latlrr0}) is equivalent to
\begin{align}
\min_{\Z,\bL,\E} \|\Z\|_*+\|\bL\|_*+\frac{\lambda}{2}\|\E\|_F^2, 
\ \text{s.t.} \ \bm{1}^T\Z=\bm{1}^T,\ \A\Z+\bL\B-\X=\E.\label{latlrr220}
\end{align}

We show the details of M-ADMM to solve (\ref{latlrr0}) with 2 blocks and Linearized ADMM with Parallel Splitting (L-ADMM-PS) \cite{LADMPS} and M-ADMM to solve (\ref{latlrr220}) with 3 blocks.

\begin{enumerate}[(a)]
	\item Solve (\ref{latlrr0}) by M-ADMM

	The augmented Lagrangian function of (\ref{latlrr0}) is
	\begin{equation*}
	\mathcal{L}(\Z,\bL,\blambda,\beta)=\|\Z\|_*+\|\bL\|_*+\frac{\lambda}{2}\|\A\Z+\bL\B-\X\|_F^2+\langle\blambda,\bm{1}^T\Z-\bm{1}^T\rangle+\frac{\beta}{2}\|\bm{1}^T\Z-\bm{1}^T\|^2.
	\end{equation*}
	It is easy to verify that $\frac{1}{2}\|\A\Z+\bL\B-\X\|_F^2$ is $\{L_1\I,L_2\I\}$-smooth, where $L_1=2\norm{\A}_2^2$, $=L_2=2\norm{\B}_2^2$, and $\frac{1}{2}\|\bm{1}^T\Z-\bm{1}^T\|^2$  is $\eta$-smooth, where $\eta>\|\bm{1}\|^2$. By using these properties, M-ADMM solves (\ref{latlrr}) by the following updating rules
	
	\begin{equation*}
	\left\{
	\begin{aligned}
	\Z^{k+1}=&\arg\min_{\Z} \|\Z\|_*+\frac{\lambda L_1+\betak\eta}{2}\left\|\Z-\Z^k+\frac{\lambda\A^T(\A\Z^k+\bL^k\B-\X)+\bm{1}\left(\betak(\bm{1}^T\Z^k-\bm{1}^T)+\lambdak\right)}{\lambda L_1+\betak\eta} \right\|^2_F,  \\
	\bL^{k+1}=&\arg\min_{\bL}\|\bL\|_*+\frac{\lambda L_2}{2}\left\|\bL-\bL^k+\frac{(\A\Z^k+\bL^k\B-\X)\B^T}{L_2}\right\|_F^2,\\
	\lambdakk=&\lambdak+\betak(\bm{1}^T\Z^{k+1}-\bm{1}^T).
	\end{aligned}
	\right.
	\end{equation*}
	Note the $\Z$- and $\bL$-updates can be obtained by the Singular Value Thresholding (SVT) \cite{cai2010singular}.
	\item Solve (\ref{latlrr220}) by L-ADMM-PS
	
	The augmented Lagrangian function of (\ref{latlrr22}) is
	\begin{align*}
	\mathcal{L}(\Z,\bL,\blambda,\beta)=&\|\Z\|_*+\|\bL\|_*+\frac{\lambda}{2}\|\E\|_F^2+\langle\blambda_1,\bm{1}^T\Z-\bm{1}^T\rangle+\frac{\beta}{2}\|\bm{1}^T\Z-\bm{1}^T\|^2\\
	&+\langle\blambda_2,\A\Z+\bL\B-\X-\E\rangle+\frac{\beta}{2}\|\A\Z+\bL\B-\X-\E\|_F^2.
	\end{align*}
	
	Note that $h(\Z,\bL,\E)=\frac{1}{2}\|\bm{1}^T\Z-\bm{1}^T\|^2+\frac{1}{2}\|\A\Z+\bL\B-\X-\E\|_F^2$ is $\{\eta_1\I,\eta_2\I,\eta_3\I\}$-smooth, where $\eta_1>\|\bm{1}\|^2+3\|\A\|_2^2$, $\eta_2>3\|\B\|_2^2$ and $\eta_3>3$.
	By using such a property, 
	L-ADMM-PS solves (\ref{latlrr22}) by the following updating rules 
	
	\begin{equation*}
	\left\{
	\begin{aligned}
	\Z^{k+1}=&\arg\min_{\Z} \|\Z\|_*+\frac{\betak\eta_1}{2}\left\|\Z-\Z^k+\frac{\bm{1}(\blambda^k_1+\betak(\bm{1}^T\Z^k-\bm{1}^T))+\A^T(\blambda_2^k+\betak(\A\Z^k+\bL^k\B-\X-\E^k))}{\betak\eta_1} \right\|_F^2,  \\
	\bL^{k+1}=&\arg\min_{\bL} \|\bL\|_*+\frac{\betak\eta_2}{2}\left\|\bL-\bL^k+\frac{(\blambda^k_2+\betak(\A\Z^k+\bL^k\B-\X-\E^k))\B^T}{\betak\eta_2} \right\|_F^2,\\
	\E^{k+1}=&\arg\min_{\E}\frac{\lambda}{2}\|\E\|_F^2-\langle\blambda^k_2+\betak(\A\Z^{k}+\bL^k\B-\X-\E^k),\E\rangle+\frac{\betak\eta_3}{2}\|\E-\E^k\|_F^2,\\
	\lambdakk_1=&\lambdak_1+\betak(\bm{1}^T\Z^{k+1}-\bm{1}^T),\\
	\lambdakk_2=&\lambdak_2+\betak(\A\Z^{k+1}+\bL^{k+1}\B-\X-\E^{k+1}),
	\end{aligned}
	\right.
	\end{equation*}

	\item Solve (\ref{latlrr220}) by M-ADMM
	
	M-ADMM divides the variables $\{\Z,\bL,\E\}$ into two super blocks, i.e., $\{\Z\}$ and $\{\bL,\E\}$. Then it solves (\ref{latlrr220}) by the following updating rules
	\begin{equation*}
	\left\{
	\begin{aligned}
	\Z^{k+1}=&\arg\min_{\Z} \|\Z\|_*+\frac{\betak\eta_1}{2}\left\|\Z-\Z^k+\frac{\bm{1}(\blambda^k_1+\betak(\bm{1}^T\Z^k-\bm{1}^T))+\A^T(\blambda_2^k+\betak(\A\Z^k+\bL^k\B-\X-\E^k))}{\betak\eta_1} \right\|_F^2,  \\
	\bL^{k+1}=&\arg\min_{\bL} \|\bL\|_*+\frac{\betak\eta_2}{2}\left\|\bL-\bL^k+\frac{(\blambda^k_2+\betak(\A\Z^{k+1}+\bL^k\B-\X-\E^k))\B^T}{\betak\eta_2} \right\|_F^2,\\
	\E^{k+1}=&\arg\min_{\E}\frac{\lambda}{2}\|\E\|_F^2+\frac{\betak}{2}\left\|\A\Z^{k+1}+\bL^k\B-\X-\E+\frac{\lambdak_2}{\betak}\right\|^2_F+\frac{\betak\eta_3}{2}\|\E-\E^k\|_F^2,\\
	\lambdakk_1=&\lambdak_1+\betak(\bm{1}^T\Z^{k+1}-\bm{1}^T),\\
	\lambdakk_2=&\lambdak_2+\betak(\A\Z^{k+1}+\bL^{k+1}\B-\X-\E^{k+1}),
	\end{aligned}
	\right.
	\end{equation*}
	where $\eta_1=\|\bm{1}\|^2+\|\A\|_2^2$, $\eta_2>2\|\B\|_2^2$ and $\eta_3>1$.

\end{enumerate}

\subsection{Latent Low-Rank Representation}
 Consider the following Latent Low-Rank Representation (LRR) problem
 \begin{equation}\label{latlrr}
 \min_{\Z,\bL} \|\Z\|_*+\|\bL\|_*+\frac{\lambda}{2}\|\X\Z+\bL\X-\X\|_F^2, \ \text{s.t.} \ \bm{1}^T\Z=\bm{1}^T,
 \end{equation}
 where $\lambda>0$ and $\|\cdot\|_F$ denotes the Frobenius norm of a matrix.
 Problem (\ref{latlrr}) is equivalent to
 \begin{align}
 \min_{\Z,\bL,\E} \|\Z\|_*+\|\bL\|_*+\frac{\lambda}{2}\|\E\|_F^2, 
 \ \text{s.t.} \ \bm{1}^T\Z=\bm{1}^T,\ \X\Z+\bL\X-\X=\E.\label{latlrr22}
 \end{align}
 
 We show the details of M-ADMM to solve (\ref{latlrr}) with 2 blocks and Linearized ADMM with Parallel Splitting (L-ADMM-PS) \cite{LADMPS} and M-ADMM to solve (\ref{latlrr22}) with 3 blocks.
 
 \begin{enumerate}[(a)]
 	\item Solve (\ref{latlrr}) by M-ADMM

 The augmented Lagrangian function of (\ref{latlrr}) is
 \begin{equation*}
 \mathcal{L}(\Z,\bL,\blambda,\beta)=\|\Z\|_*+\|\bL\|_*+\frac{\lambda}{2}\|\X\Z+\bL\X-\X\|_F^2+\langle\blambda,\bm{1}^T\Z-\bm{1}^T\rangle+\frac{\beta}{2}\|\bm{1}^T\Z-\bm{1}^T\|^2.
 \end{equation*}
It is easy to verify that $\frac{1}{2}\|\X\Z+\bL\X-\X\|_F^2$ is $\{L_1\I,L_2\I\}$-smooth, where $L_1=L_2=2\|\X\|_2^2$, and $\frac{1}{2}\|\bm{1}^T\Z-\bm{1}^T\|^2$  is $\eta$-smooth, where $\eta>\|\bm{1}\|^2$. By using these properties, M-ADMM solves (\ref{latlrr}) by the following updating rules

\begin{equation*}
\left\{
\begin{aligned}
\Z^{k+1}=&\arg\min_{\Z} \|\Z\|_*+\frac{\lambda L_1+\betak\eta}{2}\left\|\Z-\Z^k+\frac{\lambda\X^T(\X\Z^k+\bL^k\X-\X)+\bm{1}\left(\betak(\bm{1}^T\Z^k-\bm{1}^T)+\lambdak\right)}{\lambda L_1+\betak\eta} \right\|^2_F,  \\
\bL^{k+1}=&\arg\min_{\bL}\|\bL\|_*+\frac{\lambda L_2}{2}\left\|\bL-\bL^k+\frac{(\X\Z^k+\bL^k\X-\X)\X^T}{L_2}\right\|_F^2,\\
\lambdakk=&\lambdak+\betak(\bm{1}^T\Z^{k+1}-\bm{1}^T).
\end{aligned}
\right.
\end{equation*}
Note the $\Z$- and $\bL$-updates can be obtained by the Singular Value Thresholding (SVT) \cite{cai2010singular}.

 
\item Solve (\ref{latlrr22}) by L-ADMM-PS

The augmented Lagrangian function of (\ref{latlrr22}) is
\begin{align*}
\mathcal{L}(\Z,\bL,\blambda,\beta)=&\|\Z\|_*+\|\bL\|_*+\frac{\lambda}{2}\|\E\|_F^2+\langle\blambda_1,\bm{1}^T\Z-\bm{1}^T\rangle+\frac{\beta}{2}\|\bm{1}^T\Z-\bm{1}^T\|^2\\
&+\langle\blambda_2,\X\Z+\bL\X-\X-\E\rangle+\frac{\beta}{2}\|\X\Z+\bL\X-\X-\E\|_F^2.
\end{align*}

Note that $h(\Z,\bL,\E)=\frac{1}{2}\|\bm{1}^T\Z-\bm{1}^T\|^2+\frac{1}{2}\|\X\Z+\bL\X-\X-\E\|_F^2$ is $\{\eta_1\I,\eta_2\I,\eta_3\I\}$-smooth, where $\eta_1>\|\bm{1}\|^2+3\|\X\|_2^2$, $\eta_2>3\|\X\|_2^2$ and $\eta_3>3$.
By using such a property, 
L-ADMM-PS solves (\ref{latlrr22}) by the following updating rules 
 
\begin{equation*}
\left\{
\begin{aligned}
\Z^{k+1}=&\arg\min_{\Z} \|\Z\|_*+\frac{\betak\eta_1}{2}\left\|\Z-\Z^k+\frac{\bm{1}(\blambda^k_1+\betak(\bm{1}^T\Z^k-\bm{1}^T))+\X^T(\blambda_2^k+\betak(\X\Z^k+\bL^k\X-\X-\E^k))}{\betak\eta_1} \right\|_F^2,  \\
\bL^{k+1}=&\arg\min_{\bL} \|\bL\|_*+\frac{\betak\eta_2}{2}\left\|\bL-\bL^k+\frac{(\blambda^k_2+\betak(\X\Z^k+\bL^k\X-\X-\E^k))\X^T}{\betak\eta_2} \right\|_F^2,\\
\E^{k+1}=&\arg\min_{\E}\frac{\lambda}{2}\|\E\|_F^2-\langle\blambda^k_2+\betak(\X\Z^{k}+\bL^k\X-\X-\E^k),\E\rangle+\frac{\betak\eta_3}{2}\|\E-\E^k\|_F^2,\\
\lambdakk_1=&\lambdak_1+\betak(\bm{1}^T\Z^{k+1}-\bm{1}^T),\\
\lambdakk_2=&\lambdak_2+\betak(\X\Z^{k+1}+\bL^{k+1}\X-\X-\E^{k+1}),
\end{aligned}
\right.
\end{equation*}

\item Solve (\ref{latlrr22}) by M-ADMM

M-ADMM divides the variables $\{\Z,\bL,\E\}$ into two super blocks, i.e., $\{\Z\}$ and $\{\bL,\E\}$. Then it solves (\ref{latlrr22}) by the following updating rules
\begin{equation*}
\left\{
\begin{aligned}
\Z^{k+1}=&\arg\min_{\Z} \|\Z\|_*+\frac{\betak\eta_1}{2}\left\|\Z-\Z^k+\frac{\bm{1}(\blambda^k_1+\betak(\bm{1}^T\Z^k-\bm{1}^T))+\X^T(\blambda_2^k+\betak(\X\Z^k+\bL^k\X-\X-\E^k))}{\betak\eta_1} \right\|_F^2,  \\
\bL^{k+1}=&\arg\min_{\bL} \|\bL\|_*+\frac{\betak\eta_2}{2}\left\|\bL-\bL^k+\frac{(\blambda^k_2+\betak(\X\Z^{k+1}+\bL^k\X-\X-\E^k))\X^T}{\betak\eta_2} \right\|_F^2,\\
\E^{k+1}=&\arg\min_{\E}\frac{\lambda}{2}\|\E\|_F^2+\frac{\betak}{2}\left\|\X\Z^{k+1}+\bL^k\X-\X-\E+\frac{\lambdak_2}{\betak}\right\|^2_F+\frac{\betak\eta_3}{2}\|\E-\E^k\|_F^2,\\
\lambdakk_1=&\lambdak_1+\betak(\bm{1}^T\Z^{k+1}-\bm{1}^T),\\
\lambdakk_2=&\lambdak_2+\betak(\X\Z^{k+1}+\bL^{k+1}\X-\X-\E^{k+1}),
\end{aligned}
\right.
\end{equation*}
where $\eta_1=\|\bm{1}\|^2+\|\X\|_2^2$, $\eta_2>2\|\X\|_2^2$ and $\eta_3>1$.

\item Solve (\ref{latlrr}) by  ADMM

Reformulate (\ref{latlrr}) as
\begin{equation}\label{latlrr33}
\min_{\J,\bS,\Z,\bL} \|\J\|_*+\|\bS\|_*+\frac{\lambda}{2}\|\X\Z+\bL\X-\X\|_F^2, \ \st \ \onev^T\Z=\onev^T, \Z=\J, \bL=\bS.
\end{equation}
\end{enumerate}
The augmented Lagrangian function is
\begin{align*}
\mathcal{L}(\J,\bS,\Z,\bL,\blambda,\beta) = & \|\J\|_*+\|\bS\|_*+\frac{\lambda}{2}\|\X\Z+\bL\X-\X\|_F^2+\langle\blambda_1,\bm{1}^T\Z-\bm{1}^T\rangle+\frac{\beta}{2}\|\bm{1}^T\Z-\bm{1}^T\|^2\\
&+\langle\blambda_2,\Z-\J\rangle+\frac{\beta}{2}\|\Z-\J\|^2+\langle\blambda_3,\bL-\bS\rangle+\frac{\beta}{2}\|\bL-\bS\|^2.
\end{align*}
Then   ADMM solves (\ref{latlrr33}) by
\begin{equation*}
\left\{
\begin{aligned}
\Z^{k+1} =& (\lambda\X^T\X+\betak\onev\onev^T+\betak\I)^{-1}(-\lambda\X^T(\bL^k\X-\X)-\onev\blambda_1^k+\betak\onev\onev^T-\blambda_2^k+\betak\J^{k}),\\
\bS^{k+1}=&\arg\min_{\bS} \|\bS\|_*+\frac{\betak}{2}\left\|\bS-\bL^k-\frac{\blambda^k_3}{\betak} \right\|_F^2,  \\
\J^{k+1}=&\arg\min_{\J} \|\J\|_*+\frac{\betak}{2}\left\|\J-\Z^{k+1}-\frac{\blambda^k_2}{\betak} \right\|_F^2,  \\
\bL^{k+1}=& (-\lambda(\X\Z^{k+1}-\X)\X^T-\blambda_3^k+\betak\bS^{k+1})(\lambda\X\X^T+\betak\I)^{-1},\\
\lambdakk_1=&\lambdak_1+\betak(\bm{1}^T\Z^{k+1}-\bm{1}^T),\\
\lambdakk_2=&\lambdak_2+\betak(\Z^{k+1}-\J^{k+1}),\\
\lambdakk_3=&\lambdak_3+\betak(\bL^{k+1}-\bS^{k+1}),\\
\betakk=&\min(\rho\betak,\mu_{\max}),
\end{aligned}
\right.
\end{equation*}
where $\rho=1.1$, $\mu_{\max}=10^6$, $\beta^{(0)}=10^{-4}$.

\subsection{Nonnegative Sparse Coding}
We consider the following nonnegative sparse coding problem
\begin{align}
\min_{\x,\e}\sumi\|\x_i\|_1+\lambda\|\e\|_1, \ 
\text{s.t. }  \y=\sumi\A_i\x_i+\e, \x_i\geq\bzero, \alli,\label{nonsparse2}
\end{align}
 We show the details of Linearized ADMM with Parallel Splitting (L-ADMM-PS) \cite{LADMPS} and M-ADMM to solve (\ref{nonsparse2}).
 
 \begin{enumerate}[(a)]
 	\item Solve (\ref{nonsparse2}) by L-ADMM-PS
 	\begin{equation*}
 	\left\{
 	\begin{aligned}
 	\x_i^{k+1}=&\arg\min_{\x_i\geq \bm{0}} \|\x_i\|_1+\frac{\betak\eta_i}{2}\left\|\x_i-\x_i^k+\frac{\A_i^T(\blambda^k+\betak(\A\xk+\e^k-\y))}{\betak\eta_i} \right\|^2, \ i=1,\cdots,n, \\
 	\e^{k+1}=&\arg\min_{\e} \lambda\|\e\|_1+\frac{\betak\eta_{\e}}{2}\left\|\e-\e^k+\frac{\blambda^k+\betak(\A\xk+\e^k-\y)}{\betak\eta_{\e}} \right\|^2,\\
 	\lambdakk=&\lambdak+\betak(\A\xkk+\e^{k+1}-\y),
 	\end{aligned}
 	\right.
 	\end{equation*}
 	where $\eta_i>(n+1)\|\A_i\|_2^2$ and $\eta_{\e}>n+1$. Note that $\x_i$- and $\e$-updates have closed form solutions.

 	\item Solve (\ref{nonsparse2}) by M-ADMM
 	
 	First, the variables $\x_1,\cdots,\x_n$ are divided into two super blocks, i.e., $\x_i$'s for $i\in B_1$ and $\x_i$'s for $i\in B_2$. Let $|B_1|=n_1$ and $|B_2|=n_2$. The variable $\e$ is further grouped into the second block. Then M-ADMM solves (\ref{nonsparse2}) by the following updating rules
 	\begin{equation*}
 	\left\{
 	\begin{aligned}
 	\x_i^{k+1}=&\arg\min_{\x_i\geq \bzero} \|\x_i\|_1+\frac{\betak\eta_i}{2}\left\|\x_i-\x_i^k+\frac{\A_i^T(\blambda^k+\betak(\A\xk+\e^k-\y))}{\betak\eta_i} \right\|^2, \ i\in B_1, \\
 	\x_i^{k+1}=&\arg\min_{\x_i\geq \bzero} \|\x_i\|_1+\frac{\betak\eta_i}{2}\left\|\x_i-\x_i^k+\frac{\A_i^T(\blambda^k+\betak(\A_{B_1}\xkk_{B_1}+A_{B_2}\xk_{B_2}+\e^k-\y))}{\betak\eta_i} \right\|^2, \ i\in B_2, \\
 	\e^{k+1}=&\arg\min_{\e} \lambda\|\e\|_1+\frac{\betak}{2}\left\|\A_{B_1}\xkk_{B_1}+\A_{B_2}\xk_{B_2}+\e-\y+\frac{\lambdak}{\betak}\right\|^2+\frac{\betak\eta_{\e}}{2}\|\e-\e^k\|^2,\\
 	\lambdakk=&\lambdak+\betak(\A\xkk+\e^{k+1}-\y),
 	\end{aligned}
 	\right.
 	\end{equation*}
	where $\eta_i>n_1\|\A_i\|_2^2$ for $i\in B_1$, $\eta_i>(n_2+1)\|\A_i\|_2^2$ for $i\in B_2$ and $\eta_{\e}>n_2$. 
\end{enumerate}

\bibliographystyle{aaai}
\bibliography{supp}